\documentclass[seceqn,MSNbibl,number,citesort,dvips]{arxbj}
\usepackage{graphicx}

\aid{0}
\volume{20}
\issue{3}
\pubyear{2014}
\firstpage{1292}
\lastpage{1343}
\doi{10.3150/13-BEJ522} 

\makeatletter

\newcommand{\rright}{\right}
\newcommand{\lleft}{\left}

\newtheorem{theorem}{Theorem}[section]
\newtheorem{corollary}{Corollary}[section]
\newremark{Example}{Example}
\newtheorem{lemma}{Lemma}[section]
\newtheorem{proposition}{Proposition}[section]
\newremark{Remark}{Remark}

\newcommand{\A}{{\mathcal{A}}}
\newcommand{\I}{\tilde{L}}
\newcommand{\p}{\bar{p}}
\newcommand{\q}{\bar{q}}

\newcommand{\Var}{\mathrm{Var}}
\newcommand{\eqref}[1]{(\ref{#1})}
\renewcommand{\epsilon}{\varepsilon}
\makeatother

\begin{document}
\begin{frontmatter}

\title{Optimal alignments of longest common subsequences and their path
properties}
\runtitle{Path
properties}

\begin{aug}
\author[1]{\inits{J.}\fnms{J\"uri} \snm{Lember}\corref{}\thanksref{1}\ead[label=e1]{jyril@ut.ee}},
\author[2]{\inits{H.}\fnms{Heinrich} \snm{Matzinger}\thanksref{2}\ead[label=e2]{matzing@math.gatech.edu}}
\and
\author[3]{\inits{A.}\fnms{Anna} \snm{Vollmer}\thanksref{3}\ead[label=e3]{anna-lisa.vollmer@plymouth.ac.uk}}
\runauthor{J. Lember, H. Matzinger and A. Vollmer} 
\address[1]{University of Tartu,
Liivi 2-513 50409, Tartu,
Estonia. \printead{e1}}
\address[2]{Georgia Tech,
School of Mathematics,
Atlanta, GA 30332-0160, USA.\\ \printead{e2}}
\address[3]{University of Plymouth, School of Computation and Mathematics,
Plymouth, PL4 8AA, Devon, UK. \printead{e3}}
\end{aug}

\received{\smonth{12} \syear{2011}}
\revised{\smonth{11} \syear{2012}}

%
\begin{abstract}
We investigate the behavior of optimal alignment paths
for homologous (related) and independent random sequences. An alignment between
two finite sequences is optimal if it corresponds to the longest
common subsequence (LCS). We prove the existence of lowest and
highest optimal alignments and study their differences. High
differences between the extremal alignments imply the high variety of
all optimal alignments. We present several simulations indicating
that the homologous (having the same common ancestor) sequences have
typically the distance between the extremal alignments
of much smaller size than independent sequences. In particular, the
simulations suggest that for the homologous sequences, the growth
of the distance between the extremal alignments is logarithmical.
The main theoretical results of the paper prove that (under some
assumptions) this is the case, indeed. The paper suggests that the
properties of the optimal alignment paths characterize the
relatedness of the sequences.
\end{abstract}

%

\end{frontmatter}
%
\section{Introduction}\label{intro}
Let $\mathcal{A}$ be a finite alphabet. In everything that follows,
$X=X_1\ldots X_n\in\mathcal{A}^n$ and $Y=Y_1\ldots Y_n\in\mathcal{A}^n$
are two strings of length $n$. A common subsequence of $X$ and $Y$
is a sequence that is a subsequence of $X$ and at the same time of
$Y$. We denote by $L_n$ the length of the \textit{longest common
subsequence (LCS)} of $X$ and $Y$. LCS is a special case of a
sequence alignment that is a very important tool in computational
biology, used for comparison of DNA and protein sequences (see, e.g.,
\cite{watermanintrocompbio,Vingron,Watermanphase,Durbin,china}).
They are also used in computational linguistics, speech recognition
and so on. In all these applications, two strings with a relatively
long LCS, are deemed related. Hence, to distinguish related pairs of
strings from unrelated via the length of LCS (or other similar
optimality measure), it is important to have some information about
the (asymptotical) distribution of $L_n$. Unfortunately, although
studied for a relatively long time, not much about the statistical
behavior of $L_n$ is known even when the sequences $X_1,X_2,\ldots$
and $Y_1,Y_2,\ldots$ are both i.i.d. and independent of each other.
Using the subadditivity, it is easy to see the existence of a
constant $\gamma$ such that
%
\begin{equation}
\label{sub} {L_n\over n}\to\gamma\qquad \mbox{a.s. and in }L_1.
\end{equation}
(see,
e.g., \cite{Alexander,Vingron,rate}). Referring to the celebrated
paper of Chvatal and Sankoff \cite{Sankoff1}, the constant $\gamma$
is called the \textit{Chvatal--Sankoff constant}; its value is unknown
for even as simple cases as i.i.d. Bernoulli sequences. In this
case, the value of $\gamma$ obviously depends on the Bernoulli
parameter~$p$. When $p=0.5$, the various bounds indicate that
$\gamma\approx0.81$ \cite{steele86,kiwi,Baeza1999}. For a smaller
$p$, $\gamma$ is even bigger. Hence, a common subsequence of two
independent Bernoulli sequences typically makes up large part of the
total length, if the sequences are related, LCS is even larger. As
for the mean of $L_n$, not much is also known about the variance of
$L_n$. In \cite{Sankoff1}, it was conjectured that for Bernoulli
parameter
$p=0.5$, the variance is of order $\mathrm{o}(n^{2/3})$. Using an
Efron--Stein type of
inequality, Steele \cite{steele86} proved $\Var [L_n]\leq
2p(1-p)n$. In \cite{Waterman-estimation}, Waterman conjectured that
$\Var [L_n]$ grows linearly. In series of
papers, Matzinger and others prove the Waterman conjecture for
different models
\cite{BonettoLCS,periodicLCS,LCIS,notsym}.

Because of relatively rare knowledge about its asymptotics, it is
rather difficult to build any statistical test based on $L_n$ or any
other global optimality criterion. The situation is better for local
alignments (see e.g., \cite{Watermanphase,Waterman-estimation}),
because for these alignments approximate $p$-values were recently
calculated
\cite{siegmund,hansen}.

In the present paper, we propose another approach -- instead of
studying the length of LCS, we investigate the properties and
behavior of the optimal alignments. Namely, even for moderate $n$,
the LCS is hardly unique. Every LCS corresponds to an optimal
alignment (not necessarily vice versa), so in general, we have
several optimal alignments. The differences can be of the local
nature meaning that the optimal alignments do not vary much, or they
can be of global nature. We conjecture that the variation of the
optimal alignments characterizes the relatedness or homology of the
sequences. To measure the differences between various optimal
alignment, we consider so-called extremal alignments and study their
differences.
\begin{Example*} Let us consider a practical example to give an insight
in what follows. Let \mbox{$X=\,$ATAGCGT}, $Y=\,$CAACATG. There are two longest
common subsequences: {AACG} and {AACT}. Thus, $L_7=4$. To every
longest common subsequence corresponds two optimal alignments. These
optimal alignments can be presented as follows:
\vspace*{7pt}
\begin{center}

\includegraphics{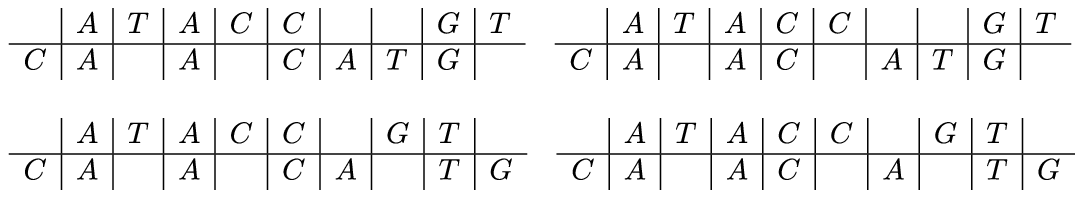}

\end{center}
\vspace*{2pt}
First, two alignments correspond to optimal subsequence AACG, the
last two correspond to AACT. In the following, we shall often
consider the alignments as the pairs
$\{(i_1,j_1),\ldots,(i_k,j_k)\}$, where $X_{i_t}=Y_{i_t}$ for every
$t=1,\ldots,k$. With this notation, the four optimal alignments
above are $\{(1,2), (3,3), (5,4), (6,7)\}$, $\{(1,2), (3,3), (4,4),
(6,7)\}$, $\{(1,2),(3,3),(5,4),(7,6)\}$ and
$\{(1,2),(3,3),(4,4),(7,6)\}$. We now represent every alignment as
two-dimensional plot. For the alignments in our example, the
two dimensional plots are
\begin{center}

\includegraphics{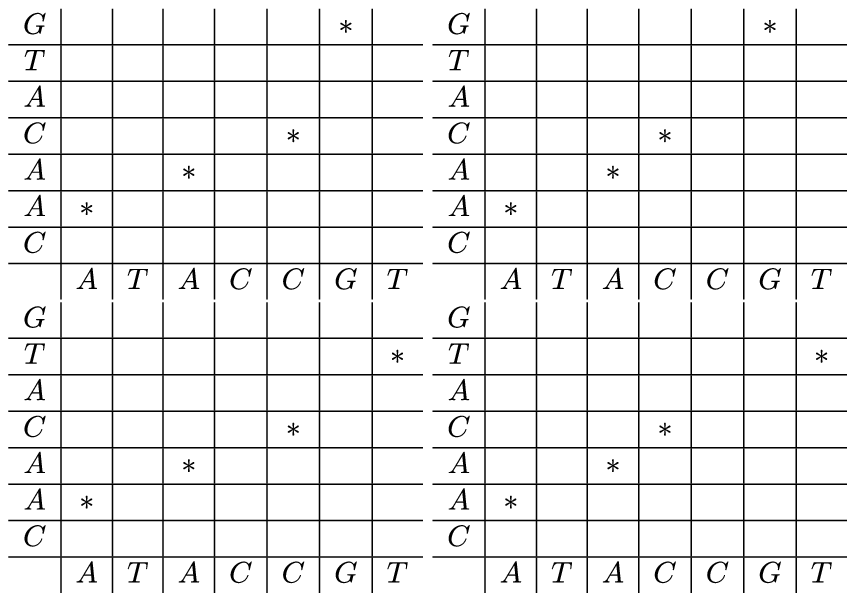}

\end{center}
\vspace*{3pt}
Putting all four alignment into one graph, we see that on some
regions all alignments are unique, but on some region, they vary:
\vspace*{7pt}
\begin{center}

\includegraphics{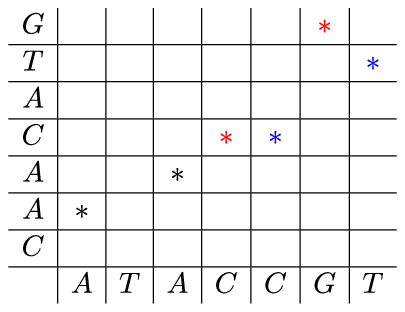}

\end{center}
\vspace*{3pt}
In the picture above, the alignment $(1, 2), (3, 3), (4, 4), (6, 7)$ (corresponding to AACG)
lies above all others. This alignment will be called \textit{highest alignment}. Similarly the
alignment $(1, 2), (3, 3), (5, 4), (7, 6)$ (corresponding to AACT) lies below all others. This
alignment will be called \textit{lowest alignment}. The highest and lowest alignment will be called
\textit{extremal alignments}.

Thus, the highest (lowest) alignment is the optimal alignment that
lies above (below) all other optimal alignments in two-dimensional
representation. For big $n$, we usually align the dots in the two
dimensional representation by lines. Then, to every alignment
corresponds a curve. We shall call this curve the \textit{alignment
graph} (when it is obvious from the context, we skip ``graph''). In
Figure~\ref{fig1}, there are extremal alignments  of two
independent i.i.d. four letter sequences (with uniform marginal
distributions) of length $n=1000$. It is visible that the extremal
alignments are rather far from each other, in particular, the
maximum vertical and horizontal distances are relatively big.
%
\begin{figure}

\includegraphics{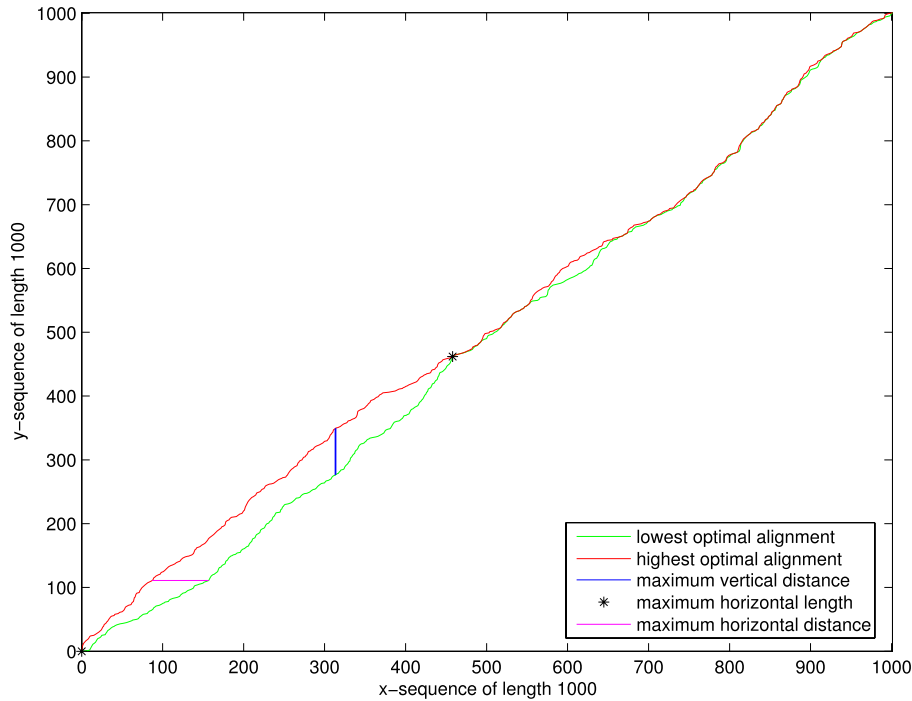}

\caption{The extremal alignments of two independent i.i.d. four letter
sequences.}\label{fig1}
\end{figure}
%

We call the sequences $X$ and $Y$ unrelated, if they are
independent. There are many ways to model the related sequences, the
model in the present paper is based on the assumption that there
exists a common ancestor, from which both sequences $X$ and $Y$ are
obtained by independent random mutations and deletions. The
sequences with common ancestor are called homologous, detecting
homology of given sequences is one of the major tasks in modern
computational molecular biology \cite{china}. In this paper, we
shall call the homologous sequences \textit{related}.

More precisely, we consider an $\mathcal{A}$-valued i.i.d. process
$Z_1,Z_2,\ldots$ that will be referred to as the common ancestor or
the ancestor process. A letter $Z_i$ has a probability to mutate
according to a transition matrix that does not depend on $i$. The
mutations of the letters are assumed to be independent. After
mutations, some letters of the mutated process disappear. The
disappearance is modeled via a deletion process $D^x_1,D^x_2,\ldots$
that is assumed to be an i.i.d. Bernoulli sequence with parameter
$p$, that is, $P(D^x_i=1)=p$. If $D^x_i=0$, then the $i$th (possibly
mutated letter) disappears. In such a way, a random sequence
$X_1,X_2,\ldots$ is obtained. The sequence $Y_1,Y_2,\ldots$ is
obtained similarly: the ancestor process $Z_1,Z_2,\ldots$ is the
same, but the mutations and deletions (with the same probabilities)
are independent of the ones used to generate $X$-sequence. The
formal definition is given in
Section~\ref{subsec:relateddef}.

Figure~\ref{fig2} presents a typical picture or extremal alignments of two
related four-letter sequences (of uniform marginal distribution) of
length 949. The sequences in Figure~\ref{fig2}, thus, have the same marginal
distribution as the ones in Figure~\ref{fig1}, but they are not independent
any more. Clearly the extremal alignments are close to each other;
in particular the maximal vertical and horizontal distance is much
smaller than these ones in Figure~\ref{fig1}.\looseness=1

\begin{figure}

\includegraphics{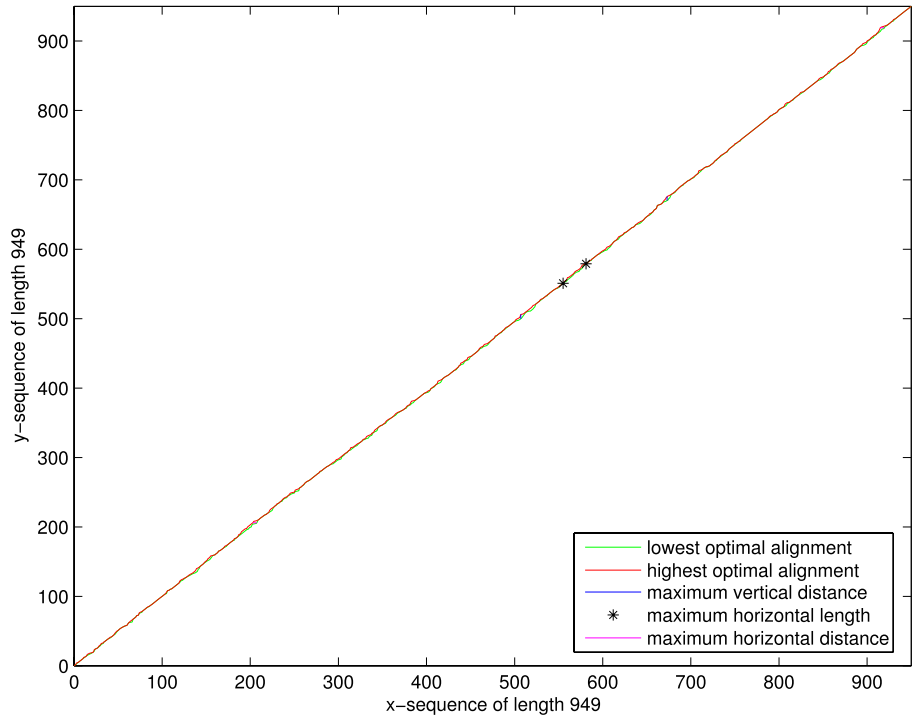}

\caption{The extremal alignments of two related four letter
sequences.}\label{fig2}
\end{figure}
%

Figures~\ref{fig1} and \ref{fig2} as well as many other similar simulations
(see \cite{VollmerReport}) clearly indicate that for related
sequences the differences of optimal alignments are of local nature,
whilst for independent sequences they vary much more. This motivates
us to find a way to quantify the non-uniqueness and use the
obtained characteristic as a measure of the relatedness. For that we
measure the differences of extremal alignments in several ways: the
maximal vertical and horizontal distance and Hausdorff's distance
(see Section~\ref{defHauss} for formal definition of Hausdorff's
distance). The simulations in Section~\ref{sec:simulations} show
that for independent sequences, the growth of both of them is almost
linear; for related sequence, however, it is logarithmic. Under some
assumptions, the latter is confirmed by the main theoretical
results about related sequences, Theorems \ref{main1}, \ref{main2}
and \ref{main3}. More specifically, Theorem~\ref{main1} states that
under some assumption that never holds for independent sequences,
there exist universal constants $C$ and $D$ so that for $n$ big
enough,\looseness=0
\[
P\bigl(h_o(H,L)>C\ln n\bigr)\leq Dn^{-2}.
\]
Here, $h_o(H,L)$ stands for a slight modification of Hausdorff's
distance between extremal alignments, which we shall call restricted
Hausdorff's distance. We conjecture the result also holds for (full)
Hausdorff's distance, denoted by $h$. Note that by Borel--Cantelli
lemma, from the inequality above, it follows that
\[
P\bigl(h_o(H,L)\leq C\ln n, \mbox{eventually }\bigr)=1,
\]
that is, the ratio $h_o(H,L)/\ln n$ is eventually bounded above by
$C$, a.s. Theorem~\ref{main2} states the similar result with maximal
vertical distance instead of Hausdorff's distance.

Theorem~\ref{main3} considers the sequences with random lengths. The
expected length of both sequences is $n$, the randomness comes from
the fact that instead of fixing the lengths of both sequences, we
fix the length of the common ancestor process. In a sense, this
situation is more realistic, since in practice the sequences are
hardly of exactly the same lengths; however, when they are related,
then the common ancestor must be of the same length for both of the
sequences. It turns out that the case of the random lengths the
statement of Theorems \ref{main1} and \ref{main2} hold with (full)
Hausdorff's distance $h$ instead of the restricted Hausdorff's
distance $h_o$. More precisely, Theorem~\ref{main3} states that
under the same assumptions as in Theorem~\ref{main1}, there exist
universal constants $C_r$ and $D_r$ so that
\[
P\bigl(h\bigl(H^r,L^r\bigr)>C_r\ln n\bigr)
\leq D_rn^{-2},
\]
where $h(H^r,L^r)$ stands for (full) Hausdorff's distance between
extremal alignments of random-length sequences.

Another measure could be the length of the biggest non-uniqueness
stretch, that is, the (horizontal) length between $*$'s. The simulations
in Section~\ref{sec:simulations} show that the length of the biggest
non-uniqueness stretch behaves similarly: the growth is almost
linear for the independent and logarithmic for the related
sequences. The latter has not been proven formally in this paper,
but we conjecture that it can be done using similar arguments as in
the proof of Theorems \ref{main1}, \ref{main2} and~\ref{main3}.
\end{Example*}
\subsection{The organization of the paper and the main
results}
\subsubsection{Preliminary results}
The paper is organized as follows. In Section~\ref{preliminaries},
the necessary notation is introduced and the extremal alignments are
formally defined and proven to exist (Proposition~\ref{high}). Also
some properties of the extremal alignments are proven. The section
also provides some
combinatorial bounds needed later.

Section~\ref{sec:ind} considers the case, where $X$ and $Y$ are
independent. The main result of the section is Theorem~\ref{thm:bound} that states for independent sequences the
Chvatal--Sankoff constant $\gamma$ satisfies the inequality
%
\begin{equation}
\label{ineq:31} \gamma\log_2p_o+2(1-\gamma)
\log_2q+2h(\gamma)\geq0,
\end{equation}
where
%
\begin{eqnarray}
\label{q-def} p_a&:=&P(X_i=a),\qquad  q:=1-\min
_{a}p_a,\nonumber\\[-8pt]\\[-8pt]
  p_o&:=&\sum
_{a\in\mathcal{A}}p^2_a,\qquad  h(p):=-p
\log_2p-(1-p)\log_2(1-p),\nonumber
\end{eqnarray}
that is, $h$ is the binary entropy function. The equality $\gamma
\log_2p_o+2(1-\gamma)\log_2q+2h(\gamma)=0$ has two solutions, hence,
as a byproduct, (\ref{ineq:31}) gives (upper and lower) bounds to
unknown $\gamma$. These bounds need not to be the best possible
bounds, but they are easy and universal in the sense that they hold
for any independent model. For example, taking the distribution of
$X_i$ and $Y_i$ uniform over the alphabet with $K$ letters (thus
$q=1-{1\over K}$ and $p_o={1\over K}$), we obtain the following
upper bounds $\bar{\gamma}$ to unknown $\gamma$. In the last row of
the table, the estimators $\hat{\gamma}$ of unknown $\gamma$ is
obtained via simulations. It is interesting to note that,
independently of $K$, the upper bound overestimates $\gamma$ about
the same amount. We also obtain the lower bounds, but for these
model the lower bounds are very close to zero and therefore not
informative.

In Section~\ref{sec:related}, the preliminary results for
related (homologous) sequences are presented. In Section~\ref{subsec:relateddef}, the formal definition of related sequences
are given. Our definition of relatedness is based on the existence
of common ancestor. Hence, our model models the homology in most
natural way. In our model, the related sequences $X_1,X_2,\ldots$
and $Y_1,Y_2,\ldots$
both consists of i.i.d. random variables, but the sequences are,
in general, not independent. Independence is a special case of the
model so that all results for
related sequences automatically hold for the independent ones. It is
also important to note that (unless the sequences are independent),
the two dimensional process $(X_1,Y_1),(X_2,Y_2),\ldots$ is not
stationary, hence also not ergodic. Hence, for the related
sequences ergodic theorems cannot be automatically applied. In
particular, Kingsman's subadditive ergodic theorem cannot be applied
any more to prove the convergence (\ref{sub}). This convergence as
well as the corresponding large deviation bound has been proven in
Section~\ref{subsec:properties}. Since we often consider the
sequences of unequal length, instead of (\ref{sub}), we prove a more
general convergence (Proposition~\ref{prop:piirv}):
%
\begin{equation}
\label{piirv-a} {L(X_1,\ldots,X_n;Y_1,\ldots,Y_{\lfloor na \rfloor})\over n}\to \gamma_{\tt{R}}(a),\qquad  \mbox{a.s.}
\end{equation}
Here $\gamma_{\tt{R}}(a)$ is a constant. We shall denote $\gamma
_{\tt{R}}(1)=:\gamma_{\tt{R}}$. From (\ref{piirv-a}), it
follows that for
any $a>0$,
\[
\gamma_{\tt{R}}(a)=a\gamma_{\tt{R}}\biggl({1\over a}
\biggr).
\]
%
If the sequences are independent and $a=1$, then $\gamma_{\tt
{R}}(a)=\gamma$. Corollary~\ref{cor} postulates the corresponding
large deviation result, stating that for every $\Delta>0$ there
exists $c>0$ such that for every $n$ big enough
%
\begin{equation}
\label{LD} P \biggl(\biggl|{L(X_1,\ldots,X_n;Y_1,\ldots,Y_{\lfloor na \rfloor})\over
n}-\gamma_{\tt{R}}(a)\biggr|>\Delta
\biggr)\leq\exp[-cn].
\end{equation}
In the \hyperref[app]{Appendix}, it is proven that $\gamma_{\tt{R}}(a)>\gamma
_{\tt{R}}$,
if $a>1$ and $\gamma_{\tt{R}}(a)<\gamma_{\tt{R}}$, if $a<1$ (Lemma~\ref{lemma:appendix}). That result together with (\ref{LD})
(obviously (\ref{piirv-a}) follows from (\ref{LD})) are the basic
theoretical tools for proving the main results of the paper,
Theorems \ref{main1}, \ref{main2} and \ref{main3}.

%
\begin{table}
\tablewidth=\textwidth
\tabcolsep=0pt
\caption{Upper bounds to Chvatal--Sankoff constant via inequality
(\protect\ref{ineq:31})}\label{tab1}
\begin{tabular*}{\textwidth}{@{\extracolsep{\fill}}llllllll@{}}
\hline
$K$ & 2 & 3 & 4 & 5 & 6 & 7 & 8 \\
\hline
$\bar{\gamma}$ & 0.866595 & 0.786473 & 0.729705 & 0.686117 & 0.650983
& 0.621719 & 0.596756 \\
$\hat{\gamma}$ & 0.81 & 0.72 & 0.66 & 0.61 & 0.57 & 0.54 & 0.52 \\
\hline
\end{tabular*}
\end{table}
%

\subsubsection{(Restricted) Hausdorff's distance and the main
results}\label{defHauss}
\paragraph*{Definition of (restricted) Hausdorff's distance}
We are interested in measuring the distance between the lowest and
highest alignment. One possible measure would be the maximum
vertical or horizontal distance (provided they are somehow defined).
However, those distances need not match the intuitive meaning of the
closeness of the alignment. For example, the following two
alignments (marked with $x$ and $o$, respectively) have a relatively
long maximal vertical distance (3), though they are intuitively
rather close:
%
\begin{equation}
\label{maxmin} %
\begin{tabular}{c}

\includegraphics{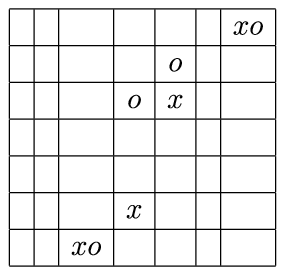}

\end{tabular}
\end{equation}
To overcome the problem, we measure the distance between two
alignments also in terms of Hausdorff's distance. More precisely,
let $U, V
\subset\{1,\ldots,n\}^2$,
be two alignments, both represented as sets of two-dimensional
points. The \textit{Hausdorff's distance between $U$ and $V$} is:
\[
h(U,V):=\max\Bigl\{\sup_{u\in U}\inf_{v\in V}d(u,v),
\sup_{v\in V}\inf_{u\in
U}d(u,v)\Bigr\},
\]
where $d$ is a distance in $\mathbb{R}^2$. In our
case, we take $d$ as the maximum-distance (but one can also consider
the usual Eucledian metric). We remark that Hausdorff's distance is
defined for any kind of sets. For the alignments in (\ref{maxmin}),
the Hausdorff's distance is obviously 1 (if $d$ were Euclidean, the
Hausdorff's distance would be $\sqrt{2}$).

Let now, for every $n$, $\alpha_n\in(0,1)$ be fixed, and we define
the subset $U_o\subseteq U$ consisting of those elements $(i,j)$ of
$U$ that have the first coordinate at least $n\alpha_n$ further from
$n$: $i\leq n(1-\alpha_n)$. Similarly, the subset $V_o\subset V$ is
defined. Formally, thus
\[
U_o:=\bigl\{(i,j)\in U: i \leq n(1-\alpha_n)\bigr\},\qquad
V_o:=\bigl\{(i,j)\in V: i \leq n(1-\alpha_n)\bigr\}.
\]
The \textit{restricted Hausdorff's distance} between
$U$ and $V$ is defined as follows:
\[
h_o(U,V):=\max\Bigl\{\sup_{u\in U_o}\inf
_{v\in V}d(u,v),\sup_{v\in
V_o}\inf
_{u\in
U}d(u,v)\Bigr\},
\]
where $d$ is a distance in $\mathbb{R}^2$. Clearly
$h_o(U,V)\leq h(U,V)$. Since in our case $U$ and $V$ are alignments
so that different points have different coordinates, the definition
of $h_o$ can be (somehow loosely) interpreted as a fraction
$\alpha_n$ of both alignments are left out when applying maximum in
Hausdorff's distance. We shall consider the case $\alpha_n\to0$.
Hence, the proportion of points left out decreases as $n$ grows.
\paragraph*{Sequences with fixed length}
We now state our main theorems for the sequences of fixed lengths.
Recall the definition of $p_o$ and $q$ from (\ref{q-def}). Let
%
\begin{equation}
\label{q} \p:=\max_{a\in\mathcal{A}}p_a,\qquad  \q:=1-\min
_{a,b\in\mathcal
{A}}P(X_1=a|Y_1=b,p=1).
\end{equation}
Here
$P(X_1=a|Y_1=b,p=1)$ is the conditional probability given that no
deletion occurs, or, in other words $X_1$ and $Y_1$ have the common
ancestor (see Section~\ref{subsec:relateddef} for formal
definition). Finally, let
\[
\rho:={p_0\q\over\p q}.
\]
%
In the following theorems, $h_o(L,H)$ stands for the restricted
Hausdorff's distance between alignments $L$ and $H$, both
represented as a set of 2-dimensional points. Recall that
Hausdorff's distance could be defined with the help in any metric in
$\mathbb{R}^2$. In the following, we shall consider both maximum and
$l_2$-norms. Throughout the paper, we shall use $\wedge$ and $\vee$
for $\min$ and $\max$, respectively.
%
\begin{theorem}\label{main1}
Let $X$ and $Y$ be related. Let $L,H$ be the (2-dimensional
representations of) lowest and highest alignments of $X$ and $Y$.
Assume
%
\begin{equation}
\label{maincondition1} \gamma_{\tt{R}}\log_2 \p+(1-
\gamma_{\tt{R}})\log_2 (q \q)+ \bigl((1-\gamma_{\tt{R}})
\wedge\gamma_{\tt{R}}\bigr)\log_2 (\rho\vee 1 )+ 2h(
\gamma_{\tt{R}})<0.
\end{equation}
Then there exist positive constants $M,C,D<\infty$ such that, for
$n$ big enough,
%
\begin{equation}
\label{inequalityC} P \bigl(h_o(L,H)>C\ln n \bigr)\leq
Dn^{-2},
\end{equation}
where $h_o$ is defined with
\[
\alpha_n:=M \sqrt{16 \ln n\over p n}
\]
and Hausdorff's distance is defined using maximum norm. If $h_o$ is
defined with respect to $l_2$ norm, then (\ref{inequalityC}) holds
with $C$ replaced by $\sqrt{2}C$.
\end{theorem}
%
For independent sequences $\q=q$, thus $\rho={p_o\over\p}\leq1$.
Then also $\gamma=\gamma_{\tt{R}}$ so that
\begin{eqnarray*}
&&\gamma_{\tt{R}}\log_2 \p+(1-\gamma_{\tt{R}})
\log_2 (q \q)+ \bigl((1-\gamma_{\tt{R}})\wedge
\gamma_{\tt{R}}\bigr)\log_2 (\rho\vee 1 )
\\
&&\quad = \gamma\log_2 \p+ 2(1-\gamma)\log_2 q \geq\gamma
\log_2 p_0 + 2(1-\gamma)\log_2 q \geq- 2h(
\gamma).
\end{eqnarray*}
The last inequality follows from (\ref{ineq:31}) (Theorem~\ref{thm:bound}). Hence, for unrelated (independent) sequences the
condition (\ref{maincondition1}) fails. It does not necessarily mean
that in this case (\ref{inequalityC}) holds not true, but based on
our simulations in Section~\ref{sec:simulations} we conjecture that
this is indeed the case.

In Theorem~\ref{main1}, we used the 2-dimensional representation of
alignments, so an alignment were identified with a finite set of
points. In the alignment graph, these points are joined by a line.
We consider the highest and lowest alignment graphs, and we are
interested in the maximal vertical (horizontal) distance between
these two piecewise linear curves. This maximum is called vertical
(horizontal) distance between lowest and highest alignment graphs.
The following theorem is stated in terms of vertical distance.
Clearly the same result holds for horizontal distance as well. In
the theorem, we shall also use the letters $L$ and $H$, but now they
stand for extremal alignment graphs rather than for the alignments
as the sets of the points. Since an alignment and the corresponding
alignment graph are very closely related, we hope that the notation
is not too ambiguous and the difference will be clear from the
context.
%
\begin{theorem}\label{main2}
Let $X$ and $Y$ be related. Let $L,H$ be the lowest and highest
alignment graphs of $X$ and $Y$. Assume (\ref{maincondition1}). Then
for $n$ big enough,
%
\begin{equation}
\label{inequality2C} P \Bigl(\sup_{x\in[0,n(1- \alpha_n)]}H(x)- L(x)>2C\ln n \Bigr)\leq
Dn^{-2},
\end{equation}
where the constants $C$, $D$ and $\alpha_n$ are the same as in
Theorem~\ref{main1}.
\end{theorem}
%
Hence, Theorems \ref{main1} and \ref{main2} state that when
$\gamma_{\tt{R}}$ is sufficiently bigger than (corresponding)
$\gamma$, then the distance between the extremal alignment (either
measured with restricted Hausdorff's metric or using alignment
graphs) grows no faster than logarithmically in $n$. Clearly,
$\gamma_{\tt{R}}$ is the bigger the more $X$ and $Y$ are related.
Hence, the inequality (\ref{maincondition1}) measures the degree of
the relatedness -- if this is big enough, then the distance between
extremal values grows (at most) logarithmically. Theorem~\ref{thm:bound} states that for independent sequence
(\ref{maincondition1}) fails, so that the assumptions of Theorems
\ref{main1} and \ref{main2} hold
for related sequences, only.

The fact that the distances between extremal alignments are measured
with respect to the restricted Hausdorff's distance, that is, so that a
small fraction of the alignments left out is obviously a bit
disappointing. Technically, this is due to the requirement that both
sequences are of the same length. As we shall see, this is not the
case when the lengths of the sequences are random. However, as also
the simulations in Section~\ref{sec:simulations} suggest, we believe
that the results of Theorems \ref{main1} and \ref{main2} hold also
when $h_o$ is replaced by the (full) Hausdorff's distance
$h$ and supremum is taken over $[0,n]$.

Theorems \ref{main1} and \ref{main2} are proven in Section~\ref{sec:proofs}. The proof is based on the observation that under
(\ref{maincondition1}) the probability that the sequences with
length about $n$ do not contain any related pairs is exponentially
small in $n$ (Lemmas \ref{related-pair} and \ref{related-pair2}).
Section~\ref{subsec:proprel} studies the location of the related
pairs in two dimensional representation. It turns out that with high
probability, the gaps between them are no longer then $A\ln n$,
where $A$ is suitable big constant. Applying these properties
together with Lemma~\ref{related-pair2}, we obtain that every
optimal alignment, including the extremal ones, cannot be far away
from the related points, since otherwise it would have a long piece
without any related pair contradicting Lemma~\ref{related-pair2}.
This argument is formalized in Lemma~\ref{abi} and Lemma~\ref{nurk}
in Section~\ref{sec:abi}. The formal proof of Theorems
\ref{main1} and \ref{main2} are given in Sections~\ref{sec:proof1} and \ref{sec:proof2}, respectively.
\paragraph*{The sequences with random length}
The related sequences are defined as follows: there is a common
ancestor process $Z_1,Z_2,\ldots$ consisting on $\mathcal{A}$-valued
i.i.d. random variables. Every letter $Z_i$ has a probability to
mutate according to a transition matrix that does not depend on $i$.
The mutations are independent of each other. After mutations, some
of the letters disappears. Thus, to every letter $Z_i$, there is
associated a Bernoulli random variables $D^x_i$ with $P(D_i^x=1)=p$.
When $D^x_i=0$, then the corresponding (mutated) letter disappears.
The deletions $D_1^x,D_2^x,\ldots$ are independent and the
remaining letters form the sequence $X_1,X_2,\ldots\,$. The $Y$
sequence is defined in the same way: every ancestor letter $Z_i$ has
another random mutation (independent of the all other mutations
including the ones that were used to define the $X$-sequence), and
independent i.i.d. deletions $D_i^y$ with the same probability. For
more detailed definition, see Section~\ref{subsec:relateddef}.

When dealing with the sequences of random length, we consider
exactly $m$ ancestors $Z_1,\ldots,Z_m$. Hence after deletions, the
length of obtained $X$-sequence is $n_x:=\sum_{i=1}^m D_i^x$ and
the length of $Y$-sequence is $n_y:=\sum_{i=1}^m D_i^y$. The
expected length of both sequences is thus $mp$ and we choose
$m(n):={n\over p}$ so that the expected length of the both sequences
is $n$. For simplicity, $m(n)$ is assumed to be integer. Thus, we
shall consider the sequences $X:=X_1,\ldots,X_{n_x}$ and
$Y:=Y_1,\ldots,Y_{n_y}$ of random lengths. It turns out that
mathematically this case is somehow easier so that the counterpart
of Theorem~\ref{main3} holds with full Hausdorff's distance $h$
instead of $h_o$.
%
\begin{theorem}\label{main3}
Let $X$ and $Y$ be the related sequences of random lengths. Let
$L,H$ be the (2-dimensional representation) of the highest and
lowest alignment. Assume (\ref{maincondition1}). Then there exist
constants $C_r$ and $D_r$ so that
%
\begin{equation}
\label{v3} P \bigl(h(H,L)>C_r\ln n\bigr)\leq D_r
n^{-2},
\end{equation}
where $h$ is the Hausdoff's distance with respect to maximum norm.
If $h$ is defined with respect to $l_2$ norm, then
(\ref{inequalityC}) holds with $C_r$ replaced by $\sqrt{2}C_r$.
\end{theorem}
%
From the proofs, it is easy to see that the random length analogue
of Theorem~\ref{main2} with $\alpha_n=0$ holds as well.
Theorem~\ref{main3} is proven in Section~\ref{random}.

Finally, in Section~\ref{sec:simulations}, some simulations about
the speed of the convergence are studied. The simulation clearly
indicate that for related sequences the growth of Hausdorff's and
vertical distance is at order of $\mathrm{O}(\ln n)$, hence the simulations
fully confirm the main results of the paper.

We would like to mention that to our best knowledge, the idea of
considering the extremal alignments as a characterization of the
homology has not been exploited, although the optimal and
sub-optimal alignments have deserved some attention before
\cite{macrounilargefluct,barder}. Therefore, the present paper as
the first step does not aim to minimize the assumptions or propose
any ready-made tests. These are the issues of the further research.
In a follow-up paper \cite{hirmo}, we apply some extremal-alignments
based characteristics to the real DNA-sequences, and compare the
results with standard sequence-alignment tools like BLAST.
\section{Preliminaries}\label{preliminaries}
Let $X_1,\ldots,X_{n_x}$ and $Y_1,\ldots,Y_{n_y}$ be two sequences
of lengths $n_x$ and $n_y$ from finite alphabet $\mathcal{A}=\{
0,1,\ldots, |\mathcal{A}|-1\}$. Let there exist two subsets of
indices $\{i_1,\ldots,i_k\}\subset\{1,\ldots,n_x \}$ and
$\{j_1,\ldots,j_k\}\subset\{1,\ldots,n_y \}$ satisfying
$i_1<i_2<\cdots<i_k$, $j_1<j_2<\cdots<j_k$ and $X_{i_1}=Y_{j_1},
X_{i_2}=Y_{j_2},\ldots, X_{i_k}=Y_{j_k}.$ Then $X_{i_1}\cdots
X_{i_k}$ is a common subsequence of $X$ and $Y$ and the pairs
%
\begin{equation}
\label{alignment} \bigl\{(i_1,j_1),\ldots,
(i_k,j_k)\bigr\}
\end{equation}
are (the 2-dimensional
representation of) the corresponding alignment. Let
\[
L(X_1,\ldots ,X_{n_x};Y_1,
\ldots,Y_{n_y})
\]
be the biggest $k$ such that there
exist such subsets of indices. The longest common subsequence is any
common subsequence with length $L(X_1,\ldots,X_{n_x};Y_1,\ldots
,Y_{n_y})$ and any alignment corresponding to a longest common
subsequence is called optimal.
In the following, we shall often consider the case, where, for some
constants $a,b>0$, $n_x=\lfloor bn \rfloor$, $n_y=\lfloor an
\rfloor$. Let us denote
\[
L_{bn,an}=L(X_1,\ldots,X_{\lfloor bn \rfloor};Y_1,
\ldots,Y_{\lfloor an
\rfloor}),\qquad  L_n:=L_{n,n}.
\]
Thus $L_n$ is the length of the
longest common sequence, when both sequences are of equal length,
$n_x=n_y=n$. The random variable $L_n$ is the main objet of
interest.
\subsection*{Extremal alignments: Definition and properties}
We now formally define the highest (optimal) alignment corresponding
to $L_n$. Let
\[
\bigl\{ \bigl(\bigl(i_1^1,j_1^1
\bigr),\ldots, \bigl(i_k^1,j_k^1
\bigr) \bigr),\ldots, \bigl(\bigl(i_1^{|A|},j_1^{|A|}
\bigr),\ldots, \bigl(i_k^{|A|},j_k^{|A|}
\bigr) \bigr) \bigr\}
\]
be the set of all optimal alignments. Hence, $k=L_n$ and $A=\{1,\ldots
, |A|\}$ is
the index set so that the elements of $A$ will be identified with
optimal alignments. For every $i^{\alpha}_l$ (resp.,~$j^{\alpha}_l$),
where $\alpha\in A$\vadjust{\goodbreak} and $l\in\{1,\ldots,k\}$, we shall denote
$j(i^{\alpha}_l):=j^{\alpha}_l$ (resp.,~$i(j^{\alpha}_l):=i^{\alpha}_l$). We define
\[
J:=\bigl\{j_l^{\alpha}:\alpha\in A, l=1,\ldots,k\bigr\},\qquad  I:=
\bigl\{i_l^{\alpha
}:\alpha\in A, l=1,\ldots,k\bigr\}.
\]
Let $j^h_k:=\max_{\alpha}j^{\alpha}_k=\max J$.
There might be many alignments $\alpha$ such that
$j_k^{\alpha}=j^h_k$. Among such alignments take $i_k^h$ to be
minimum. Formally,
$i_k^h=\min\{i^{\alpha}_k: j^{\alpha}_k=j^h_k\}.$ After fixing
$(i_k^h,j_k^h)$, we take $j_{k-1}^h$ as the biggest $j\in J$ such
that the corresponding $i$ is smaller than $i^h_k$. Formally,
\[
j^h_{k-1}:=\max\bigl\{j_l^{\alpha}: i
\bigl(j_l^{\alpha}\bigr)<i_k^h,\alpha
\in A, l=1,\ldots,k \bigr\}.
\]
There might be several $i$'s
such that corresponding $j$ is $j^h_{k-1}$. Amongst them, we choose
the minimum. Thus,
\[
i^h_{k-1}:=\min\bigl\{i_l^{\alpha}:j
\bigl(i_l^{\alpha}\bigr)=j^h_{k-1},
\alpha\in A, l=1,\ldots,k \bigr\}.
\]
Proceeding so, we obtain an alignment. We call this the \textit
{highest alignment procedure}. We now prove that the procedure can be
repeated $k$ times, that is, the obtained alignment is optimal.
%
\begin{proposition}\label{high}
The highest alignment procedure produces an optimal alignment
\[
\bigl\{\bigl(i_1^h,j_1^h\bigr),
\ldots,\bigl(i_k^h,j_k^h\bigr)
\bigr\},
\]
where $(i_t^h,j_t^h)$ can
be obtained as follows
%
\begin{equation}
\label{l} j^h_{t}:=\max\bigl\{j^{\alpha}_{t}:
\alpha\in A\bigr\},\qquad  i^h_{t}:=\min\bigl
\{i^{\alpha}_{t}: j\bigl(i^{\alpha}_{t}\bigr)=
j^h_{t}, \alpha \in A\bigr\},\qquad  t=1,\ldots,k.
\end{equation}
\end{proposition}
%
%
\begin{pf}
Clearly the pair $(i^h_{k},j^h_{k})$ is the last pair of an optimal
alignment, that is, there exists $\alpha\in A$ such that
$(i^h_{k},j^h_{k})=(i_k^{\alpha},j_k^{\alpha})$. So (\ref{l}) holds
with $t=k$. Similarly, there exists a $\beta\in A$ such that
$j^h_{k-1}=j^{\beta}_{k-1}$. Let us show this. There exists a
$\beta$ such that $j^h_{k-1}=j^{\beta}_l$, we have to show that
$l=k-1$. Note that $l$ cannot be $k$, since otherwise
$(i_1^{\beta},j_1^{\beta}),\ldots,(i_k^{\beta},j_k^{\beta
}),(i_k^{h},j_k^{h})$
would be an alignment of length $k+1$. Suppose $l\leq k-2$. Since
$j^{\beta}_{l}<j^{\beta}_{k-1}<j^{\beta}_{k}\leq j^{h}_k=\max J$, by
definition of $j^{h}_{k-1}$, it must be that $i^{h}_k \leq
i^{\beta}_{k-1}.$ Since $i^{h}_k=i^{\alpha}_k>i^{\alpha}_{k-1}$, we
have that $i^{\alpha}_{k-1}<i^{\beta}_{k-1}<i^{\beta}_{k}$. On the
other hand, $j_{k-1}^{\alpha}\leq
j_{k-1}^{h}=j_{l}^{\beta}<j_{k-1}^{\beta}$ implying that
$j_{k-1}^{\alpha}<j^{\beta}_{k-1}<j^{\beta}_{k}$. Hence,
$(i_1^{\alpha},j_1^{\alpha}),\ldots,(i_{k-1}^{\alpha
},j_{k-1}^{\alpha}),(i_{k-1}^{\beta},j_{k-1}^{\beta}),
(i_{k}^{\beta},j_{k}^{\beta})$ would be an alignment of length
$k+1$. Therefore, $j^h_{k-1}=\max\{j^{\alpha}_{k-1}:
i^{\alpha}_{k-1}<i^h_k, \alpha\in A\}$. Let us now prove that
(\ref{l}) with $t=k-1$ holds. If this were not the case, then
$j^h_{k-1}<\max\{j^{\alpha}_{k-1}: \alpha\in A\}$. This implies the
existence of $\beta$ so that $j^{\beta}_{k-1}>j^h_{k-1}$ and
$i^{\beta}_{k-1}\geq i^h_k>i^h_{k-1}$. But as we saw, those
inequalities would give an alignment with the length $k+1$. This
concludes the proof of (\ref{l}) with $t=k-1$. For $t=k-2,\ldots,1$
proceed similarly.
\end{pf}
%
%
\begin{figure}

\includegraphics{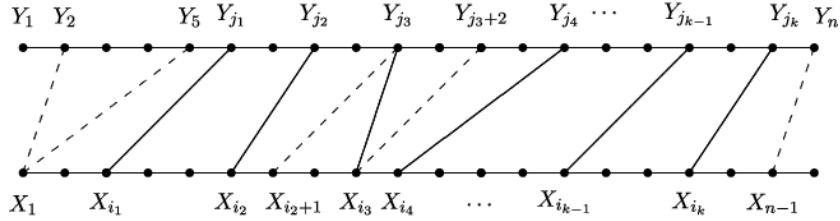}

\caption{An example of the highest
alignment.}\label{fig3}
\end{figure}
%

Figure~\ref{fig3} is an example of an highest alignment. The solid
lines are aligned pairs (the upper-index $h$ is dropped from the
notation). If $Y_{j_3+2}=Y_{j_3}$, then, as showed by dashed line,
$X_{i_3}$ could be aligned with $Y_{j_3+2}$ that contradicts the
highest alignment procedure. Thus, every $Y_{j_t}$ in the highest
alignment is different from all $Y_j$ that are right after $Y_{j_t}$
and before $Y_{j_{t+1}}$. This observation is postulated as
statements (\ref{cor1}) and (\ref{cor2}) in the following corollary.
Similarly, if $X_{i_2+1}=X_{i_3}$, then, as showed by dashed line,
$X_{i_2+1}$ could be aligned with $Y_{j_3}$ that also contradicts
the highest alignment procedure. Thus, in the highest alignment all
$X_i$-s right after $X_{i_{t-1}}$ and before $X_{i_{t}}$ must differ
from $X_{i_t}$. This observation is formulated as the statements
(\ref{cor1a}) and (\ref{cor2a}) in the following corollary. In the
highest alignment, typically, $i_1<j_1$ and $j_k>i_k$. If $X_1$ is
not aligned, then it must be that $Y_i\ne X_1$ for
$i=1,\ldots,j_1-1$, otherwise they could be aligned (as showed by
dashed line) contradicting the optimality. These observations are
statements (\ref{cor4a}) and (\ref{cor4b}) in the following
corollary. Similarly, if $Y_n$ is not aligned, it should be
different from all $X_{i_k+1},\ldots,X_n$. These observations are
statements (\ref{cor5a}) and (\ref{cor5b}) in the following
corollary.\looseness=-1
%
\begin{corollary}\label{corI} The highest alignment has the following
properties:\vspace*{-1pt}
%
\begin{eqnarray}
\label{cor3}  &&X_{i^h_t}=Y_{j^h_t}, \qquad t=1,\ldots,k;
\\[-1pt]
\label{cor1} &&Y_{j_t^h}\ne Y_j,\qquad  j_t^h<j<j_{t+1}^h,\qquad
t=1,\ldots,k-1;
\\
\label{cor2} &&Y_{j_k^h}\ne Y_j,\qquad  j_k^h<j
\leq n;
\\
\label{cor1a} &&X_{i_t^h}\ne X_i,\qquad  i_t^h>i>i_{t-1}^h,\qquad
t=2,\ldots,k;
\\
\label{cor2a} &&X_{i_1^h}\ne X_i, \qquad 1\leq
i<i^h_1;
\\
\label{cor4a} && \mbox{if }i^h_1>1, \mbox{ then }
X_1\ne Y_j, \qquad 1\leq j<j^h_1;
\\
\label{cor4b} && \mbox{if } j^h_1>1, \mbox{ then }
Y_1\ne X_i,\qquad  1\leq i <i^h_1;
\\
\label{cor5a} && \mbox{if }n>j^h_k, \mbox{ then }
Y_n\ne X_i, \qquad i^h_k<i\leq n;
\\
\label{cor5b} && \mbox{if } n>i^h_k, \mbox{ then }
X_n\ne Y_j,\qquad  j^h_k<j\leq n.
\end{eqnarray}
\end{corollary}
\begin{pf} The equalities (\ref{cor3}) are
obvious. Suppose that for a $t=1,\ldots,k-1$ there exists an index
$j$ such that $j_t^h<j<j_{t+1}^h$ and $Y_{j^h_t}=Y_j$. Then the
pairs
\[
\bigl\{\bigl(i^h_1,j_1^h\bigr),
\ldots ,\bigl(i^h_{t-1},j_{t-1}^h
\bigr),\bigl(i^h_{t},j\bigr),\bigl(i^h_{t+1},j_{t+1}^h
\bigr),\ldots ,\bigl(i^h_k,j_k\bigr)\bigr\}
\]
would correspond to an optimal alignment, say $\beta$, satisfying
\[
j_t^{\beta}=j>j_t^h=\max\bigl
\{j^{\alpha}_{t}: \alpha\in A\bigr\}.\vadjust{\goodbreak}
\]
Thus, (\ref{cor1}) holds. The same argument proves (\ref{cor2}),
(\ref{cor1a}) and (\ref{cor2a}). If one of the inequalities in
(\ref{cor4a})--(\ref{cor5b}) is not fulfilled, then it would be
possible to align one more pair without disturbing already existing
aligned pairs. This contradicts the optimality.
\end{pf}
%


One can also think of the left-most or nord-west alignment. It could
be defined as an alignment
$\{(i_1^w,j_1^w),\ldots,(i_k^w,j_k^w)\}$, where $i_1^w=\min I$,
$j_1^w=\max\{j_l^{\alpha}: i(j_l^{\alpha})=i_1^w, \alpha\in
A, l=1,\ldots,k\}$ and for every $t=2,\ldots,k$,
\begin{eqnarray*}
i_t^w&:=&\min\bigl\{i_l^{\alpha}:j
\bigl(i_l^{\alpha}\bigr)>j^w_{t-1},
\alpha\in A, l=1,\ldots,k\bigr\},\\
 j_t^w&:=&\max\bigl
\{j_l^{\alpha}:i\bigl(j_l^{\alpha}
\bigr)=i_t^w, \alpha\in A, l=1,\ldots,k\bigr\}.
\end{eqnarray*}
Here the
superscript ``$w$'' stands for west. By the analogue of Proposition~\ref{high},
%
\begin{equation}
\label{l-right} i^w_{t}=\min\bigl\{i^{\alpha}_{t}:
\alpha\in A\bigr\},\qquad  j^w_{t}=\max\bigl
\{j^{\alpha}_{t}: i\bigl(j^{\alpha}_{t}\bigr)=
i^w_{t},\alpha\in A\bigr\},\qquad  t=1,\ldots,k.
\end{equation}
Using (\ref{l}) and (\ref{l-right}), it
is easy to see that the left-most and highest alignments actually
coincide. Indeed, by (\ref{l}) and (\ref{l-right}), $j_t^h\geq
j_t^w$ and $i_t^w\leq i_t^h$, $\forall t$. If, for a $t$,
$(i^h_t,j^h_t)\ne(i^w_t,j^w_t)$, then, by the definitions, both
inequalities have to be strict, that is, $i_t^w<i_t^h$ and $j_t^w<j_t^h$.
To see this, suppose $i_t^w=i_t^h$. This means that there exists an
alignment $\alpha$, such that $i^{\alpha}_t=i_t^w$ and
$j^{\alpha}_t=j_t^h$. This, in turn, implies that
\[
\max\bigl\{j^{\alpha}_t: i\bigl(j_t^{\alpha}
\bigr)=i_t^w, \alpha\in A \bigr\}=\max\bigl\{
j^{\alpha}_t, \alpha\in A \bigr\},
\]
that is, $j_t^w=j_t^h$. The same argument shows that if
$j_t^w=j_t^h$, then also $i_t^w=i_t^h$. Thus $(i^h_t,j^h_t)\ne
(i^w_t,j^w_t)$ implies that $i_t^w<i_t^h$ and $j_t^w<j_t^h$. These
inequalities, however, would imply the
existence of an alignment with the length $k+1$.

The lowest (the right-most) alignment
$\{(i_1^l,j_1^l),\ldots,(i_k^l,j_k^l)\}$ will be defined similarly:
%
\begin{equation}
\label{lowest} j^l_{t}:=\min\bigl\{j^{\alpha}_{t}:
\alpha\in A\bigr\}, \qquad i^l_{t}:=\max\bigl
\{i^{\alpha}_{t}: j\bigl(i^{\alpha}_{t}\bigr)=
j^l_{t}, \alpha \in A \bigr\}, \qquad t=1,\ldots,k.
\end{equation}
%
\begin{Remark*} Note that the left-most alignment equals the lowest
alignment of $(Y_n,\ldots,Y_1)$ and $(X_n,\ldots,X_1)$ implying that
the highest alignment of $(X_1,\ldots,X_n)$ and
$(Y_1,\ldots,Y_n)$ equals to the lowest alignment of
$(Y_n,\ldots,Y_1)$ and $(X_n,\ldots,X_1)$. Thus, the lowest
alignment between $(X_1,\ldots,X_n)$ and $(Y_1,\ldots,Y_n)$ can be
defined as the highest alignment between $(Y_n,\ldots,Y_1)$ and
$(X_n,\ldots,X_1)$.
\end{Remark*}

\subsection*{Combinatorics} Another way to study an alignment of $
X_1,\ldots, X_{n_x}$ and $Y_1,\ldots, Y_{n_y}$ is to present it as a
strictly increasing mapping
%
\begin{equation}
\label{kujutis} v: \{1,\ldots,n_x\} \hookrightarrow\{1,
\ldots,n_y\}.
\end{equation}
Notation (\ref{kujutis}) means: There exists $I(v)\subset
\{1,\ldots,n_x\}$ and a mapping
\[
v: I \to\{1,\ldots,n_y\}
\]
such that $Y_{v(i)}=X_i,$ $\forall i\in I$ and $v$ is strictly
increasing: $v(i_2)>v(i_1)$, if $i_2 > i_1$. The length of $v$ is
denoted as $|v|$. In the notation of previous sections, thus,
$j_t:=v(i_t)$, $t=1,\ldots, |v|$.

Consider now the case $n_x=n_y=n$, that is, both sequences are of length
$n$. Let then $V_k$ be the set of all alignments with length $k$.
Formally,
\[
V_k:= \bigl\{v: \{1,\ldots,n\} \hookrightarrow\{1,\ldots,n\} :
|v|=k \bigr\}.
\]
Fix $\Delta>0$, $\gamma\in(0,1)$ and let
%
\begin{equation}
\label{wdef} W_n(\gamma,\Delta):=\bigcup
_{k=(\gamma-\Delta)n}^{(\gamma+\Delta)n} V_k.
\end{equation}
Hence, $W_n$ consists of these alignments that
have length not smaller that $(\gamma-\Delta)n$ and not bigger that
$(\gamma+\Delta)n$. In the subsequent sections, we shall show that
there exists a constant $\gamma$ (depending on the model) so that
for $n$ big enough all optimal alignments belong to $W_n$ with high
probability. Thus, in a sense the set $W_n$ contains all alignments
of interest. We are interested in bounding the size of that set. For
that, we use the binary entropy function
\[
h(p):=-p\log_2p-(1-p)\log_2(1-p).
\]
Let, for $\gamma,\Delta\in(0,1)$ such that
$0<\gamma-\Delta,\gamma+\Delta<1$
%
\begin{equation}
\label{gamma1} H(\gamma,\Delta):=\max_{\alpha\in[\gamma-\Delta,\gamma+\Delta]} h(\alpha).
\end{equation}
%
Since
\[
{{n}\choose{pn}}\leq2^{h(p)n},
\]
for every
%
\begin{equation}
\label{k} (\gamma-\Delta)n\leq k\leq(\gamma+\Delta)n,
\end{equation}
it holds
\[
|V_k|={{n}\choose{{k\over n}n}}^2
\leq2^{2H(\gamma,\Delta)n}.
\]
Hence, the number of alignments in $W_n$ can be bounded as follows:
%
\begin{equation}
\label{boundVn} \bigl|W_n(\gamma,\Delta)\bigr|\leq2\Delta n 2^{2H(\gamma,\Delta)n}.
\end{equation}
Let us consider now a more general case $n_y>n_x$. Denote
$m=n_y>n_x=n$. Assume that $m\leq n(1+\Delta)$. Then
\[
|V_k|={{n}\choose{{k\over n}n}} {{m}\choose{
{k\over m}m}}\leq 2^{h({k/
n})n+h({k/ m})n(1+\Delta)}.
\]
Instead of (\ref{k}), we assume
$k$ to satisfy
%
\begin{equation}
\label{kk} \gamma-\Delta\leq{k\over n}\leq\gamma+2\Delta.
\end{equation}
Then
\[
\gamma-2\Delta\leq{\gamma-\Delta\over1+\Delta}\leq{k\over m} \leq
{k\over n} \leq\gamma+2\Delta
\]
and
\[
2^{h({k/
n})n+h({k/ m})n(1+\Delta)}\leq 2^{H(\gamma,2\Delta)n+H(\gamma,2\Delta)n(1+\Delta)}= 2^{H(\gamma,2\Delta)n(2+\Delta)}.
\]
In this case, defining
\[
W_{n,m}(\gamma,\Delta):=\bigcup_{k=(\gamma-\Delta)n}^{(\gamma
+2\Delta)n}
V_k,
\]
it holds
\[
|W_{n,m}|\leq 3\Delta n 2^{(2+\Delta)H(\gamma,2\Delta)n}.
\]
%
\section{Independent sequences}\label{sec:ind} In this
section, only, let $X=X_1,\ldots,X_n$ and $Y=Y_1,\ldots,Y_n$ be two
independent i.i.d. sequences from the alphabet $\mathcal{A}$. Recall
that for any $a>0$, $L_{an,n}=L(X_1,\ldots,X_{\lfloor an
\rfloor};\break
Y_1,\ldots,Y_n)$ and $L_n=L_{n,n}$. By the Kingman's subbadditive
ergodic theorem, there exists a constant $\gamma(a)\in(0,1]$ so that\vspace*{-1pt}
\[
{L_{an,n}\over n}\to\gamma(a) \qquad \mbox{a.s. and in }L_1.
\]
We
shall denote $\gamma:=\gamma(1)$, the constant $\gamma$ is often
called the \textit{Chvatal--Sankoff} constant. In the \hyperref[app]{Appendix}, it will
be shown that when $a<1$, then $\gamma(a)<\gamma$ (Lemma~\ref{lemma:appendix}).

Note that $L_{an,n}$ is a function of $n(1+a)$ i.i.d. random
variables. Clearly, changing one of the variables changes the value
of $L_n$ at most by one, so that by McDiarmid inequality (see, e.g.,
\cite{rate}), for every $\Delta>0$\vspace*{-1pt}
%
\begin{equation}
\label{McD} P \bigl(|L_{an,n}-E L_{an,n}|>n\Delta \bigr)\leq2\exp
\biggl[-{2 \Delta^2
\over(1+a)}n\biggr].
\end{equation}
Take $n_o(\Delta,a)$ so big that $|{EL_{an,n} \over
n}-\gamma(a)|<{\Delta\over2}$. Then\vspace*{-1pt}
%
\begin{eqnarray}\label{LDV2ind}
&&P \bigl(\bigl|L_{an,n}-\gamma(a) n\bigr|\geq n\Delta \bigr)\nonumber\\[-0.5pt]
&&\quad \leq P
\bigl(|L_{an,n}-EL_{an,n}|+\bigl|EL_{an,n}-\gamma(a) n\bigr|\geq
n\Delta \bigr)\\[-0.5pt]
&&\quad  \leq P \biggl(|L_{an,n}-EL_{an,n}|\geq n
{\Delta\over
2} \biggr)\leq 2\exp\biggl[-{\Delta^2\over2(1+a)} n
\biggr], \qquad n>n_o.\nonumber
\end{eqnarray}
Taking $a=1$, we see the existence of $n_o(\Delta)$ so that
for $n>n_o$ with high probability all optimal alignments
are contained in the set $W_n(\gamma,\Delta)$ as defined
in (\ref{wdef}).

Recall that for any optimal alignment $v$, $(i^h_1, j_1^h)$ and
$(i^h_{|v|}, j_{|v|}^h)$ are the first and last pairs of indexes of
the highest alignment of $X$ and $Y$. We consider the random
variables $S:=j^h_1-1, T:=n-i^h_{|v|}$. Clearly $S$ and $T$ have the
same law. The following proposition states that for any $c\in
(0,1)$, the probabilities $P(S>cn)=P(T>cn)$ decrease exponentially
fast.
%
\begin{proposition}\label{STind}
Let $c\in(0,1)$. Then there exists constant $d(c)>0$, so that, for
$n$ big enough, $P(T>cn)=P(S>cn)\leq\exp[-dn]$.
\end{proposition}
\begin{pf}
Note that $\{T>cn\}\subset\{L_{(1-c)n,n}=L_{n}\}$ and for any
$\bar{\gamma}$,\vspace*{-1pt}
\[
\{L_{(1-c)n,n}=L_{n}\}\subset\{L_{(1-c)n,n}\geq \bar{
\gamma}n\}\cup\{L_{n}\leq\bar{\gamma}n\}.
\]
Let $a:=1-c$. By
Lemma~\ref{lemma:appendix}, $\gamma>\gamma(a)$. Let\vspace*{-1pt}
\[
\bar{\gamma}:={\gamma_+\gamma(a)\over2}, \qquad \Delta:=\gamma-\bar {\gamma}=\bar{
\gamma}-\gamma(a).
\]
Use
(\ref{LDV2ind}) to see that for $n$ big enough,\vspace*{-1pt}
\begin{eqnarray*}
P(T>cn)&\leq& P (L_{an,n}\geq\bar{\gamma}n )+P (L_{n}\leq
\bar{\gamma }n )\\
&=&P \bigl(L_{an,n}\geq\bigl(\gamma(a)+\Delta\bigr)n
\bigr)+P \bigl(L_{n}\leq({\gamma}-\Delta)n \bigr)
\\
&\leq&2\exp\biggl[-{\Delta^2\over2(1+a)} n\biggr]+2\exp\biggl[-
{\Delta^2\over4} n\biggr].
\end{eqnarray*}
This concludes the proof.
\end{pf}

Recall the definition of $q$ and $p_o$ in (\ref{q-def}). Note that
for independent sequences, $p_o=P(X_i=Y_i)$. The following lemma
bounds the probability that an alignment $v\in V_k$ is the highest
optimal alignment.
%
\begin{lemma}\label{lemmaB} Let $v\in V_k$. Let $B(v)$ be the event
that $v$ is
the highest optimal alignment of $X$ and~$Y$. Then
%
\begin{equation}
\label{Bhinnang} P\bigl(B(v)\bigr)\leq p_o^kq^{2(n-k)-(j_1-1)-(n-i_k)}.
\end{equation}
\end{lemma}
\begin{pf}
Let $v\in V_k$ be an alignment. We denote by $i_1,\ldots, i_k$ the
elements of $I(v)$ and we define $j_t:=v(i_t)$, $t=1,\ldots,k$.

Since all random variables $X_1,\ldots, X_n, Y_1,\ldots, Y_n$ are
independent, by Corollary~\ref{corI} the probability of $B(v)$ could
be estimated as follows
\begin{eqnarray*}
\label{Bprob} P \bigl(B(v) \bigr)\leq&\displaystyle \prod_{t=1}^k
P\bigl(B_t(v)\bigr),
\end{eqnarray*}
where, for $t=2,\ldots,k-1$
\[
B_t(v):=\{X_{i_t}=Y_{j_t};Y_j\ne
Y_{j_t}, j_t<j<j_{t+1}; X_i\ne
X_{i_t}, i_{t-1}<i<i_t\}
\]
and
\begin{eqnarray*}
&&B_1(v):=\lleft\{ %
\begin{array} {l@{\qquad}l}
\{X_{i_1}=Y_{j_1};X_i\ne X_{i_1}, i
<i_1;X_1\ne Y_j, \\
\hphantom{\{}j <j_1;
Y_j\ne Y_{j_1}, j_1<j<j_2\}, &
\mbox{if $i_1>1$;}
\\
\{X_{i_1}=Y_{j_1}; Y_j\ne Y_{j_1},
j_1<j<j_2\}, & \mbox{if $i_1=1$.} \end{array}
\rright.
\\
&&B_k(v):=\lleft\{ %
\begin{array} {l@{\qquad}l}
\{X_{i_k}=Y_{j_k};X_i\ne X_{i_k},
i_{k-1}<i<i_k; \\
\hphantom{\{}Y_j\ne Y_{j_k},
j_k<j; X_i\ne Y_n, i_k<i\}, &
\mbox{if $j_k>n$;}
\\
\{X_{i_k}=Y_{j_k};X_i\ne X_{i_k},
i_{k-1}<i<i_k\}, & \mbox{if $j_k=n$.}
\end{array} %
\rright.
\end{eqnarray*}
By independence, clearly for $t=2,\ldots,k-1$,
\[
P\bigl(B_t(v)\bigr)=\sum_a
p^2_a(1-p_a)^{i_t-i_{t-1}-1+j_{t+1}-j_t-1}\leq
p_oq^{i_t-i_{t-1}+j_{t+1}-j_t-2}.
\]
For the
events $B_1(v)$ and $B_k(v)$, we estimate
\begin{eqnarray*}
P\bigl(B_1(v)\bigr)&\leq&\lleft\{ %
\begin{array} {l@{\quad}l}
p_oq^{j_2-j_1-1}, & \mbox{if $i_1=1$;}
\\
p_oq^{j_2-j_1-1+j_1-1+i_1-1}, & \mbox{if $i_1>1$.} \end{array}
\rright.
\\
P\bigl(B_k(v)\bigr)&\leq&\lleft\{ %
\begin{array} {l@{\quad}l}
p_oq^{i_k-i_{k-1}-1}, & \mbox{if $j_k=n$;}
\\
p_oq^{i_k-i_{k-1}-1+n-j_k+n-i_k}, & \mbox{if $j_k<n$.} \end{array}
\rright.
\end{eqnarray*}
These equations yield (\ref{Bhinnang}). Note that in
(\ref{Bhinnang}), the term $(n-i_k)$ disappears when $j_k<n$ and the
term $(j_1-1)$ disappears, when $i_1>1$.
\end{pf}
%

Our first main result is a bound to the unknown Chvatal--Sankoff
constant $\gamma$.
%
\begin{theorem}\label{thm:bound}
Let $X_1,X_2,\ldots$ and $Y_1,Y_2,\ldots$ be two independent i.i.d.
sequences with the same distribution. Let $\gamma$ be the
corresponding Chvatal--Sankoff constant. Then the following condition
holds
%
\begin{equation}
\label{cond:ind+} \gamma\log_2p_o+2(1-\gamma)
\log_2q+2h(\gamma)\geq0.
\end{equation}
\end{theorem}
%
%
\begin{pf}
The proof is based on the contradiction: assuming that
(\ref{cond:ind+}) fails leads to the existence of constants
$c>0,b>0$ (independent of $n$) such that for $n$ big enough,
$P(S+T\leq2cn)\leq\exp[-bn].$ Then, for big $n$,
\[
1-\exp[-bn]\leq P(S+T > 2cn)\leq P(S>cn)+P(T>cn)
\]
contradicting Proposition~\ref{STind}.

If (\ref{cond:ind+}) is not fulfilled, then it is possible to find
constants $\Delta>0, c >0$ so small that
%
\begin{equation}
\label{cond:ind2} -b_1:=(\gamma-\Delta)\log_2p_o+2(1-
\gamma-\Delta-c)\log _2q+2H(\gamma, \Delta)<0.
\end{equation}
Fix now $\Delta>0,c>0$ so small that
(\ref{cond:ind2}) holds. Let
\[
E_{\Delta}:=\bigl\{|L_n-n\gamma|<n\Delta\bigr\}.
\]
When $E_{\Delta}$ holds, then all optimal alignments belong to the
set $W_n:=W_n(\gamma,\Delta)$. By Lemma~\ref{lemmaB}, for every
$v\in W_n$
%
\begin{equation}
P\bigl(B(v)\bigr)\leq p_o^{n(\gamma-\Delta)}q^{2n(1-\gamma-\Delta)-(n-i_k)-(j_1-1)}.
\end{equation}
Note that $\bigcup_{v\in W_n}B(v)=E_{\Delta}$. Let, for every $v$,
$s(v):=j_1-1$ and $t(v):=n-i_ {|v|}$. Then by (\ref{LDV2ind}) and
(\ref{boundVn}), there exists $b>0$ (independent of $n$) so that for
$n$ big enough
\begin{eqnarray*}
P(S+T\leq2cn)&\leq&\sum_{v\in W_n: s(v)+t(v)\leq
2cn}P\bigl(B(v)\bigr)+P
\bigl(E_{\Delta}^c\bigr)
\\
&\leq&2\Delta n 2^{n (2H(\gamma,\Delta)+(\gamma-\Delta)\log_2
p_o+2(1-\gamma-\Delta-c)\log_2q )}+P\bigl(E_{\Delta}^c\bigr)
\\
&\leq& 2\Delta n 2^{-b_1n}+2\exp\biggl[-{\Delta^2\over4}n\biggr]
\leq\exp[-b n].
\end{eqnarray*}
\upqed
\end{pf}
%
\section{Related sequences: Definition and theory}\label{sec:related}
\subsection{Definition of relatedness}\label{subsec:relateddef}
Let us now define the relatedness of the sequences $(X,Y)$. Our
concept of relatedness is based on the assumption that there exists
a common ancestor, from which both sequences $X$ and $Y$ are
obtained by independent random mutations and deletions. In the
following, the common ancestor is an $\mathcal{A}$-valued i.i.d.
process $Z_1,Z_2,\ldots\,$. We could imagine that $X$ and $Y$ is the
genome of two species whilst $Z$ is the genome of a common ancestor.
In computational linguistics, $X$ and $Y$ could be words from two
languages which both evolved from the word $Z$ in an ancient
language.

A letter $Z_i$ has a probability to mutate according to a transition
matrix that does not depend on~$i$. Hence, a mutation of the letter
$Z_i$ can be formalized as $f(Z_i, \xi_i)$, where $f\dvtx  \A\times
\mathbb{R}\to\A$ is a mapping and $\xi_i$ is a standard normal
random variable. The mapping $f_i(\cdot):=f(\cdot, \xi_i)$ from $\A$
to $\A$ will be referred as the random mapping. The mutations of the
letters are assumed to be independent. This means that the random
variables $\xi_1,\xi_2,\ldots$ or the random mappings
$f_1,f_2,\ldots$ are independent (and identically distributed).
After mutations, the sequence is $f_1(Z_1),f_2(Z_2),\ldots.$ Some of
its elements disappear. This is modeled via a deletion process
$D^x_1,D^x_2,\ldots$ that is assumed to be an i.i.d. Bernoulli
sequence with parameter $p$ that is, $P(D^x_i=1)=p$. If $D^x_i=0$, then
$f_i(Z_i)$ is deleted. The resulting sequence, let it be $X$, is,
therefore, the following: $X_i=f_j(Z_j)$ if and only if $D^x_j=1$
and $\sum_{k=1}^{j}D^x_k=i$. We call the index $j$ the \textit{ancestor
of $i$}, it shall be denoted by $a^x(i)$. The mapping $a^x$ depends
on the deletion process $D^x$, only. Now
\[
X_i=f_{a^x(i)}(Z_{a^x(i)}), \qquad i=1,\ldots,n.
\]
Similarly, the sequence $Y$ is obtained from $Z$. For mutations, fix
an i.i.d. standard normal sequence $\eta_1,\eta_2,\ldots$ so that
the mutated sequence is $ g_1(Z_1), g_2(Z_2),\ldots$ with $
g_i(\cdot):=f(\cdot, \eta_i).$ Note that the transition matrix
corresponding to $Y$-mutations equals the one corresponding to
$X$-mutations implying that the random mappings $g_i$ and $f_i$ have
the same distribution. Since the mutations of $X$ and $Y$ are
supposed to be independent, we assume the sequences $\xi$ and $\eta$
or the random mappings sequences $f_1, f_2,\ldots$ and
$g_1,g_2,\ldots$ are independent. Note that then the pairs
$(f_1(Z_1),g_1(Z_1)),(f_2(Z_2),g_2(Z_2)),\ldots$ are independent,
but $f_i(Z_i)$ and $g_i(Z_i)$, in general, are not. Finally,
\[
Y_i=f_{a^y(i)}(Z_{a^y(i)}),
\]
where, as previously, $a^y(i)=j$ if and only if $D^y_j=1$ and
$\sum_{k=1}^{j}D^y_k=i$. Here, $D^y_1,D^y_2,\ldots$ is an i.i.d.
Bernoulli sequence with the same parameter as $D^x$ but independent
of $D^x.$ Hence, the deletions of $Y$ and $X$ are independent.

In the following, we shall call the sequences $X=X_1\ldots X_n$ and
$Y=Y_1\ldots Y_n$ \textit{related}, if they follow the model described
above. Note that for the related sequences, the random variables
$X_1,X_2,\ldots$ as well as $Y_1,Y_2,\ldots$ are still i.i.d., but
these two sequences are, in general, not independent any more. As
mentioned above, the process $(X_1,Y_1),(X_2,Y_2),\ldots$ is not
stationary, hence also not ergodic. It is, however, a regenerative
process. We shall also call the random variables $X_i$ and $Y_j$
\textit{related}, if they have the same ancestor. However, the
definition of the related sequences does not exclude the case, when
the functions $f$ and $g$ do not depend on $Z_i$ so that the
sequences $X$ and $Y$ are independent. Thus, in what follows, all
results for related sequences automatically hold for independent
sequences as well.

With this notation (recall (\ref{q})), $\q=1-\min_{a,b\in\mathcal
{A}}P(f(\xi,Z)=a|g(\eta,Z)=b).$ Note that
$P(f(\xi,Z)=a|g(\eta,Z)=b)=P(X_i=a|Y_j=b)$ given $X_i$ and $Y_j$ are
related.
\subsection{Limits and large deviation inequalities for related
sequences}\label{subsec:properties}
In this subsection, we consider the random variables $L_{n,an}$,
where $a>0$. By symmetry, for any $n$, the random variable
$L_{n,an}$ has the same law as $L_{an,n}$; moreover, the processes
$\{L_{an,n}\}$ and $\{L_{n,an}\}$ have the same distribution so that
in what follows, everything holds for $L_{an,n}$ as well.
\subsubsection*{The existence of $\gamma_{\tt{R}}(a)$}
At first, we
shall prove the
convergence (\ref{piirv-a}). As mentioned in the Section~\ref{sec:ind}, for independent sequences, this follows from
subadditive ergodic theorem. The same holds, if the sequences are
related, but no deletion occurs, that is, $p=1$. In the presence of
deletion, however, an additional argument is needed.
%
\begin{proposition}\label{prop:piirv}
Let $a>0$. Then there exists a constant $\gamma_{\tt{R}}(a)$ such
that (\ref{piirv-a}) holds.
\end{proposition}
%
%
\begin{pf}
At first note that without loss of generality, we may assume $a\leq
1$. Indeed, with $m:=\lfloor na \rfloor$,
\begin{eqnarray*}
L(X_1,\ldots,X_n;Y_1,\ldots,Y_{\lfloor na \rfloor})&=&
L(X_1,\ldots,X_{\lceil{m/
a}\rceil};Y_1,\ldots,Y_m)\\
&=&L(X_1,
\ldots,X_m;Y_1,\ldots,Y_{\lceil
{m/
a}\rceil}),
\end{eqnarray*}
where the last equality follows from the symmetry of
the model. Hence, the limit in (\ref{piirv-a}) exists if and only if
the limit of ${1\over m}
L(X_1,\ldots,X_m;Y_1,\ldots,Y_{\lceil{m/ a}\rceil})$ exists. The
latter is equivalent to the existence of limit ${1\over
m}L_{m,{m/ a}}$.
Hence, to the end of the proof, let $0<a\leq1$.

We consider the sequence of i.i.d. random vectors $U_1,U_2,\ldots$,
where
%
\begin{equation}
\label{vectors} U_i:=\bigl(f_i(Z_i),g_i(Z_i),D^x_i,D^y_i
\bigr).
\end{equation}
Let, for any
positive integer $m$, $n_x(m):=\sum_{i=1}^mD_i^x$ and
$n_y(m):=\sum_{i=1}^{\lfloor am \rfloor}D_i^y$.
Thus $X_1,\ldots,X_{n_x}$ and
$Y_1,\ldots,Y_{n_y}$ are both determined by i.i.d. random vectors
$U_1,\ldots, U_m$. Let
\[
L(U_1,\ldots,U_m):=L(X_1,
\ldots,X_{n_x};Y_1,\ldots ,Y_{n_y}).
\]
By subadditivity, there exists constant $\gamma_{\tt{U}}$ such that
%
\begin{equation}
\label{limesU} \lim_{m\to\infty}{L(U_1,\ldots,U_m)\over m}=
\gamma_{\tt{U}}, \qquad \mbox{a.s. and in }L_1.
\end{equation}
Let $\underline{n}(m):=n_x(m)\wedge{n_y(m)\over a}$ and
$\overline{n}(m):=n_x(m)\vee{n_y(m)\over a}$. Thus,
%
\begin{equation}
\label{limsup} {\underline{n}\over m}{L(X_1,\ldots,X_{\underline{n}};Y_1,\ldots,
Y_{\lfloor a\underline{n} \rfloor})\over\underline{n}}\leq
{L(U_1,\ldots,U_m)\over m}\leq{\overline{n}\over m}{L(X_1,\ldots,
X_{\overline{n}}, Y_1,\ldots, Y_{\lceil\overline{n}a\rceil})\over
\overline{n}}.
\end{equation}
%
By SLLN, ${\overline{n}(m)\over m}\to p$, a.s. and
${\underline{n}(m)\over m}\to p$, a.s. Since
\begin{eqnarray*}
\limsup_n {L_{n,an}\over n}=\limsup
_m {L_{\underline{n}(m),a\underline{n}(m)}\over\underline{n}(m)},\qquad  \liminf_n
{ L_{n,an} \over n}=\liminf_m {L_{\overline{n}(m), a
\overline{n}(m)}\over\overline{n}(m)}
\end{eqnarray*}
%
and
\begin{eqnarray*}
\liminf_m {L_{\overline{n}(m), a \overline{n}(m)}\over
\overline{n}(m)}=\liminf
_m {1\over
\overline{n}(m)}L(X_1,
\ldots,X_{\overline{n}(m)};Y_1,\ldots, Y_{\lceil a \overline{n}(m) \rceil}),
\end{eqnarray*}
from (\ref{limsup}), it follows
\begin{eqnarray*}
&&\limsup_n {1\over n}L(X_1,
\ldots,X_n;Y_1,\ldots, Y_{\lfloor a n
\rfloor})p\\
&&\quad \leq
\gamma_{\tt{U}}\leq\liminf_n {1\over
n}L(X_1,
\ldots,X_n;Y_1,\ldots, Y_{\lfloor a n \rfloor})p,\qquad  \mbox{a.s.}
\end{eqnarray*}
This is the a.s. convergence in (\ref{piirv-a}) with
$\gamma_{\tt{R}}(a)={\gamma_{\tt{U}}\over p}.$ The
convergence in
$L_1$ follows by dominated convergence theorem.
\end{pf}
%
\subsubsection*{Large deviation inequalities}
Next, we prove a large deviation lemma for related sequences.
%
\begin{lemma} Assume $X$ and $Y$ are related.
Then, for every $\Delta>0$ and $0<a\leq1$,
%
\begin{equation}
\label{LDV} P\bigl(|L_{n,an}-EL_{n,an}|\geq n\Delta\bigr)\leq4\exp
\biggl[-{p\over8}\Delta^2 a n\biggr].
\end{equation}
\end{lemma}
\begin{pf} As we saw in Section~\ref{sec:ind}, for independent
sequence, this type of inequality (\ref{McD}) trivially follows from
McDiarmid inequality. In the present case, we have
to add an extra control over the deletion process.

Fix positive integer $m$ and consider the vectors $U_1,\ldots,U_m$
defined in (\ref{vectors}). Recall $n_x(m)$ and $n_y(m)$. Fix $n$
and let
\[
\I_m:=L(X_1,\ldots, X_{n\wedge n_x};Y_1,
\ldots, Y_{\lfloor an
\rfloor\wedge n_y}).
\]
Note that $\I_m$ is a function of $5m$ independent random variables:
\[
\I_m=\I_m\bigl(Z_1,\ldots,Z_m,
\xi_1,\ldots,\xi_m,\eta_1,\ldots,\eta
_m,D_1^x,\ldots,D_m^x,D_1^y,
\ldots,D_m^y\bigr).
\]
Changing $Z_i$ (given all other variables are fixed) corresponds to
possible change of an element of $X$ and an element of $Y$. A
change of one element of $X$ (or $Y$) causes the change of $\I_m$ at
most by~1. Hence, the maximum change of $\I_m$ induced by changing of
$Z_i$ (given all other variables are fixed) is~2. Similarly, the
maximum change of $\I_m$ due to the change of $\xi_i$ or $\eta_i$
(given all other variables are fixed) is 1. Changing $D^x_i$ from 1
to 0 corresponds to removing one element of $X$-side and, in the
case $n_x>n$ adding one more $X$ to the end. Changing $D^x_i$ from 0
to 1 corresponds to adding one element to $X$-side and, perhaps,
removing the last $X$ (when $n_x\geq n$). This, again, changes the
value of $\I_m$ at most by 1. Any change of $\eta_i$ has the same
effect. Denoting by $r_i$, $i=1,\ldots, 5m$ the maximum change of
$\I_m$ induced by the $i$th variable, we have that $r_i=2$ if
$i=1,\ldots, m$ and $r_i=1$ for $i=m+1,\ldots, 5m$ so that
$\sum_{i=1}^{5m}r_i^2=8m$. Therefore, by McDiarmid inequality,
%
\begin{equation}
\label{cond} P \bigl(|\I_m-E \I_m|\geq m\Delta \bigr)\leq2\exp
\biggl[-{\Delta
^2\over
4}m\biggr].
\end{equation}
Let $E_m$ be the event that $n_x\geq n$ and $n_y\geq an$. Formally,
$E_m:=E_y(m)\cap E_x(m)$, where $E_x(m):=\{\sum_{i=1}^{m}D^x_i\geq
n\}$ and $E_y(m):=\{\sum_{i=1}^{\lfloor am \rfloor}D^y_i \geq an
\}$. When $E_m$ holds, then $\I_m=L_{n,an}$, so that
\[
\bigl\{|L_{n,an}-EL_{n,an}| < m\Delta \bigr\}\supset\bigl\{|\I_m-E
\I_m| < m\Delta\bigr\}\cap E_m
\]
and
%
\begin{equation}
\label{LI} P \bigl(|L_{n,an}-EL_{n,an}| \geq m\Delta \bigr)\leq P \bigl(|
\I_m-E\I_m| \geq m\Delta \bigr)+P\bigl(E^c_m
\bigr).
\end{equation}
Take $m={2\over p}n$. Then (\ref{cond}) is
%
\begin{equation}
\label{cond2} P \biggl(|\I_m-E\I_m| \geq
{2\over p}n\Delta \biggr)\leq 2\exp\biggl[-{\Delta^2\over2p}n
\biggr].
\end{equation}
To estimate $P(E^c_m)\leq P(E_x^c)+P(E_y^c)$, use Hoeffding
inequality (with $m={2n\over p}$)
\begin{eqnarray*}
P\bigl(E^c_y\bigr)&=&P \Biggl(\sum
_{i=1}^{am}D^y_i<an
\Biggr)=P \Biggl(\sum_{i=1}^{am}D^y_i-amp<an-amp
\Biggr)
\\
&\leq& P \Biggl(\sum_{i=1}^{am}
D^y_i-amp<-{amp\over2}\Biggr)\leq \exp
\biggl[-{p^2\over2}am\biggr]=\exp[-pan],
\\
P\bigl(E^c_x\bigr)&=&P \Biggl(\sum
_{i=1}^{m}D^x_i<n \Biggr)
\leq\exp[-pn]\leq\exp[-pan].
\end{eqnarray*}
Thus, with $m={2n\over
p}$, $P(E_m^c)\leq2\exp[-pan]$ and plugging it together with
(\ref{cond2}) into (\ref{LI}) entails
%
\begin{equation}
\label{LDP1} P \biggl(|L_{n,an}-EL_{n,an}| \geq
{2n\over p}\Delta \biggr)\leq 2\exp\biggl[-{\Delta^2\over2p}n
\biggr] +2\exp[-pan].
\end{equation}
Take $\Delta'={2\Delta\over p}$. Then (\ref{LDP1}) is
\begin{eqnarray*}
P \bigl(|L_{n,an}-EL_{n,an}|\geq\Delta' n \bigr)
\leq2 \exp\biggl[-{(\Delta')^2p\over8}n\biggr]+2\exp[-pan].
\end{eqnarray*}
If
$\Delta'\leq1$, then $2\exp[-apn]\leq2\exp[-(\Delta')^2 apn]$,
implying that the right-hand side is bounded by $4\exp[-{(\Delta
')^2\over
8} apn]$. This proves (\ref{LDV}) for $\Delta\leq1$. Since
$L_{n,an}\leq n$, for $\Delta>1$, (\ref{LDV}) trivially
holds.
\end{pf}
%

The following corollary states an inequality similar to that of
(\ref{LDV2ind}) for related sequences.
%
\begin{corollary}\label{cor}
Assume $X$ and $Y$ are related, $0<a\leq1$. Then, for every
$\Delta>0$ there exists $n_o(\Delta, a)$ big enough so that
%
\begin{equation}
\label{LDV2} P\bigl(\bigl|L_{n,an}-\gamma_{\tt{R}}(a)n\bigr|\geq n\Delta
\bigr)\leq4\exp\biggl[-{p\over
32}\Delta^2 an\biggr],\qquad
n>n_o.
\end{equation}
\end{corollary}
\begin{pf} Let $n$ be so big that $|{EL_{n,an}/ n}-\gamma_{\tt
{R}}(a)|<{\Delta/ 2}.$ Then $|EL_{n,an}-\gamma_{\tt{R}}(a)n|\leq
(\Delta/ 2)n$ and
\begin{eqnarray*}
P \bigl(\bigl|L_{n,an}-\gamma_{\tt{R}}(a)n\bigr|\geq n\Delta \bigr)&\leq& P
\bigl(|L_{n,an}-EL_{n,an}|+\bigl|EL_{n,an}-
\gamma_{\tt{R}}(a)n\bigr|\geq n\Delta \bigr)
\\
&\leq &P \biggl(|L_{n,an}-EL_{n,an}|\geq n{\Delta\over2}
\biggr)\leq 4\exp\biggl[-{p\over32}a\Delta^2 n\biggr],
\end{eqnarray*}
where the last inequality follows from (\ref{LDV})
\end{pf}
%
\section{Proofs of main results for related sequences}\label{sec:proofs}
\subsection{Every highest alignment contains a related
pair}\label{subsec:corner}
\subsubsection{The key lemma}
The following lemma is the cornerstone of what follows.
%
\begin{lemma}\label{related-pair} Assume that $X=X_1\ldots X_n$ and
$Y=Y_1\ldots Y_n$ are related
and satisfy (\ref{maincondition1}). Then there exists a constant
$b_2>0$ such that for every $n$ big enough,
%
\begin{equation}
\label{eq:lemma} P(\mbox{highest alignment of } X \mbox{ and } Y\mbox{ alignes no
related letters})\leq \mathrm{e}^{-nb_2}.
\end{equation}
\end{lemma}
%
%
\begin{pf}
Let $v\in V_k$ be an alignment. Let $I=I(v)=\{i_1,\ldots,i_k\}$ and
let $j_t:=v(i_t)$. Hence $X_{i_t}=Y_{j_t}$, for every
$t=1,\ldots,k$. We denote by $J$ the set $\{j_1,\ldots, j_k\}$.

We are bounding the probability that $v$ is the highest optimal
alignment of $X$ and $Y$ and that the random variables $X_{i_t}$ and
$Y_{j_t}$ are not related for every $t=1,\ldots, k$. Let us
introduce some notations and events. Let, for every $j=j_1+1,\ldots,
n$, $b(j)$ be the last element of $J$ strictly smaller than $j$.
Formally, $b(j):=\max\{j_t: j_t<j\}$. Similarly, for every
$i=1,\ldots, i_k-1$, let $c(i)$ be the first element of $I$ strictly
larger than $i$. Formally, $c(i):=\min\{i_t: i_t>i\}$. Also denote
\[
a^x= \bigl(a^x(i_1),\ldots,a^x(i_k)
\bigr),\qquad  a^y= \bigl(a^y(j_1),
\ldots,a^y(j_k) \bigr)
\]
and let $a^x\ne a^y$ be
$a^x(j_t)\ne a^y(j_t)$ for every $t=1,\ldots,k$. We now define the
following events
\begin{eqnarray*}
A(v)&:=&\{X_{i_1}=Y_{j_1},\ldots,X_{i_k}=Y_{j_k}
\},
\\
B(v)&:=&\bigl\{Y_j\ne Y_{b(j)}, j\in\{j_1+1,
\ldots, n\}\setminus J\bigr\},
\\
C(v)&:=&\bigl\{X_i\ne X_{c(i)}, i\in\{1,\ldots,
i_k-1\}\setminus I\bigr\},
\\
D(v)&:=&\bigl\{a^x\ne a^y\bigr\},\qquad  E(v):=A(v)\cap B(v)
\cap C(v) \cap D(v).
\end{eqnarray*}
By Corollary~\ref{corI}, it holds
\[
\{ \mbox{$v$ is the highest optimal alignment and no aligned pair of $v$ is
related}\}\subset E(v).
\]
Note that the vectors $a^x$
and $a^y$ depend on the deletion processes $D^x$ and $D^y$, only.
Thus, given $a^x$ and $a^y$, the events $A(v)$, $B(v)$ and $C(v)$
depend on the ancestor process $Z$ and on the random mappings $g$
and $f$, only. In particular, given $a^x$ and $a^y$, the dependence
structure (related pairs) is fixed as well. In the following, we
shall consider the case $a^x\ne a^y$. This means that there exists
no $t=1,\ldots, k$ such that $X_{i_t}$ is related to $Y_{j_t}$.

We shall bound the probability $P(E(v)|a^x,a^y)$.

At first, let us bound the probability $P(A(v)|a^x,a^y)$, $a^x\ne
a^y$. Thus, in what follows, we assume $a^x$ and $a^y$ satisfying
$a^x\ne a^y$ are fixed. For any two indexes $s,t\in
\{1,\ldots,k\}$, let $s \leftrightarrow t$ denote that either
$X_{i_t}$ and $Y_{j_s}$ are related (i.e., they have the same
ancestor) or $X_{i_s}$ and $Y_{j_t}$ are related. We call a subset
$G=\{t_1,\ldots,t_l\}\subset\{1,\ldots, k\}$ a \textit{dependence
group}, if:
\begin{enumerate}[2.]
\item$t_i \leftrightarrow t_{i+1}$ for every $i=1,\ldots,l-1$;
\item there is no index in $\{1,\ldots, k\}\setminus G$ that is related
to $t_1$ or $t_l$.
\end{enumerate}
Note that a group with $|G|$ elements contains $|G|-1$ related
pairs. Let $\{t_1,\ldots,t_l\}$ be a dependence group. Without loss
of generality, assume that $X_{i_{t_1}}$ is related to $Y_{j_{t_2}}$.
Then $X_{i_{t_2}}$ is related to $Y_{j_{t_3}}$ and so on. In
particular, $X_{i_{t_k}}$ is independent of $Y_{j_{t_l}}$, $l\leq
k$. Recall the definition of $p_o$, $\overline{p}$ $q$ and $\q$ from
(\ref{q-def}) and (\ref{q}). Hence,
\begin{eqnarray*}
&&P(X_{i_t}=Y_{j_t}; t\in G)\\
&&\quad =P(X_{i_{t_1}}=Y_{j_{t_1}})
\prod_{k=2}^l P(X_{i_{t_k}}=Y_{j_{t_k}}|
X_{i_{t_1}}=Y_{j_{t_1}},\ldots,X_{i_{t_{k-1}}}=Y_{j_{t_{k-1}}})
\\
&&\quad =p_o\prod_{k=2}^l \biggl(
\sum_{a\in\mathcal
{A}}P(X_{i_{t_k}}=a)P(Y_{j_{t_k}}=a|X_{i_{t_1}}=Y_{j_{t_1}},
\ldots ,X_{i_{t_{k-1}}}=Y_{j_{t_{k-1}}}) \biggr)\leq p_o({
\overline p})^{l-1}.
\end{eqnarray*}
By 2., the random variables $\{X_{i_t},Y_{j_t}: t\in G\}$ are all
independent of the random variables $\{X_{i_t},Y_{j_t}: t\in
\{1,\ldots, k\}\setminus G\}$. Let $G_1,\ldots, G_u$ be all
dependence groups. Let $G=G_1\cup\cdots\cup G_u$. Thus, $r:=|G|-u$
is the number of related pairs amongst $X_{i_1},\ldots, X_{i_k}$ and
$Y_{j_1},\ldots, Y_{j_k}$. By independence of the groups,
%
\begin{eqnarray}
\label{koer1} P\bigl(A(v)|a^x,a^y\bigr)&=&\prod
_{s=1}^uP(X_{i_t}=Y_{j_t}; t
\in G_s)\prod_{t\notin G} P(X_{i_t}=Y_{j_t})\nonumber\\[-8pt]\\[-8pt]
&\leq& p_o^u({\overline p})^{|G|-u}p_o^{k-|G|}=p_o^{k-r}({
\overline p})^r.\nonumber
\end{eqnarray}
Let us
now bound the probability $P(B(v)|A(v),a^x,a^y)$, where, as
previously, $a^x \ne a^y$. Recall the sets $I$ and $J$. Let
\[
I^c:=\{1,\ldots, i_k-1\}\setminus I,\qquad
J^c:=\{j_1+1,\ldots, n\} \setminus J.
\]
%
Let $J^c_1$ be the set of indexes in $J^c$ with the property that
the corresponding $Y$-s are related to an element in $X_{i}, i\in
I$. Formally, $j\in J^c_1$ if and only if there exists a $i\in I$ so
that $Y_j$ is related to $X_i$. Let $J^c_2=J^c\setminus J^c_1$. It
means, if $j\in J_2^c$, then $Y_j$ is either related to an $X_i$
with $i\notin I$ or not related to any other random variable at
all. In particular, the random variables $\{Y_j: j\in J_2^c\}$ are
independent of the event $A(v)$. Since $Y_1,Y_2,\ldots$ are
independent, we obtain (let us omit the fixed $a^x$ and $a^y$ from
the notations)
\begin{eqnarray*}
P\bigl(B(v)|A(v)\bigr)&=&P\bigl(Y_j\ne Y_{b(j)}, j\in
J_2^c\cup J^c_1|A(v)\bigr)
\\
&=&P\bigl(Y_j\ne Y_{b(j)}, j\in J_2^c
\bigr)P\bigl(Y_j\ne Y_{b(j)}, j\in J^c_1|A(v)
\bigr).
\end{eqnarray*}
Clearly,
%
\begin{equation}
\label{koer2} P\bigl(Y_j\ne Y_{b(j)}, j\in
J_2^c\bigr)\leq q^{|J_2^c|}.
\end{equation}
Let us estimate
$P(Y_j\ne Y_{b(j)}, j\in J^c_1|A(v))$. Note
\begin{eqnarray*}
&&P\bigl(Y_j\ne Y_{b(j)}, j\in J^c_1|A(v)
\bigr)\\
&&\quad =\sum_{(y_1,\ldots, y_k)\in
\mathcal{A}^k}P\bigl(Y_j\ne
Y_{b(j)}, j\in J^c_1| Y_{j_t}=X_{i_t}=y_{t},
\forall t\bigr)P \bigl(X_{i_t}=y_{t}, \forall t| A(v)
\bigr).
\end{eqnarray*}
Given $(y_1,\ldots,y_k)$, let $y_{b(j)}$ be the value
of $Y_{b(j)}$. Let us estimate
\begin{eqnarray*}
P\bigl(Y_j\ne Y_{b(j)}, j\in J^c_1|
Y_{j_t}=X_{i_t}=y_{t}, \forall t\bigr)&=&P
\bigl(Y_j\ne y_{b(j)}, j\in J^c_1|
Y_{j_t}=X_{i_t}=y_t, \forall t\bigr)
\\
&=&P\bigl(Y_j\ne y_{b(j)}, j\in J^c_1|
X_{i_t}=y_t, \forall t\bigr)
\\
&=&{P(Y_j\ne y_{b(j)}, j\in J^c_1;X_{i_t}=y_t, \forall t)\over\prod_{t=1}^k
P(X_{i_t}=y_t)}.
\end{eqnarray*}
The last two equalities follow from the fact that $Y_1,Y_2,\ldots$
are independent and $X_1,X_2,\ldots$ are independent. When $j\in
J_1^c$, then $Y_j$ is related to a $X_{i_t}$. Denote
$J_1^c:=\{j^1,\ldots,j^s\}$. Clearly,
%
\begin{equation}
\label{s-ineq} s:=|J_1^c| \leq|J^c|\wedge k
= (n-j_1+1-k)\wedge k.
\end{equation}
Without loss of generality, assume that the random variables in
$J_1^c$ are related to the $X_{i_1},\ldots,X_{i_s}$. Then the pairs
of related random variables $(Y_{j^1},X_{t_1}),\ldots,
(Y_{j^s},X_{t_s})$ are independent so that
\begin{eqnarray*}
{P(Y_j\ne y_{b(j)}, j\in J^c_1;X_{i_t}=y_t, \forall t)\over
\prod_{t=1}^k P(X_{i_t}=y_t)}&=&{\prod_{t=1}^s P(Y_{j^t}\ne
y_{b(j^t)},X_{i_t}=y_t)\over\prod_{t=1}^s P(X_{i_t}=y_t)} \\
&=&\prod
_{t=1}^s P(Y_{j^t}\ne
y_{b(j^t)}|X_{i_t}=y_t)\leq (\bar{q})^s.
\end{eqnarray*}
Therefore,
%
\begin{equation}
\label{barq} P\bigl(Y_j\ne Y_{b(j)}; j\in
J^c_1|A(v)\bigr)\leq(\bar{q})^s.
\end{equation}
%
By entirely similar argument, we estimate $P(C(v)|A(v)\cap B(v))$.
Indeed, given $(y_1,\ldots,\break y_k)\in\mathcal{A}^k$,
\begin{eqnarray*}
&&P\bigl(C(v)|A(v)\cap B(v),Y_{j_t}=y_t, \forall t\bigr)\\
&&\quad =
P\bigl(X_i\ne X_{c(i)}, i=I^c|
Y_{j_t}=X_{i_t}=y_t, \forall t;
Y_j\ne Y_{b(j)}, j=J^c\bigr)
\\
&&\quad = P\bigl(X_i\ne a_{c(i)}, i=I^c|
Y_{j_t}=y_t, \forall t; Y_j\ne
y_{b(j)}, j=J^c\bigr).
\end{eqnarray*}
Every $X_i$ is related to at most one $Y_j$. Let $I_0^c,
I_1^c,I_2^c$ be mutually exclusive set of indexes so that:
\begin{itemize}[$\bullet$]
\item If $i\in I_0^c$, then $X_i$ is not related to any $Y_j$ from
$J\cup
J^c$.
\item If $i\in I_1^c$, then $X_i$ is related to a $Y_j$ so
that $j\in J$. Let $t(i)\in\{1,\ldots,k\}$ be the
corresponding index.
\item If $i\in I_2^c$, then $X_i$ is related to a $Y_j$ so
that $j\in J^c$. Let $j_r(i)\in J$ be the corresponding index.
\end{itemize}
Then, just like previously, using the independence of related pairs,
we obtain
\begin{eqnarray*}
&&P\bigl(X_i\ne y_{c(i)}, i=I^c|
Y_{j_t}=y_t, \forall t; Y_j\ne
y_{b(j)}, j=J^c\bigr)
\\
&&\quad =\prod_{i\in I_0^c} P (X_i\ne
y_{c(i)} ) \prod_{i\in
I_1^c}P (X_i
\ne y_{c(i)}|Y_{t(i)}=y_{t(i)} )\prod
_{i\in
I_2^c}P (X_i\ne y_{c(i)}|Y_{j_r(i)}
\ne y_{b(j_r(i))} )
\\
&&\quad \leq q^{|I_0^c|}({\bar q})^{|I_1^c|+|I_2^c|}\leq({\bar q})^{|I^c|},
\end{eqnarray*}
where the second last inequality follows from the fact that given
$X_i$ and $Y_j$ are related, for any $a,b\in\mathcal{A}$
\begin{eqnarray*}
P(X_i\ne a|Y_j\ne b)&=&\sum
_{c\ne b}P(X_i\ne a|Y_j=c)P(Y_j=c|Y_j
\ne b)
\\
&=&\sum_{c\ne b}\bigl(1-P(X_i=a|Y_j=c)
\bigr)P(Y_j=c|Y_j\ne b)\leq\bar{q}
\end{eqnarray*}
and the last inequality follows from the fact that $q\leq{\bar q}$.
Therefore,
%
\begin{equation}
\label{barq2} P\bigl(C(v)|A(v)\cap B(v)\bigr)\leq(\q)^{|I^c|}=(
\q)^{i_k-k}.
\end{equation}
%
By (\ref{koer1}), (\ref{koer2}), (\ref{barq}), (\ref{barq2}) with
$\rho=(p_o \q) /(\p q)$ and $r+s\leq k$, we have
\begin{eqnarray*}
P\bigl(E(v)|a^x,a^y\bigr)&\leq& p_o^{k-r}(
\p)^r q^{n-j_1+1-k-s} (\q)^{s+i_k-k}=p_o^k
\biggl({\p\over p_o} \biggr)^r \biggl(
{\q\over q
} \biggr)^s q^{n-j_1+1-k}{(
\q)}^{i_k-k}
\\
&\leq& p_o^k \biggl({\p\over p_o}
\biggr)^{k-s} \biggl({\q\over q
} \biggr)^s
q^{n-j_1+1-k}{(\q)}^{i_k-k} \leq(\p)^k \rho^s
q^{n-j_1+1-k}({\overline q})^{i_k-k}.
\end{eqnarray*}
By (\ref{s-ineq}), it holds $0\leq s\leq k \wedge(n-j_1+1-k)\leq k
\wedge(n-k),$ so that
\[
\max_s \rho^s \leq\lleft\{ %
\begin{array} {l@{\qquad}l} \rho^{k \wedge(n-k)}, & \mbox{if $\rho\geq1$};
\\
1, & \mbox{if $\rho< 1$}. \end{array} %
\rright.
\]
Hence,
%
\begin{eqnarray}
\label{viiskaks} P\bigl(E(v)\bigr)&\leq&\sum_{a^x,a^y: a^x\ne
a^y}P
\bigl(E(v)|a^x,a^y\bigr)P\bigl(D^x=a^x,D^y=a^y
\bigr)\nonumber\\[-8pt]\\[-8pt]
&\leq& (\p)^k (\rho\vee 1)^{k\wedge(n-k)} (q\q)^{n-k}
q^{1-j_1}({\overline q})^{i_k-n}.\nonumber
\end{eqnarray}
%
Recall that (\ref{maincondition1}) is
\[
\gamma_{\tt{R}}\log_2 \p+(1-\gamma_{\tt{R}})
\log_2 (q \q)+ \bigl((1-\gamma_{\tt
{R}})\wedge
\gamma_{\tt{R}}\bigr)\log_2 (\rho\vee1 )+ 2h(\gamma
_{\tt{R}})<0.
\]
%
When this holds, then it is possible to find $\Delta>0$ so small
that
\begin{eqnarray*}
-b&:=&(\gamma_{\tt{R}}-\Delta)\log_2 \p+(1-
\gamma_{\tt
{R}}-\Delta)\log_2 (q \q)\\
&&{}+ \bigl((1-
\gamma_{\tt{R}}+\Delta )\wedge (\gamma_{\tt{R}}+\Delta)\bigr)
\log_2 (\rho\vee1 )-\Delta \log_2(q \q)+ 2H(
\gamma_{\tt{R}},\Delta)<0.
\end{eqnarray*}
Let
\[
E_{\Delta}:=\bigl\{|L_n-n\gamma_{\tt{R}}|<n\Delta\bigr\}.
\]
When $E_{\Delta}$ holds, then all optimal alignments belong to the
set $W_n:=W_n(\gamma_{\tt{R}},\Delta)$. For every $v\in W_n$, with
$|v|=k$, it holds
%
\begin{equation}
\label{condk} n(\gamma_{\tt{R}}-\Delta)\leq k \leq n(
\gamma_{\tt
{R}}+\Delta).
\end{equation}
Let, for every $v$, $s(v)=j_1-1$ and $t(v)=n-i_k$. Let
\[
U_n(\gamma_{\tt{R}},\Delta):=\bigl\{v\in W_n:
s(v)\leq\Delta n, t(v)\leq\Delta n\bigr\}.
\]
Using
these two inequalities together with (\ref{condk}), we have that for
every $v\in U_n$,
\begin{eqnarray*}
\log_2 P\bigl(E(v)\bigr)&\leq& n \biggl[(\gamma_{\tt{R}}-
\Delta)\log_2 \p+(1-\gamma_{\tt{R}}-\Delta)\log_2
(q \q)
\\
&&{}\hphantom{n \biggl[} +\bigl((1-\gamma_{\tt{R}})\wedge\gamma_{\tt{R}}+\Delta\bigr)\log
_2(\rho\vee 1)-{s(v)\over n}\log_2 q-
{t(v)\over n}\log_2 \q \biggr]
\\
&\leq &\bigl(-b-2H(\gamma_R,\Delta) \bigr)n.
\end{eqnarray*}
Let
\[
E:=\{\exists\mbox{ highest alignment of} X \mbox{ and } Y\mbox{ alignes no
related letters}\}.
\]
Recall that
$S=j_1^h-1,T=n-i_k^h$ and, by (\ref{boundVn}), it holds
\[
|U_n|\leq\bigl|W_n(\gamma_{\tt{R}},\Delta)\bigr|\leq 2
\Delta n2^{2H(\gamma_{\tt{R}},\Delta)}.
\]
Then by
Corollary~\ref{cor:st} and Corollary~\ref{cor}, for $n$ big enough
\begin{eqnarray*}
P(E)&\leq&\sum_{v\in U_n}P\bigl(E(v)\bigr)+P(S>\Delta
n)+P(T>\Delta n)+P\bigl(E^c_{\Delta}\bigr)
\\
&\leq&|U_n|2^{(-b-2H(\gamma_{\tt{R}},\Delta))}+P(S>\Delta n)+P(T>\Delta n)+P
\bigl(E^c_{\Delta}\bigr)
\\
&\leq&2\Delta n 2^{-bn}+4\exp\biggl[-{\Delta^2\over32}n\biggr]+2
\exp\bigl[-d(\Delta)n\bigr].
\end{eqnarray*}
Hence, for big $n$, the inequality \eqref{eq:lemma}
holds.
\end{pf}
%
\paragraph*{Sequences with unequal lengths}
In the previous lemma, $X$ and $Y$ were of the same length, $n$.
This lemma can be generalized for the case $X$ and $Y$ are of
different length, provided that the difference is not too big. Let
$X^n:=X_1\ldots X_n$, $Y^m:=Y_1\ldots Y_m$. Without loss of
generality, let us assume $m\geq n$. We know that if
(\ref{maincondition1}) holds, then there exists $\Delta>0$ so small
that
\begin{eqnarray}
\label{mainconditiondelta} &&(\gamma_{\tt{R}}-\Delta)\log_2 \p+(1-
\gamma_{\tt
{R}}-2\Delta)\log_2 (q \q)+ \bigl((1-
\gamma_{\tt{R}})\wedge\gamma_{\tt{R}}+2\Delta \bigr)
\log_2 (\rho\vee1 )\nonumber\\[-8pt]\\[-8pt]
&&\quad {}-2\Delta\log_2(q \q)+ 2H(
\gamma_{\tt{R}},2\Delta)<0.\nonumber
\end{eqnarray}
The restriction for $m$ is: $m\leq(1+\Delta)n$.
%
\begin{lemma}\label{related-pair2}
Let $n\leq m\leq(1+\Delta)n$, where $\Delta>0$ satisfies
(\ref{mainconditiondelta}). Assume that $X^n$ and $Y^m$ are related.
Then there exists a constant $b_3(\Delta)>0$ such that for every
$n>n_o$,
\[
P\bigl( \mbox{the highest alignment of } X^n \mbox{ and }
Y^m\mbox{ aligns no related letters}\bigr)\leq \mathrm{e}^{-nb_3}.
\]
\end{lemma}
\begin{pf} The proof follows the one of Lemma~\ref{related-pair};
$\Delta$ is now taken from the assumptions, so it satisfies (\ref
{mainconditiondelta}). This
$\Delta$ defines the set $E_{\Delta}$ as in the previous lemma.
However, by definition, $L_n$ is the length of the LCS between $X^n$
and $Y^n$, whilst in the present case we are dealing with the LCS
between $X^n$ and $Y^m$. Clearly $L_n\leq L_{n,m}\leq L_n+n\Delta$.
Hence, if $E_{\Delta}$ holds, then
\[
\gamma_{\tt{R}}-\Delta\leq{L_n\over n}\leq
{L_{n,m}\over n}\leq{L_n\over n}+\Delta\leq
\gamma_{\tt
{R}}+2\Delta,
\]
that is, all optimal alignments belong to the set
$W_{n,m}(\gamma_{\tt{R}},\Delta)$. The set $U_{n,m}$ is defined as
follows
\[
U_{n,m}(\gamma_{\tt{R}},\Delta):=\bigl\{v\in W_{n,m}(
\gamma_{\tt
{R}},\Delta): s(v)\leq2\Delta n, t(v)\leq 2\Delta n\bigr\}.
\]
The upper bound (\ref{viiskaks}) holds with $n$
replaced by $m$:\vspace*{1pt}
\[
P\bigl(E(v)\bigr)\leq(\p)^k (\rho\vee1)^{k\wedge(m-k)} (q
\q)^{m-k} q^{1-j_1}({\overline q})^{i_k-m}.
\]
%
Using the bounds $(\gamma_{\tt{R}}-\Delta)n\leq|u| \leq
(\gamma_{\tt{R}}+2\Delta)n$ and $n\leq m\leq n(1+\Delta)$, for every
$v\in U_{n,m}$, we obtain the following estimate
\begin{eqnarray*}
\log_2 P\bigl(E(v)\bigr)&\leq& n \biggl[(\gamma_{\tt{R}}-
\Delta)\log_2 \p+\bigl((1-\gamma_{\tt{R}})\wedge
\gamma_{\tt{R}}+2\Delta\bigr)\log _2(\rho\vee 1)
\\
&&\hphantom{n \biggl[}{}+(1-\gamma_{\tt{R}}-2\Delta)\log_2 (q \q)-
{s(v)\over n}\log q_2-{t(v)\over n}
\log_2 \q \biggr]\\
&\leq& - \bigl(b+2H(\gamma_{\tt
{R}},\Delta)
\bigr)n,
\end{eqnarray*}
where $b>0$ by the assumption (\ref{mainconditiondelta}) on
$\Delta$. The rest of the proof goes as the one of Lemma~\ref{related-pair} with $P(S>2\Delta n)$ and $P(T>2\Delta n)$
instead of $P(S>\Delta n)$ and $P(T>\Delta n)$ and $3\Delta
n2^{-bn}$ instead of $2\Delta n2^{-bn}$.
\end{pf}
%
\subsubsection{Applying Lemma \texorpdfstring{\protect\ref{related-pair2}}{5.2} repeatedly: The
$B$-events}
\paragraph*{Regenerativity} Let $\tau^x_0=\tau^y_0=0$ and let $\tau^x_k$
($\tau^y_k$), $k=1,2,\ldots$ be the indexes of the $k$th related
pair. So, $(X_{\tau^x_1},Y_{\tau^y_1})$ is the first related pair,
$(X_{\tau^x_2},Y_{\tau^y_2})$ is the second related pair and so on.
Let $a_0=0$ and $a_k$ be the common ancestor of the $k$th related
pair, that is,
\[
a_k=a^x\bigl(\tau^x_k
\bigr)=a^y\bigl(\tau^y_k\bigr).
\]
We shall use the
fact that the process $(X_1,Y_1),(X_2,Y_2),\ldots$ is regenerative
with respect to the times $(\tau_k^x,\tau_k^y)$, i.e.
%
\begin{equation}
\label{afterstopping} (X_{\tau_k^x+1},Y_{\tau_k^y+1}),(X_{\tau_k^x+2},Y_{\tau
_k^y+2}),
\ldots
\end{equation}
has the same law as $(X_1,Y_1),(X_2,Y_2),\ldots\,$. The $Z$-process
for (\ref{afterstopping}) is $Z_{a_k+1},Z_{a_k+2},\ldots\,$.
\paragraph*{Definition of $B$-events}
In what follows, let $\Delta>0$ and $0<A<\infty$. Denote $n':=A \ln
n$. We shall consider the following events:
\begin{eqnarray*}
B_k(\tilde{n},\tilde{m})&:=&\{\mbox{the highest alignment of }
X_{\tau_k^x+1},\ldots,X_{\tau_k^x+\tilde{n}} \mbox{ and } Y_{\tau_k^y+1},
\ldots,Y_{\tau_k^y+\tilde{m}}
\\
&&\hphantom{\{} \mbox{contains a related pair} \},
\\
B^1_k\bigl(n',\Delta\bigr)&:=&\bigcap
_{n' \leq\tilde{n}\leq\tilde{m}\leq
\tilde{n}(1+\Delta)}B_k(\tilde{n},\tilde{m} ),\qquad
B^2_k\bigl(n',\Delta\bigr):=\bigcap
_{n' \leq\tilde{m}\leq\tilde{n}\leq
\tilde{m}(1+\Delta)}B_k(\tilde{n},\tilde{m} ),
\\
B^{h}_k\bigl(n',\Delta
\bigr)&:=&B^1_k\bigl(n',\Delta\bigr)\cap
B^2_k\bigl(n',\Delta\bigr).
\end{eqnarray*}
Let $B^l_k(n',\Delta)$ be defined similarly, with ``lowest'' instead
of ``highest'' in the definition of $B_k(\tilde{n},\tilde{m})$.
Finally, let
\[
B\bigl(k,n',\Delta\bigr):=B^l_k
\bigl(n',\Delta\bigr)\cap B^h_k
\bigl(n',\Delta\bigr).
\]
%
Let $B_n(n',\Delta)$ be the event that for every $k$ that satisfies
$\max\{\tau_k^x,\tau_k^y\}\leq n$, $B(k,n',\Delta)$ holds. Formally,
\[
B_n\bigl(n',\Delta\bigr):=\bigcup
_{i=0}^{n} \biggl(\{K=i\}\cap \Biggl(\bigcap
_{k=0}^i B\bigl(k,n',\Delta\bigr) \Biggr)
\Biggr),
\]
where
%
\begin{equation}
\label{K} K:=\arg\max_{k=0,1,\ldots} \bigl\{\max\bigl\{
\tau_k^x,\tau_k^y\bigr\}\leq n
\bigr\}.
\end{equation}
%
\paragraph*{The bound on $P(B_n(n',\Delta))$ for small $\Delta$}
We aim to bound $P(B_n(n',\Delta))$ from below. We use the
regenerativity described above: for every $k$, the event
$B_k(\tilde{n},\tilde{m})$ has the same probability as
$B_0(\tilde{n},\tilde{m})$ so that for every $k$, Lemma~\ref{related-pair2} applies:
\[
P\bigl(B_k(\tilde{n},\tilde{m})\bigr)\geq 1-\exp[-
\tilde{n}b_3],
\]
provided $\tilde{n}(1+\Delta)\geq
\tilde{m}\geq\tilde{n}\geq n_o$ and $\Delta>0$ is small enough to
satisfy the assumptions (\ref{mainconditiondelta}). Thus, for small
enough $\Delta$ and big enough $n'$, we have
%
\begin{eqnarray}
\label{b3} P \bigl(B^2_k\bigl(n'\bigr)
\bigr)=P \bigl(B^1_k\bigl(n'\bigr) \bigr)
\geq 1-\sum_{\tilde{n}\geq n'}\sum_{\tilde{n}(1+\Delta)\geq\tilde
{m}\geq
\tilde{n}}\mathrm{e}^{-b_3\tilde{n}}=1-
\sum_{\tilde{n}\geq n'} (\Delta \tilde{n}+1)\mathrm{e}^{-b_3\tilde{n}}.
\end{eqnarray}
Clearly, for every
$0<b_4<b_3$, for every $n$ big enough,
$(n+\Delta^{-1})\mathrm{e}^{-b_3n}\leq \mathrm{e}^{-b_4n}$ for every $n>n_1$. Let
$0<b_4<b_3$ and without loss of generality assume $n_o$ being so big
that for every $n>n_o$ the inequality above holds. Then (\ref{b3})
can be bounded as follows
%
\begin{eqnarray}
\label{b4} P \bigl(B^1_k\bigl(n',
\Delta\bigr) \bigr)\geq1- \Delta\sum_{\tilde{n}\geq n'}
\mathrm{e}^{-b_4\tilde{n}}\geq 1-{B\over4}\mathrm{e}^{- b_4 n'},\qquad  n'
\geq n_o,
\end{eqnarray}
where $B$ is a constant depending on $\Delta$. Hence,
\[
P \bigl(B\bigl(k,n',\Delta\bigr) \bigr)\geq1-B\mathrm{e}^{-b_4n'},\qquad
n'\geq n_o.
\]
Finally, since $\bigcap_{k=0}^n B(k,n')\subset B_n(n'),$ we have that
(recall $n'=A \ln n$)
%
\begin{eqnarray}
\label{B-est} P \bigl(B^c_n\bigl(n',
\Delta\bigr) \bigr)&\leq&(n+1)P\bigl(B_k^c
\bigl(n',\Delta\bigr)\bigr)\nonumber\\[-8pt]\\[-8pt]
&\leq& B(n+1)\exp\bigl[-b_4n'
\bigr]\leq 2Bn\exp\bigl[-b_4n'\bigr]=2Bn^{1-b_4A}.\nonumber
\end{eqnarray}
%
\subsection{The location of the related pairs}\label{subsec:proprel}
We consider the related sequences $X_1,X_2,\ldots$ and
$Y_1,Y_2,\ldots\,$. Recall the definition of $\tau^x_k$, $\tau^y_k$
and $a_k$. As previously, we take $n'=A\ln n$.
\subsubsection{The location of the first and last related pair: $G$-events}
\paragraph*{The location of last related pair: Definition of
$G_n(\Delta)$}
Let $i(n)$ and $j(n)$ be the biggest $\tau_k^x$ and $\tau_k^y$
before $n$, that is,
\[
i(n):=\max\bigl\{\tau^x_k: \tau^x_k
\leq n\bigr\}, \qquad j(n):=\max\bigl\{\tau^y_k:
\tau^y_k\leq n\bigr\}.
\]
Clearly $i(n)=n$ if and
only if the ancestor of $X_n$ is also an ancestor of a $Y_j$ that is,
$D^y_{a^x(n)}=1$. Similarly, $i(n)=u<n$ if and only if
\[
D^y_{a^x(u)}=1,\qquad  D^y_{a^x(u+1)}=\cdots=
D^y_{a^x(n)}=0.
\]
Since the process $D^y$ is independent of $D^x$ and, therefore, also
independent of the random variables $a^x(i)$, $i=1,2\ldots\,$, we have
that for every $u=1,2,\ldots,n$
\[
P \bigl(i(n)=u \bigr)=(1-p)^{n-u}p,\qquad  P \bigl(i(n)=0
\bigr)=(1-p)^n.
\]
Hence, for any $\Delta>0$,
%
\begin{eqnarray}
P \bigl(n-i(n)\geq\Delta n \bigr)&=&P \bigl(n-i(n)\geq\lceil\Delta n \rceil
\bigr)=(1-p)^{\lceil\Delta n \rceil}\nonumber\\[-8pt]\\[-8pt]
&\leq&(1-p)^{\Delta n}=\exp \bigl[\Delta n \ln(1-p)
\bigr].\nonumber
\end{eqnarray}
Hence, the probability that the last related $X_i$ before $n$ is
further that $\Delta n$ from $n$ is exponentially small in $n$. The
same obviously holds for $j(n)$ so that
%
\begin{equation}
P \bigl(n-i(n)<\Delta n, n-j(n)<\Delta n \bigr)\geq1-2\exp\bigl[\ln (1-p)\Delta
n\bigr].
\end{equation}
However, the event $\{(n-j(n))\vee(n-i(n))<\Delta n\}$ does not
necessarily imply that the last related pair, let that be $(i',j')$
is necessarily such that $\{(n-j')\vee(n-i')<\Delta n\}$. Indeed, if
$(i',j')$ is the last related pair, then either $i'=i(n)$ or
$j'=j(n)$ but the both inequalities need not hold simultaneously. We
shall now show that also the event $\{(n-j')\vee(n-i')\}\leq\Delta
n\}$ holds with great probability. Let us first define the last
related pair formally as follows
%
\begin{eqnarray}
\label{lastrelatedpair} &&i'(n):=\tau^x_{l(n)},\qquad
j'(n):=\tau^y_{l(n)},\nonumber \\[-8pt]\\[-8pt]
&&\quad  \mbox{where }l(n):=\max
\bigl\{l=0,1,2,\ldots: \tau_k^x\leq n,
\tau_k^y\leq n\bigr\}.\nonumber
\end{eqnarray}
Let $0<\Delta<1$, $r(n):=(1-{3\over4}\Delta){n\over p}$ and
consider the event
\[
G_n^x:= \Biggl\{n(1-\Delta)\leq\sum
_{j=1}^rD_j^x\leq n
\biggl(1-{\Delta\over
2}\biggr) \Biggr\},\qquad  G_n^y:=
\Biggl\{n(1-\Delta)\leq\sum_{j=1}^rD_j^y
\leq n\biggl(1-{\Delta\over2}\biggr) \Biggr\}.
\]
To simplify the calculations, let us
assume without the loss of generality that $r$ is an integer. Assume
that $G_n^x\cap G_n^y$ holds,
$i(n)> n(1-{\Delta\over2})$ and let $Y_j$ be related to
$X_{i(n)}$. Let $a$ be their common ancestor, that is,
$a:=a^y(j)=a^x(i(n))$. Since $i(n)>n(1-{\Delta\over2})$, the event
$G_n^x$ guarantees that $a>r$. Indeed, if $a\leq r$, then we reach
to the contradiction, since
\[
i(n)=\sum_{j=1}^aD_j^x
\leq\sum_{j=1}^r D_j^x
\leq n\biggl(1-{\Delta\over
2}\biggr).
\]
Thus $a>r$ and because of $G_n^y$, it holds
\[
j=\sum_{j=1}^aD_j^y
\geq\sum_{j=1}^rD_j^y
\geq n(1-\Delta).
\]
Hence, if $\tau^x_{l(n)}\geq\tau^y_{l(n)}$ (i.e., $i(n)=i'(n)\geq
j$), then we have that $\tau^y_{l(n)}=j'(n)\geq n(1-\Delta)$. The
roles of $X$ and $Y$ can be changed so that
%
\begin{eqnarray}
&&\biggl\{\bigl(n-j(n)\bigr)\vee\bigl(n-i(n)\bigr)<{\Delta\over2} n
\biggr\}\cap G_n^x\cap G_n^y\nonumber\\[-8pt]\\[-8pt]
&&\quad \subset \bigl\{ \bigl(n-j'(n)\bigr)\vee\bigl(n-i'(n)
\bigr) \leq\Delta n \bigr\} =:G_n(\Delta).\nonumber
\end{eqnarray}
%
When $G_n(\Delta)$ holds, then the last related pair, say $(i',j')$,
satisfies: $(i',j')\in[(1-\Delta)n,n]\times[(1-\Delta)n,n]$. In
2-dimensional representation, this means that the last related pair
is located in a square of size $\Delta n$ in the upper-right corner.
\paragraph*{The bound on $P(G_n(\Delta) )$}
By Hoeffding's inequality,
\begin{eqnarray*}
P \bigl(\bigl(G_n^x(\Delta)\bigr)^c
\bigr)&=&P \Biggl(\Biggl|\sum_{j=1}^r
D_j^x - rp\Biggr|>{\Delta\over4}n \Biggr)=P
\Biggl(\Biggl|\sum_{j=1}^r D_j^x
- rp\Biggr|>{p\Delta
\over4-3\Delta}r \Biggr)
\\
&\leq&2\exp\biggl[-2 \biggl({p\Delta\over4-3\Delta} \biggr)^2r
\biggr]=2\exp\biggl[-{p
\Delta^2\over2(4-3\Delta)}n\biggr].
\end{eqnarray*}
Therefore, for $\Delta$
small enough,
\begin{eqnarray*}
P\bigl(G^c_n(\Delta)\bigr)&\leq& P \bigl(
\bigl(G_n^x(\Delta)\bigr)^c \bigr)+P \bigl(
\bigl(G_n^y(\Delta)\bigr)^c \bigr)+P
\biggl(n-i(n)\geq{\Delta\over2}n\biggr)+P\biggl(n-j(n)\geq
{\Delta\over2}n\biggr)
\\
&\leq& 4\exp\biggl[-{p \Delta^2\over2(4-3\Delta)}n\biggr]+2\exp\biggl[
{\Delta\over
2}n \ln(1-p)\biggr]\\
&\leq&4\exp\biggl[-{p \Delta^2\over8}n
\biggr]+2\exp\biggl[{\Delta\over
2}n \ln(1-p)\biggr]
\\
&\leq&6\exp\biggl[-{p \Delta^2\over8}n\biggr].
\end{eqnarray*}
%
Finally, we shall apply the event $G_n(\Delta_n)$ with
$\Delta_n:={\sqrt{ 16 \ln n \over p n}}.$ Then
%
\begin{eqnarray}
\label{Ggn} P \bigl(G^c_n(\Delta_n)
\bigr)\leq6\exp\biggl[-{p \over8} \Delta_n^2
n\biggr]=6n^{-2}.
\end{eqnarray}
%
\subsubsection{The location of the rest of the related pairs: $F$-events}
Fix $\Delta>0$ and denote $\alpha:=1+{\Delta\over2}$,
$\beta:=1+\Delta$ and $l(n):={\alpha\over p}n$. Again, to simplify
the technicalities, let us assume that $l(n)$ is an integer.
\paragraph*{$F$-events: The definition}
At first, we consider the events
\[
F_n^x:=\Biggl\{n < \sum_{i=1}^l
D^x_i \leq\beta n\Biggr\},\qquad  F_n^y:=
\Biggl\{n < \sum_{i=1}^l
D^y_i \leq\beta n\Biggr\}, \qquad F_n:=F_n^x
\cap F_n^y.
\]
These events are similar to the events $G_n^x$ and $G_n^y$
defined in the previous section and we shall argue similarly.
Suppose $X_i$ and $Y_j$ are related and $i\leq n$. When $F^x_n$
holds, then the ancestor of $X_i$ is at most $l$, that is, $a^x(i)\leq
l$. Since $X_i$ and $Y_j$ are related, $a^x(i)=a^y(j)=:a$. If
$F^y_n$ holds, we have $\sum_{i=1}^aD^y_i\leq\sum_{i=1}^lD^y_i\leq
\beta n$, implying that $j\leq\beta n=(1+\Delta)n$. By symmetry,
the roles of $i$ and $j$ can be changed. Thus, when the event $F_n$
holds and $(i,j)$ is a related pair, then the following implication
holds true: if $\min\{i,j\}\leq n$, then $\max\{i,j\}\leq
(1+\Delta)n$.

We now consider more refined events
\[
F\bigl(k,n'\bigr):=\bigcap_{m\geq n'}
\Biggl\{a_k+m < \sum_{i=1}^{l(m)}
D^x_{a_k+i},\sum_{i=1}^{l(m)}
D^y_{a_k+i}\leq a_k+(1+\Delta)m \Biggr\},\qquad
k=0,1,2,\ldots.
\]
%
The event $F(k,n')$ states that for any other related pair
$X_{\tau^x_l}$, $Y_{\tau^y_l}$, $l>k$, the following holds: if
$\tau^x_l-\tau^x_k \leq n'$, then $\tau^y_l-\tau^y_k \leq
n'(1+\Delta)$. If $\tau^x_l-\tau^x_k=m>n'$, then $\tau^y_l-\tau^y_k
\leq m(1+\Delta)=(\tau^x_l-\tau^x_k)(1+\Delta)$. The roles of $X$
and $Y$ can be changed, so that the statements above can be restated
as follows:
%
\begin{equation}
\label{statement} \max\bigl\{\bigl(\tau^x_l-
\tau^x_k\bigr)\vee n', \bigl(
\tau^y_l-\tau^y_k\bigr)\vee
n'\bigr\}\leq \min\bigl\{\bigl(\tau^x_l-
\tau^x_k\bigr)\vee n', \bigl(
\tau^y_l-\tau^y_k\bigr)\vee
n'\bigr\}(1+\Delta).
\end{equation}
%
Finally, let $F_n(n',\Delta)$ denote the event that for every $k$
that satisfies $\max\{\tau_k^x,\tau_k^y\}\leq n$, $F(k,n')$ holds.
Formally,
\[
F_n\bigl(n',\Delta\bigr):=\bigcup
_{i=0}^{n} \Biggl(\{K=i\}\cap \Biggl(\bigcap
_{k=0}^iH\bigl(k,n'\bigr) \Biggr) \Biggr),
\]
where $K$ is as in (\ref{K}). The event $F_n(n',\Delta)$ ensures
that (\ref{statement}) holds for every $k\leq K$. In particular, if
$(i,j)$ is a related pair such that $i\leq n$ and $j\leq n$ and
$(i',j')$ is another related pair, then
%
\begin{equation}
\label{statement1} \max\bigl\{\bigl|i-i'\bigr|\vee n',
\bigl|j-j'\bigr|\vee n'\bigr\}\leq\min\bigl
\{\bigl|i-i'\bigr|\vee n', \bigl|j-j'\bigr|\vee
n'\bigr\}(1+\Delta).
\end{equation}
%
\paragraph*{The bound on $P(F_n(n',\Delta))$} Let us first estimate
from below the probability of $F_n$. Since
\[
\bigl(F_n^x\bigr)^c=\Biggl\{\sum
_{i=1}^l D^x_i\leq n
\Biggr\}\cup\Biggl\{\sum_{i=1}^l
D^x_i>\beta n\Biggr\},
\]
by Hoeffding's inequality (recall that $l={\alpha\over p}n$)
\begin{eqnarray*}
P \Biggl(\sum_{i=1}^{l}
D^x_i-pl\leq n-pl \Biggr)&\leq& \exp\biggl[-2p
{(1-\alpha)^2\over\alpha}n\biggr],\\
 P \Biggl(\sum_{i=1}^{l}
D^x_i-pl>\beta n-pl \Biggr)&\leq&\exp\biggl[-2p
{(\beta- \alpha)^2\over
\alpha}n\biggr].
\end{eqnarray*}
Since ${(\beta-\alpha)^2\over\alpha}={(1-\alpha)^2\over
\alpha}={\Delta^2\over2(2+\Delta)}=:{d(\Delta)\over2},$ it holds
%
\begin{equation}
\label{kahepoolne} P\bigl(F_n^c\bigr)\leq2\exp[-pdn].
\end{equation}
%
For estimating $P(F_n(n',\Delta))$, we use the regenerativity
argument to see that for every $k$, the event $F(k,n')$ has the same
probability as $\bigcap_{m\geq n'}F_{m}$ so that by (\ref{kahepoolne}),
there exist constant $R(\Delta,p)<\infty,b_6(\Delta,p)>0$
\[
P\bigl(F^c\bigl(k,n'\bigr)\bigr)\leq\sum
_{m\geq n'}P\bigl(F^c_m\bigr)\leq2\sum
_{m\geq n'}\exp [-pdm] \leq R\exp\bigl[-b_6n'
\bigr].
\]
Finally, since $\bigcap_{k=0}^{n}F(k,n')\subset F_n(n',\Delta)$, we
have ($n'=A\ln n$)
%
\begin{eqnarray}
\label{H-est} P \bigl(F^c_n(n,\Delta) \bigr)&\leq&(n+1)P
\bigl(F^c\bigl(k,n'\bigr)\bigr)\leq M(n+1)\exp
\bigl[-b_6n'\bigr]\nonumber\\[-8pt]\\[-8pt]
&\leq&2Rn\exp\bigl[-b_6n'
\bigr]=2Rn^{1-Ab_6}.\nonumber
\end{eqnarray}
%
\subsection{The related pairs in extremal
alignments}\label{sec:abi} In previous subsection, we showed that
with the high probability the related pairs are rather uniformly
located almost in the main diagonal of the two-dimensional
representation of alignments (the $F$-event). We also know that with
high probability every piece of length $A\ln n$ of extremal
alignments contains at least one related pair (the $B$-event). Hence,
both extremal alignments cannot diverge from the main diagonal too
much and therefore they cannot be too far from
each other. The following lemma postulates this observation.

In the following, let $K^h$ and $K^l$ be the random number of
related pairs of the highest and lowest alignment, respectively. We
shall denote by $(i^{*h}_1,j^{*h}_1),\ldots,
(i^{*h}_{K^h},j^{*h}_{K^h})$ the related pairs of the highest
alignment and $(i^{*l}_1,j^{*l}_1),\ldots,
(i^{*l}_{K^l},j^{*l}_{K^l})$ the related pairs of the lowest
alignment. Let
\[
\overline{i}:=i^{*h}_{K^h}\wedge i^{*l}_{K^l}.
\]
%
We also agree that $i^{*h}_0:=:j^{*h}_0:=:i^{*l}_0:=:j^{*l}_0:=0$
and  with some abuse of terminology, we shall call also the pair
$(0,0)$ related of both highest and lowest alignments.
%
\begin{lemma}\label{abi} Let
$\Delta>0$ and assume that $B_n(n',\Delta)\cap F_n(n',{\Delta\over
2})$ holds. Let $(i^h,j^h)$ be a pair of  the highest alignment
of $X$ and $Y$ such that $i^h\leq\overline{i}$. Then there exists a
related pair $(i_u^{*l},j_u^{*l})$, $u\in\{0,\ldots,K^l\}$ of the
lowest alignment such that
%
\begin{eqnarray}
\label{max1} \bigl|i^h-i_u^{*l}\bigr|
\vee\bigl|j^h-j_u^{*l}\bigr| &\leq n'(1+
\Delta).
\end{eqnarray}
Moreover, there exists a related pair $(i_l^{*l},j_l^{*l})$, $l\in
\{0,\ldots,K^l\}$ of the lowest alignment such that
%
\begin{eqnarray}
\label{max2} i_l^{*l} \leq i^h \quad \mbox{and}\quad \bigl|j^h-j_l^{*l}\bigr| \leq2n'(1+
\Delta).
\end{eqnarray}
Similarly, for every pair $(i^l,j^l)$ of the lowest alignment of $X$
and $Y$ such that $i^l\leq \overline{i}$, there exists a related
pair $(i_u^{*h},j_u^{*h})$, $u\in\{0,\ldots,K^h\}$ such that
%
\begin{eqnarray}
\label{max3} \bigl|i^l-i_u^{*h}\bigr|
\vee\bigl|j^l-j_u^{*h}\bigr|\leq n'(1+
\Delta).
\end{eqnarray}
Moreover, there exists a related pair $(i_l^{*h},j_l^{*h})$ of the
highest alignment such that
%
\begin{eqnarray}
\label{max4} i^{*h}_l\leq i^l \quad \mbox{and}\quad
|j^l-j_l^{*h}| &\leq2n'(1+
\Delta).
\end{eqnarray}
\end{lemma}
%
%
\begin{pf} At first, we shall see that for every $0\leq t \leq
K^h-1$,
%
\begin{eqnarray}
\label{vahe1} \bigl(i^{*h}_{t+1}-i^{*h}_t
\bigr)\wedge\bigl(j^{*h}_{t+1}-j^{*h}_t
\bigr)&\leq& n',
\\
\label{vahe2} \bigl(i^{*h}_{t+1}-i^{*h}_t
\bigr)\vee \bigl(j^{*h}_{t+1}-j^{*h}_t
\bigr) &\leq& n'(1+\Delta).
\end{eqnarray}
Suppose there
exists $t$ such that (\ref{vahe1}) fails. The pairs
$(i^{*h}_t,j^{*h}_t)$ and $(i^{*h}_{t+1},j^{*h}_{t+1})$ are both in
the highest alignment, let it be $v$. Since $v$ is highest, the
restriction of $v$ between
\[
X_{i^{*h}_t+1},\ldots,X_{i^{*h}_{t+1}-1}, \quad \mbox{and}\quad  Y_{j^{*h}_t+1},
\ldots,Y_{j^{*h}_{t+1}-1}
\]
must be highest as well.
Denote $\tilde{n}=i^{*h}_{t+1}-1-i^{*h}_t$ and
$\tilde{m}=j^{*h}_{t+1}-1-j^{*h}_t$. If (\ref{vahe1}) does not hold,
then $\tilde{m},\tilde{n}\geq n'$. Suppose, without loss of
generality that $\tilde{m}\geq\tilde{n}$. Since
$F_n(n',{\Delta\over2})$ holds, then (\ref{statement1}) states that
$(\tilde{m}+1)\leq(\tilde{n}+1)(1+{\Delta\over2})$ implying that
$\tilde{m}\leq\tilde{n}(1+\Delta)$. Therefore, we have that the
sequences
\[
X_{i^{*h}_t+1},\ldots,X_{i^{*h}_{t}+\tilde{n}}, \quad \mbox{and}\quad  Y_{j^{*h}_t+1},
\ldots,Y_{j^{*h}_{t}+\tilde{m}}
\]
with $n'\leq
\tilde{n}\leq\tilde{m}\leq\tilde{n}(1+\Delta)$ have an optimal
alignment that contains no related pair. This contradicts
$B_n(n',\Delta)$. Hence, (\ref{vahe1}) holds. Since $t<K^h$, then
(\ref{statement1}) proves (\ref{vahe2}) (recall that
(\ref{statement1}) also holds for $i=j=0$).\vadjust{\goodbreak}

Consider an arbitrary (not necessarily related) pair
$(i^h,j^h)$ of the highest alignment so that $i^h\leq\overline{i}\leq
i^{*h}_{K^h}$. By (\ref{vahe2}), there
exists $0\leq k < K^h$ such that $i^{*h}_k\leq
i^h\leq i^{*h}_{k+1}$ and
%
\begin{equation}
\label{vahe3} \bigl(i^{*h}_{k+1}-i^{*h}_k
\bigr)\vee\bigl(j^{*h}_{k+1}-j^{*h}_k
\bigr)\leq n'(1+\Delta).
\end{equation}
Similarly, since $i^h\leq
\overline{i}\leq i^{*l}_{K^l}$, by applying (\ref{vahe2}) to the
lowest alignment, there exists $0\leq l < K^l$ such that
$i^{*l}_l\leq i^h\leq i^{*l}_{l+1}$ and
%
\begin{equation}
\label{vahe4} \bigl(i^{*l}_{l+1}-i^{*l}_l
\bigr)\vee\bigl(j^{*l}_{l+1}-j^{*l}_l
\bigr)\leq n'(1+\Delta).
\end{equation}
%
Hence $i^{*l}_l\leq i^h$ and $|i^{*l}_u-i^h|\leq n'(1+\Delta)$, for
$u=l,l+1$. For (\ref{max2}), it suffices to show that
$|j^h-j^{*l}_l|\leq2n'(1+\Delta)$. For (\ref{max1}), it suffices to
show that $\min_{u=l,l+1}|j^h-j^{*l}_u|\leq n'(1+\Delta)$. For that,
we consider three
cases separately:

(1) Suppose $i^{*h}_k\leq i_l^{*l}$. Because
$(i^{*h}_k,j^{*h}_k),(i^{*l}_l,j^{*l}_l),
(i^{*h}_{k+1},j^{*h}_{k+1})$ are related pairs and $i^{*l}_l\leq
i^h\leq i^{*h}_{k+1}$, we have $j^{*l}_l\leq j^{*h}_{k+1}$ so that
by $i^{*h}_k\leq i_l^{*l}$, it holds $j^{*h}_k\leq j^{*l}_l\leq
j^{*h}_{k+1}$. Clearly at least one inequality is strict. Since
$(i^{*h}_k,j^{*h}_k),(i^h,j^h), (i^{*h}_{k+1},j^{*h}_{k+1})$ are
aligned pairs, we have $j^{*h}_k\leq j^{h} \leq j^{*h}_{k+1}$ (with
at least one of the inequalities being strict). By
(\ref{vahe3}), we have $j^{*h}_{k+1}-j^{*h}_k\leq
n'(1+\Delta)$, implying that $|j^h-j^{*l}_l|\leq n'(1+\Delta)$. Thus,
(\ref{max1}) holds with $u=l$ and then
(\ref{max2}) trivially holds.

(2) Suppose $i_{l+1}^{*l}\leq i^{*h}_{k+1}$. The pairs
$(i^{*h}_k,j^{*h}_k),(i^{*l}_{l+1},j^{*l}_{l+1}),
(i^{*h}_{k+1},j^{*h}_{k+1})$ are related. Since $i^{*h}_{k}\leq
i^h\leq i_{l+1}^{*l}\leq i^{*h}_{k+1}$, we have that $j^{*h}_{k}\leq
j_{l+1}^{*l}\leq j_{k+1}^{*h}$ (again, at least one inequality is
strict). Since $(i^{*h}_k,j^{*h}_k),(i^h,j^h),
(i^{*h}_{k+1},j^{*h}_{k+1})$ are aligned pairs, we have
$j^{*h}_k\leq j^{h} \leq j^{*h}_{k+1}$. By
(\ref{vahe3}), we have $j^{*h}_{k+1}-j^{*h}_k\leq
n'(1+\Delta)$, implying that $|j^h-j^{*l}_{l+1}|\leq n'(1+\Delta)$.
Therefore, (\ref{max1}) holds for $u=l+1$. For (\ref{max2}), use the
inequalities (\ref{vahe4}) together with the inequalities
$|j^h-j^{*l}_{l}|\leq|j^h-j^{*l}_{l+1}|+|j^{*l}_l-j^{*l}_{l+1}|\leq
2n'(1+\Delta).$

(3) Suppose $i_l^{*l}<i^{*h}_k$ and $i^{*h}_{k+1}<i_{l+1}^{*l}$.
Since all pairs, except perhaps $(i^h,j^h)$ are related, we have
that $i_l^{*l}<i^{*h}_k\leq i^h\leq i^{*h}_{k+1}<i_{l+1}^{*l}$ and
$j_l^{*l}<j_k^{*h}\leq j^h\leq j^{*h}_{k+1}<j_{l+1}^{*l}$. By
(\ref{vahe4}), $|j^h-j_{l}^{*l}|\leq j_{l+1}^{*l}-j_l^{*l}\leq
n'(1+\Delta)$. Hence, (\ref{max1}) holds with $u=l$ and then
(\ref{max2}) trivially holds.

By symmetry, the second statement of the lemma holds by the same
argument.
\end{pf}
%

Recall the definition of $\Delta_n:=\sqrt{{16 \over p}{\ln n \over
n}}$.
%
\begin{lemma}\label{nurk} Let $1>\Delta>0$ and assume that
$B_n(n',\Delta)\cap G_n(\Delta_n) \cap
F_n(n',{\Delta\over2})$ holds. Then there exists
$n_1(\Delta)<\infty$ and $M(\Delta)<\infty$ so that for every
$n>n_1$, $n-\overline{i}\leq M\Delta_n n$.
\end{lemma}
%
%
\begin{pf} Let $(i^*,j^*):= (i^{*h}_{K^h},j^{*h}_{K^h})$. Since
$G_n(\Delta_n)$ holds,
there exists a related pair $(i',j')$ so that $i',j'\geq
(1-\Delta_n)n$. Without loss of generality, we can take $(i',j')$
the last related pair satisfying $i'\leq n$ and $j'\leq n$ so that
$i'\geq i^*$ and $j'\geq j^*$. Let now $M(\Delta)$ be so big that
%
\begin{equation}
\label{Mdef} 1<(M-1){\Delta\over2+\Delta}.
\end{equation}
First, we shall show that $n-i^*\leq M \Delta_n n$. If not, then for
$n$ big enough,
%
\begin{equation}
\label{suur} i'-i^*> n(M \Delta_n-
\Delta_n)=n \Delta_n(M-1)>n' .
\end{equation}
Then, by the definition of $M$
%
\begin{equation}
\label{id} n \Delta_n\leq n\Delta_n (M-1)
{{\Delta/2}\over1+{\Delta/
2}}\leq\bigl(i'-i^*\bigr){{\Delta/2}\over1+{\Delta/2}}
\leq \bigl(i'-i^*\bigr){\Delta\over2} .
\end{equation}
We shall now show that when
(\ref{suur}) holds, then
%
\begin{equation}
\label{vord} \bigl(n-i^*\bigr)\leq\bigl(n-j^*\bigr) (1+\Delta),\qquad  \bigl(n-j^*
\bigr)\leq\bigl(n-i^*\bigr) (1+\Delta).
\end{equation}
Consider two cases separately:

(a) $i'-i^*\leq j'-j^*$. Since $F_n(n',{\Delta\over2})$ holds, we
have that $j'-j^*\leq(i'-i^*)(1+{\Delta\over2})$ so that
\[
n-j^*=\bigl(n-j'\bigr)+\bigl(j'-j^*\bigr)\leq n
\Delta_n+\bigl(i'-i^*\bigr) \biggl(1+
{\Delta\over2}\biggr)\leq \bigl(i'-i^*\bigr) (1+{\Delta})
\leq\bigl(n-i^*\bigr) (1+\Delta),
\]
where the second last inequality holds due to (\ref{id}). We also
have that
\[
n-i^*=\bigl(n-i'\bigr)+\bigl(i'-i^*\bigr)\leq n
\Delta_n+\bigl(j'-j^*\bigr)\leq\bigl(i'-i^*
\bigr){\Delta
\over2}+\bigl(j'-j^*\bigr)\leq
\bigl(j'-j^*\bigr) \biggl(1+{\Delta\over2}\biggr).
\]
%

(b) $i'-i^*\geq j'-j^*$. By $F_n(n',{\Delta\over2})$ we have again
that $i'-i^*\leq(j'-j^*)(1+{\Delta\over2})$ so that by (\ref{id}),
we have
%
\begin{equation}
\label{id2} n\Delta_n\leq\bigl(i'-i^*\bigr)
{{\Delta/2}\over1+{\Delta/2}}\leq \bigl(j'-j^*\bigr){\Delta\over2}
\end{equation}
and arguing similarly as in
the case (a), we now obtain
\[
n-j^*\leq\bigl(i'-i^*\bigr) \biggl(1+{\Delta\over2}
\biggr),\qquad  n-i^*\leq\bigl(j'-j^*\bigr) (1+{\Delta}).
\]
We are now applying the same argument as in the previous lemma.
Recall that $(i^*,j^*)$ belongs to the highest alignment. The
restriction of the highest alignment between
\[
X_{i^*+1},\ldots,X_n\quad  \mbox{and}\quad  Y_{j^*+1},
\ldots,Y_n
\]
must be highest as well. Moreover, the restriction contains no
related pairs. The lengths of $X_{i^*+1},\ldots,X_n$ and
$Y_{j^*+1},\ldots,Y_n$ are $n-i^*$ and $n-j^*$, respectively. By
(\ref{suur}) and (\ref{vord}), for $n$ big enough, both lengths are
bigger than $n'$; by (\ref{vord}) their lengths are comparable, so
that by event $B_n(n',\Delta)$, they have to contain a related pair.
This contradicts the assumption that $(i^*,j^*)$ is the last related
pair of the highest alignment. The contradiction is due to
assumption $n-i^*> M \Delta_n n$. Hence, $n-i^* \leq M \Delta_n n$,
eventually. The same argument holds for the lowest alignment, hence
$\overline{i} \geq n- M \Delta_n n$, eventually.
\end{pf}
%
\subsection{Proof of Theorem \texorpdfstring{\protect\ref{main1}}{1.1}}\label{sec:proof1}

Choose $1>\Delta>0$ so small that (\ref{mainconditiondelta}) holds.
Let $M(\Delta)$ be defined as in (\ref{Mdef}) and $\alpha_n=M
\Delta_n$. Clearly $\alpha_n\to0$, in particular, $\alpha_n<1$ for
$n$ big enough. Recall the definition of Hausdorff's distance
between extremal alignments both represented as a set of
2-dimensional points. More precisely, let $H$ and $L$ be the
highest and lowest alignments, both represented as the set of
two-dimensional points. Clearly, $|H|=|L|=L_n$. In the statement of
Theorem~\ref{main1}, the subsets of $H$ and $L$, where the pairs
$(i,j)$ satisfying $i>n-\alpha_n n$ are left out, are considered.
More precisely, we consider the consider the points
\[
H_o:=\bigl\{\bigl(i^h,j^h\bigr)\in H:
i^h\leq n(1- \alpha_n)\bigr\},\qquad  L_o:=\bigl\{
\bigl(i^l,j^l\bigr)\in L: i^l\leq n(1-
\alpha_n)\bigr\}.
\]
If for an arbitrary element $(i^h,j^h)$ of $H_o$,
there exists an element $(i^l,j^l)$ of $L$ such that $|i^h-i^l|\vee
|j^h-j^l|\leq(1+\Delta)n'$, then $\max_{(i,j)\in H_o}
\min_{(i^l,j^l)\in L}|i-i^l|\vee|j-j^l|\leq(1+\Delta)n'$. If, in
addition, for an arbitrary element $(i^l,j^l)$ of $L_o$, there
exists an element $(i^h,j^h)$ of $H$ such that $|i^h-i^l|\vee
|j^h-j^l|\leq(1+\Delta)n'$, then the restricted Hausdorff's
distance with respect to the maximum norm between $H$ and $L$ is at
most $(1+\Delta)n'$. The restricted Hausdorff's distance between $H$
and $L$ with
respect to the $l_2$-norm is then $\sqrt{2}n'(1+\Delta)$.

\begin{pf*}{Proof of Theorem \protect\ref{main1}}
Choose $1>\Delta>0$ so small that (\ref{mainconditiondelta}) holds.
Now, let $M:=M(\Delta)$ as in (\ref{Mdef}), $b_4:=b_4(\Delta)>0$ and
$n_o(\Delta)$ be as in the bound (\ref{B-est}),
$b_6:=b_6({\Delta\over2})>0$ be as in (\ref{H-est}). Let, moreover,
$n>n_1\vee n_o$, where $n_1(\Delta)$ is as in Lemma~\ref{nurk} and
let
$h_o$ be
the restricted Hausdorff's distance with respect to the maximum norm
and $\alpha_n:=M\Delta_n=M\sqrt{16 \ln n\over p n}$. Since $n>n_1$,
$\alpha_n<1$ so that $h_o$ is correctly defined. Finally, choose $A$
so big that $\min\{Ab_4,Ab_6\}\geq3$. We aim to bound the
probability of the event $E_n:=\{h_o(L,H)\leq2n'\}$, where $n'=A\ln
n$. If the event $B_n(n',\Delta)\cap G_n(\Delta_n)\cap
F_n(n',{\Delta\over2})$ holds, then by Lemma~\ref{nurk},
$\overline{i}\geq n(1- \alpha_n)$ so that for every $(i^h,j^h)\in
H_o$ and $(i^l,j^l)\in L_o$ Lemma~\ref{abi} applies. Since
$(1+\Delta)<2$, (\ref{max1}) and (\ref{max3}) of Lemma~\ref{abi}
ensure that $h_o(H,L)\leq2n'$. Therefore,
\[
B_n\bigl(n',\Delta\bigr)\cap G_n(
\Delta_n)\cap F_n\biggl(n',
{\Delta\over2}\biggr)\subset E_n.
\]
Hence, from
(\ref{B-est}), (\ref{Ggn}) and (\ref{H-est}), for $n>n_o\vee n_1$,
\begin{eqnarray*}
P\bigl(E^c_n\bigr)&\leq& P \bigl(B^c_n
\bigl(n',\Delta\bigr) \bigr)+P \biggl(G^c_n
\biggl(n',{\Delta
\over
2}\biggr) \biggr)+P
\biggl(F^c_n\biggl(n',{\Delta\over
2}
\biggr) \biggr)\\
&\leq&2Bn^{1-Ab_4}+6n^{-2}+ 2Rn^{1-Ab_6}
\\
&\leq&2(R+B+3)n^{-2}.
\end{eqnarray*}
Thus, the theorem holds with
$D=2(R+B+3)$ and $C=2A$.
\end{pf*}
%
\subsection{Proof of Theorem \texorpdfstring{\protect\ref{main2}}{1.2}}\label{sec:proof2} In
Theorem~\ref{main1},
we used the 2-dimensional representation of alignments, so an
alignment were identified with a finite set of points. In the
alignment graph, these points are joined by a line. We consider the
highest and lowest alignment graphs, and we are interested in the
maximal vertical (horizontal) distance between these 2 piecewise
linear curves. This maximum is called vertical (horizontal) distance
between lowest and highest alignment
graphs.
\begin{pf*}{Proof of Theorem \protect\ref{main2}}
From Lemma~\ref{nurk} and (\ref{max2}) of Lemma~\ref{abi}, it
follows that on the event $F_n(n',{\Delta\over2})\cap
G(\Delta_n)\cap B_n(n',\Delta)$ the following holds: for every pair
$(i^h,j^h)$ of the highest alignment such that $i^h\leq
\overline{i}$, there exists a pair $(i_k^l,j_k^l)$ (including the
possibility that $k=0$) of the lowest alignment such that $i^l_k\leq
i^h$ and $|j^h-j^l_k|\leq2n'(1+\Delta)$. Recall that $L$ and $H$
are the lowest and highest alignment graphs, respectively. Since $L$
in non-decreasing, it follows that $H(i^h)-L(i^h)=j^h-L(i^h)\leq
j^h-j^l_k\leq2n'(1+\Delta).$ By (\ref{max4}) of Lemma~\ref{abi}, we
obtain (using the same argument) that for every pair $(i^l,j^l)$ of
the lowest alignment such that $i^l\leq\overline{i}$, the following
inequality holds: $H(i^h)-L(i^h)\leq2n'(1+\Delta)$. Since the
function $H-L$ is piecewise linear, we obtain that $\sup_{x\in
[0,\overline{i}]} (H(x)-L(x) )\leq2n'(1+\Delta)$.

The rest of the proof is the same as the one of Theorem~\ref{main1}.
\end{pf*}
%
\section{Proof of Theorem \texorpdfstring{\protect\ref{main3}}{1.3}}\label{random}
When dealing with the sequences of random lengths, it is more
convenient to consider the locations of ancestors. Recall the i.i.d.
vectors $U_i$ as defined in (\ref{vectors}). Thus, given $k, l\in
\mathbb{N}$ $(k<l)$, with some abuse of terminology, we shall call
the highest (lowest) alignment of $U_{k+1},\ldots,U_{k+l}$ the
highest (lowest) alignment between these $X$ and $Y$ sequences that
have the ancestors in the interval $[k+1,k+l].$ Note that these
sequences as well as corresponding optimal alignments are all
functions of $U_{k+1},\ldots,U_{k+l}$, only. This justifies the
terminology. Hence, the highest alignment of the random lengths
sequences $X$ and~$Y$ (as defined above) is the
highest alignment of $U_{1},\ldots,U_{m(n)}$.

Let now $k=0$ and let $l\geq1$ be fixed. We shall consider the
vectors $U_{1},\ldots,U_{l}$ and the corresponding $X$ and $Y$
sequences. Thus, $n_x(l):=\sum_{j=1}^l D_j^x$, $n_y(l):=\sum_{j=1}^l
D_j^y$ are their lengths. Let us define the events
\[
A_l^x:=\biggl\{\bigl|n_x(l)-lp\bigr|<
{p\over2}l\biggr\},\qquad  A_l^y:=\biggl
\{\bigl|n_y(l)-lp\bigr|<{p\over
2}l\biggr\},\qquad  A_l:=A_l^x
\cap A_l^y.
\]
For a fixed $\Delta>0$, let
\[
C_l(\Delta):=A_l\cap\biggl\{\bigl|n_x(l)-n_y(l)\bigr|<
{p\over2}\Delta l\biggr\}.
\]
Hence, on $C_l$ the following inequalities hold true:
\[
l{p\over2}<n_x(l)\wedge n_y(l)\leq
n_x(l)\vee n_y(l)<{3p\over
2}l,\qquad
n_x(l)\vee n_y(l)\leq \bigl( n_x(l)\wedge
n_y(l) \bigr) (1+\Delta).
\]
Let
\[
\mathcal{B}_l(\Delta):= \biggl\{(\tilde{n},\tilde{m})\in\mathbb
{N}^2: l{p\over2}< \tilde{n} \wedge\tilde{m}\leq
\tilde{n}\vee \tilde{m}<{3p\over
2}l, \tilde{n}\vee \tilde{m}\leq(
\tilde{n}\wedge \tilde{m}) (1+\Delta) \biggr\}.
\]
Thus,
\begin{equation}
\label{kapp} C_l\subset\bigcup_{(\tilde{n},\tilde{m})\in\mathcal{B}_l(\Delta
)}
\bigl\{n_x(l)=\tilde{n},n_y(l)=\tilde{m}\bigr\}.
\end{equation}
%
With applying Hoeffding's inequality three times, it is easy to see
the existence of a constant $c_1$ (depending on $\Delta$ and $p$) so
that
%
\begin{equation}
\label{CC} P(C_l)\geq1-6\exp[-c_1l].
\end{equation}
%
\subsection*{The $B$-event for the sequences of random lengths} We
shall now study the random lengths analogue of the $B$-events.
Recall that the event $B_0(\tilde{n},\tilde{m})$ states that the
highest alignment between the sequences $X_1,\ldots,X_{\tilde{n}}$
and $Y_1,\ldots,Y_{\tilde{m}}$ contains a related pair. We shall
define now the event
\[
E_k(l):=\{\mbox{the highest alignment of } U_{k+1},
\ldots,U_{k+l} \mbox{ contains a related pair}\}.
\]
%
We shall bound the probability of $P(E_k(l))$. Clearly
$P(E_k(l))=P(E_0(l))$ for every $k=1,2,\ldots\,$, hence we shall
consider the event $E_0(l)$. Obviously,
\[
\bigcup_{(\tilde{n},\tilde{m})\in\mathcal{B}_l(\Delta)} \bigl(B_0(\tilde{n},
\tilde{m})\cap \bigl\{n_x(l)=\tilde{n},n_y(l)=\tilde{m}
\bigr\} \bigr)\subset E_0(l).
\]
Since
\[
P \bigl(B_0(\tilde{n},\tilde{m})\cap \bigl\{n_x(l)=
\tilde{n},n_y(l)=\tilde{m}\bigr\} \bigr)\geq P\bigl(n_x(l)=
\tilde{n},n_y(l)=\tilde{m}\bigr)-P \bigl(B_0^c(
\tilde{n},\tilde {m}) \bigr).
\]
Since $\tilde{n}$ and $\tilde{m}$ belong to $\mathcal{B}_l(\Delta
)$, by
Lemma~\ref{related-pair2}, we obtain
\[
P \bigl(B_0^c(\tilde{n},\tilde{m}) \bigr)\leq\exp
\biggl[-b_3 {p\over
2}l\biggr],
\]
provided that $l$ is big and $\Delta$ small enough. Thus, by (\ref
{kapp}) and (\ref{CC}), we obtain
\begin{eqnarray*}
P\bigl(E_0(l)\bigr)&\geq&\sum_{(\tilde{n},\tilde{m})\in\mathcal{B}_l(\Delta
)}
\bigl(P\bigl(n_x(l)=\tilde{n},n_y(l)=\tilde{m}\bigr)-P
\bigl(B_0^c(\tilde {n},\tilde{m}) \bigr) \bigr)
\\
&\geq& P(C_l)-\bigl|\mathcal{B}_l(\Delta)\bigr| \exp
\biggl[-b_3 {p\over2}l\biggr]\geq 1-4
\exp[-c_1 l]-(pl)^2\exp\biggl[-b_3
{p\over2}l\biggr].
\end{eqnarray*}
Hence, there
exists $l_o$ and a constant $c_2>0$ (both depending on $\Delta$ and
$p$) so that for any $l>l_o$,
%
\begin{equation}
\label{El} P\bigl(E_0(l)\bigr)\geq1-\mathrm{e}^{-c_2 l}.
\end{equation}
%
Let now $\underline{l}(n):={2\over p}A \ln n$ and we assume $n$ to
be fixed and so big that $\underline{l}>l_o$ so that (\ref{El})
holds for any $l\geq\underline{l}$. Let, for any $k=0,1,\ldots$
\[
E_k:=\bigcup_{l\geq\underline{l}}E_k(l),\qquad
E^h:=\bigcup_{k=0}^{m-\underline{l}}
E_k.
\]
When the event $E^h$ holds, then the following is true: the highest
alignment of $U_{k+1},\ldots,U_{k+l}$ contains a related pair
whenever $l\geq\underline{l}$ and $k+l\leq m$. By (\ref{El}), we
obtain the estimate
%
\begin{equation}
\label{E} P\bigl(E_k^c\bigr)\leq\sum
_{l\geq\underline{l}}P\bigl(E^c_k(l)\bigr)\leq\sum
_{l\geq
\underline{l}}\mathrm{e}^{-c_2 l}\leq K \mathrm{e}^{-c_2 \underline{l}}=K
n^{-2c_2
A/ p},
\end{equation}
where $K$ is a constant. Thus,
\[
P \bigl(\bigl(E^h\bigr)^c \bigr)\leq
m(n)Kn^{-2c_2 A/ p}={n\over p}Kn^{-2c_2
A/
p}=
{K\over p}n^{1-{2c_2 A/ p}}.
\]
The event $E^h$ was defined
for the highest alignment. Similar event, let it be $E^l$ can be
defined for the lowest alignment. The bound (\ref{E}) holds also for
$E^l$. Hence, with $E:=E^h\cap E^l$, we obtain that $P(E)\geq
1-2Kp^{-1} n^{-2c_2 A/ p}$. Now we are ready to prove Theorem~\ref{main3}.
%
\begin{pf*}{Proof of Theorem \protect\ref{main3}} Choose $1>\Delta>0$ so small that (\ref{mainconditiondelta})
holds. Now let $c_2(\Delta)$ be as in (\ref{El}) and choose $A$ so
big that ${2c_2(\Delta)A\over p}>3$. By (\ref{E}), the event $E$
holds then with probability at least $1-2Kp^{-1}n^{-2}$. Now proceed
as in the proof of Lemma~\ref{abi}. Let $a^h_1,\ldots, a^h_{K^h}$
the ancestors of all related pairs in the highest alignment. Let
$a^h_0:=0$ and $a^h_{K^h+1}:=m+1$. Assume that $E$ holds. Then we
have that for every $k=0,\ldots,K^h$, $a^h_{k+1}-a^h_k <
\underline{l}$. Hence, for every pair of the highest alignment
$(i^h,j^h)$, there exists $k\in\{0,\ldots, K^h\}$ such that
$a^h_k\leq a^x(i^h)\leq a^h_{k+1}$ so that $|a^x(i^h)-a^h_k|\vee
|a^x(i^h)-a^h_{k+1}|\leq\underline{l}={2\over p}A \ln n$. Clearly,
\[
\bigl(i^{*h}_{k+1}-i^{*h}_k\bigr)\vee
\bigl(j^{*h}_{k+1}-j^{*h}_k\bigr)\leq
\bigl|a^h_{k+1}-a^h_k\bigr|\leq
{2\over p}A \ln n
\]
and also
$|i_k^{*h}-i^h|\vee
|i_{k+1}^{*h}-i^h|\leq{2\over p}A \ln n$.

Similarly, by $E^l$, there exists $0\leq l\leq K^l$ so that
$a^l_l\leq a^x(i^h)\leq a^l_{l+1}$, where $a^l_l$, $k=1,\ldots,K^l$
are the ancestors of the related pairs in the lowest alignment,
$a^l_0:=0$ and $a^l_{K^l+1}:=m+1.$ Thus,
\[
\bigl(i^{*l}_{l+1}-i^{*l}_l\bigr)\vee
\bigl(j^{*l}_{l+1}-j^{*l}_l\bigr)\leq
\bigl|a^l_{l+1}-a^l_l\bigr|\leq
{2\over p}A \ln n
\]
and also
$|i_l^{*l}-i^h|\vee|i_{l+1}^{*l}-i^h|\leq{2\over p}A \ln n$. Hence,
the inequalities (\ref{vahe3}) and (\ref{vahe4}) hold. Now proceed
as in the proof of Lemma~\ref{abi} and Theorem~\ref{main1} to see
that
\[
P \biggl(h(H,L)>{2\over p}A \ln n\biggr)\leq P
\bigl(E^c\bigr)\leq2Kp^{-1} n^{-2}
\]
so that (\ref{v3}) holds with $C_r:={2\over p}A$ and
$D_r:=2Kp^{-1}$, where $K$ is as in (\ref{E}).
\end{pf*}

\section{Simulations}\label{sec:simulations}
We now present some simulations about the growth of the distance
between the extremal alignments as well as another statistics. In
simulations, for different $n$-s up to 10\,000, 100 pairs of i.i.d.
sequences of length $n$ with were generated. Half of them were
independent i.i.d. sequences with $X_1$ and $Y_1$ distributed
uniformly over four letter alphabet. Another half of the sequences
were related with following parameters: the common ancestor process
$Z_1,Z_2,\ldots$ is i.i.d. with $Z_1$ being uniformly distributed
over four letter. The mutation matrix for generating $X$ and $Y$
sequences were the following:
\[
\bigl(P\bigl(f_1(Z_1)=a_j|Z_1=a_i
\bigr) \bigr)_{i,j=1,\ldots,4}=\lleft( %
\begin{array} {c@{\quad}c@{\quad}c@{\quad}c} 0.9 & 0.02
& 0.02 & 0.06
\\
0.02 & 0.9 & 0.06 & 0.02
\\
0.02 & 0.06 & 0.9 & 0.02
\\
0.06 & 0.02 & 0.02 & 0.9
\\
\end{array} %
\rright).
\]
The deletion probability $1-p=0.05$. Thus, the mutation matrix is
such that $X_1,X_2,\ldots$ and $Y_1,Y_2,\ldots$ were, as for
unrelated case, i.i.d. sequences with $X_1$ and $Y_1$ distributed
uniformly over four letter alphabet, but the sequences $X$ and $Y$
are clearly not independent any more. The same models were used in
generating Figures~\ref{fig1} and \ref{fig2}. Since, the marginal distributions of $X$
and $Y$ are uniform, we have $p_o=\overline p={1\over4}$ and
$q={3\over4}$. From the mutation matrix, it follows that for
related sequences $\overline{q}=1-0.02=0.98$. Hence, for related
sequences $\rho={\overline{q}\over q}>1$ and the left-hand side of
(\ref{maincondition1}) is (clearly $\gamma_{\tt{R}}>0.5$)
\begin{eqnarray*}
&&\gamma_{\tt{R}}\log_2 \overline{p}+(1-\gamma_{\tt{R}})
\bigl(\log _2(\overline{q}q)+\log_2\bigl(
\overline{q}q^{-1}\bigr)\bigr)+2h(\gamma_{\tt
{R}})\\
&&\quad = -2
\gamma_{\tt{R}}+2(1-\gamma_{\tt{R}})\log _2(0.98)+2h(
\gamma_{\tt{R}}).
\end{eqnarray*}
Hence, (\ref{maincondition1}) in this case is
%
\begin{equation}
\label{maincondition2} h(\gamma_{\tt{R}})<\gamma_{\tt{R}}-(1-
\gamma_{\tt
{R}})\log_2(0.98).
\end{equation}
The condition (\ref{maincondition2}) holds, if $\gamma_{\tt{R}}$ is
big enough. Since the solution of $h(x)=x-(1-x)\log_2(0.98)$ is
about 0.770481, (\ref{maincondition2}) holds if and only if
$\gamma_{\tt{R}}> 0.770481$. From Figure~\ref{fig2}, we estimate $\gamma
_{\tt{R}}$ as ${747\over949}=0.787$. Another simulations confirm that
$\gamma_{\tt{R}}$ is somewhere around 0.79. Thus, it is reasonable
believe that for our model, the condition (\ref{maincondition1})
holds true. Recall that for independent sequences,
(\ref{maincondition1}) always fails.

In both cases -- related and unrelated sequences -- the average of
following statistics were found: $L_n$, the horizontal length of
the maximum non-uniqueness stretch, the maximum vertical distance
and the maximum (full) Hausdorff's distance.\vadjust{\goodbreak}

\begin{figure}

\includegraphics{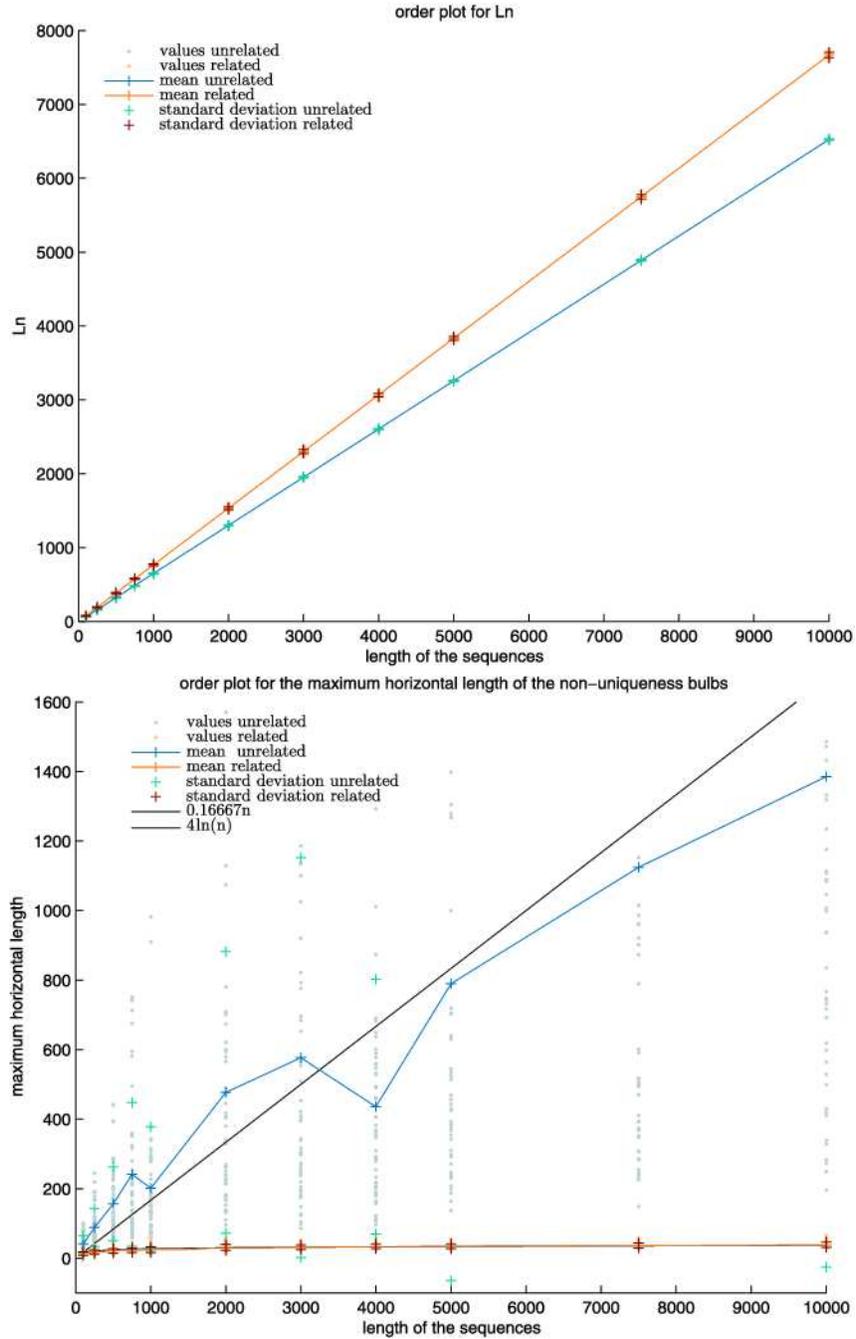}

\caption{Growth of $L_n$ (top). Growth of non-uniqueness
stretch (bottom).}\label{fig4}
\end{figure}

The top-plot in Figure~\ref{fig4} shows the growth of $L_n$ as $n$ grows.
The standard deviation around the means are marked with crosses. For
independent case, the crosses are almost overlapping implying that
the deviation is relatively small. As the picture shows, the growth
of $L_n$ is linear in both cases, the slope, however, is different:
the upper line corresponds to the related sequences, the lower line
is for independent sequences.
%

The bottom-plot in Figure~\ref{fig4} shows the horizontal length of maximum
non-uniqueness stretch. For independent sequences (upper curve), the
growth is, perhaps, smaller than linear but considerably faster than
logarithmic. The straight line is, in some sense, the best linear
approximation. The $+$-signs mark the standard deviation
around the mean that in this case is rather big, meaning that these
simulations do not give enough evidence to conclude the non-linear
growth. For related sequences (lower curve), the growth is clearly
logarithmic because it almost overlaps with the $4\ln n$-curve. We
also point out that the standard deviation for this case is
remarkably smaller and this only confirms the logarithmic growth.

\begin{figure}

\includegraphics{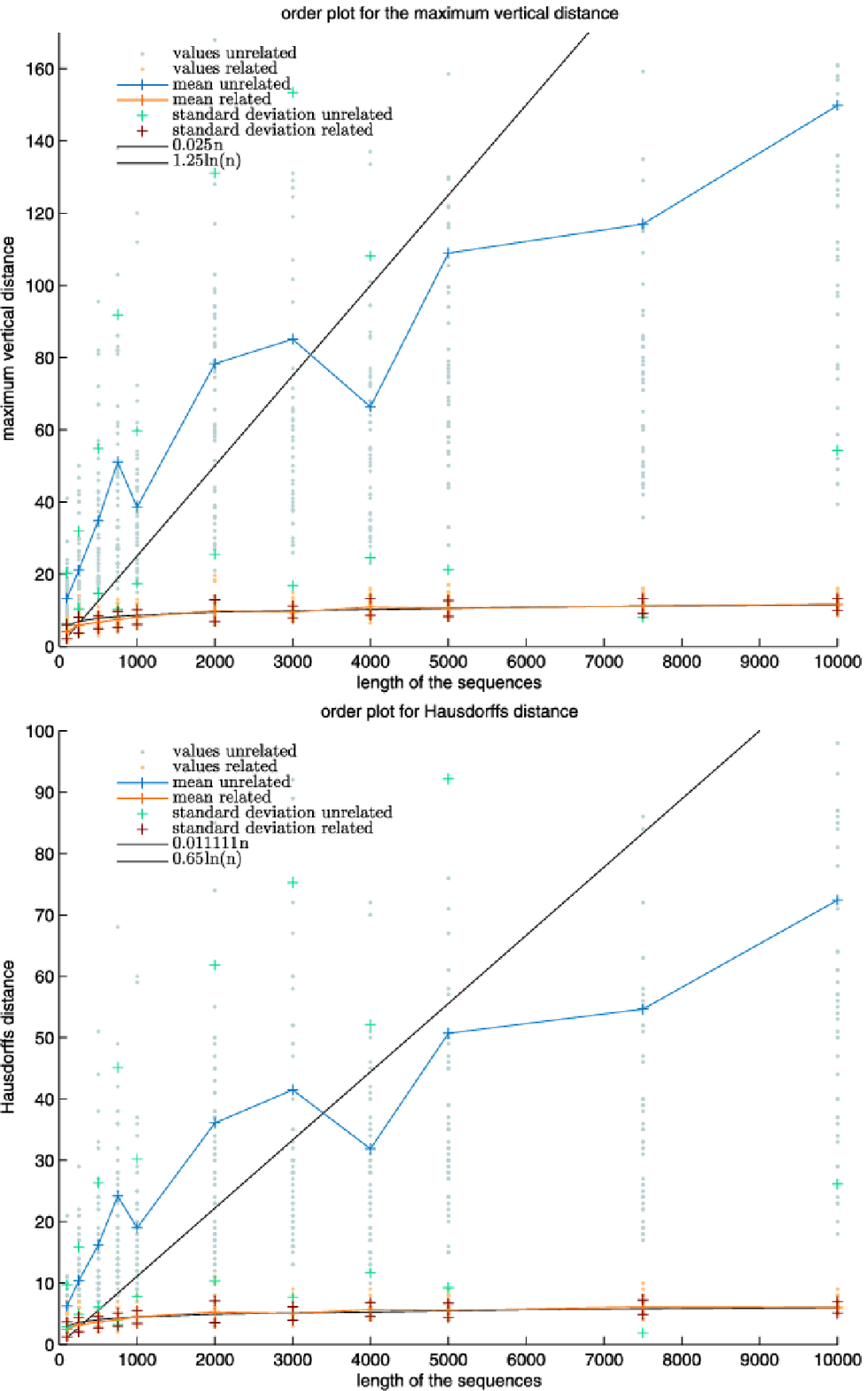}

\caption{Growth of maximal vertical distance (top).
Growth of Hausdorff's distance (bottom).}\label{fig5}
\end{figure}

In Figure~\ref{fig5}, the maximum vertical distance (top) and (full)
Hausdorff's distance with respect to the maximum-norm (bottom) are
plotted. Both pictures are similar to the bottom picture of Figure~\ref{fig4}
and can be interpreted analogously. For the related case, the growth
is clearly logarithmic (the best approximation is $1.25 \ln n$ for
maximum vertical distance, and $0.65 \ln n$ for Hausdorff's
distance) and that is a full correspondence with Theorems
\ref{main1}, \ref{main2} and \ref{main3}. Note that we have used the
full Hausdorff's distance instead of the restricted one so that the
simulations confirm the conjecture that Theorems \ref{main1} and
\ref{main2} also hold with $h$ instead of $h_o$.
%

\begin{figure}

\includegraphics{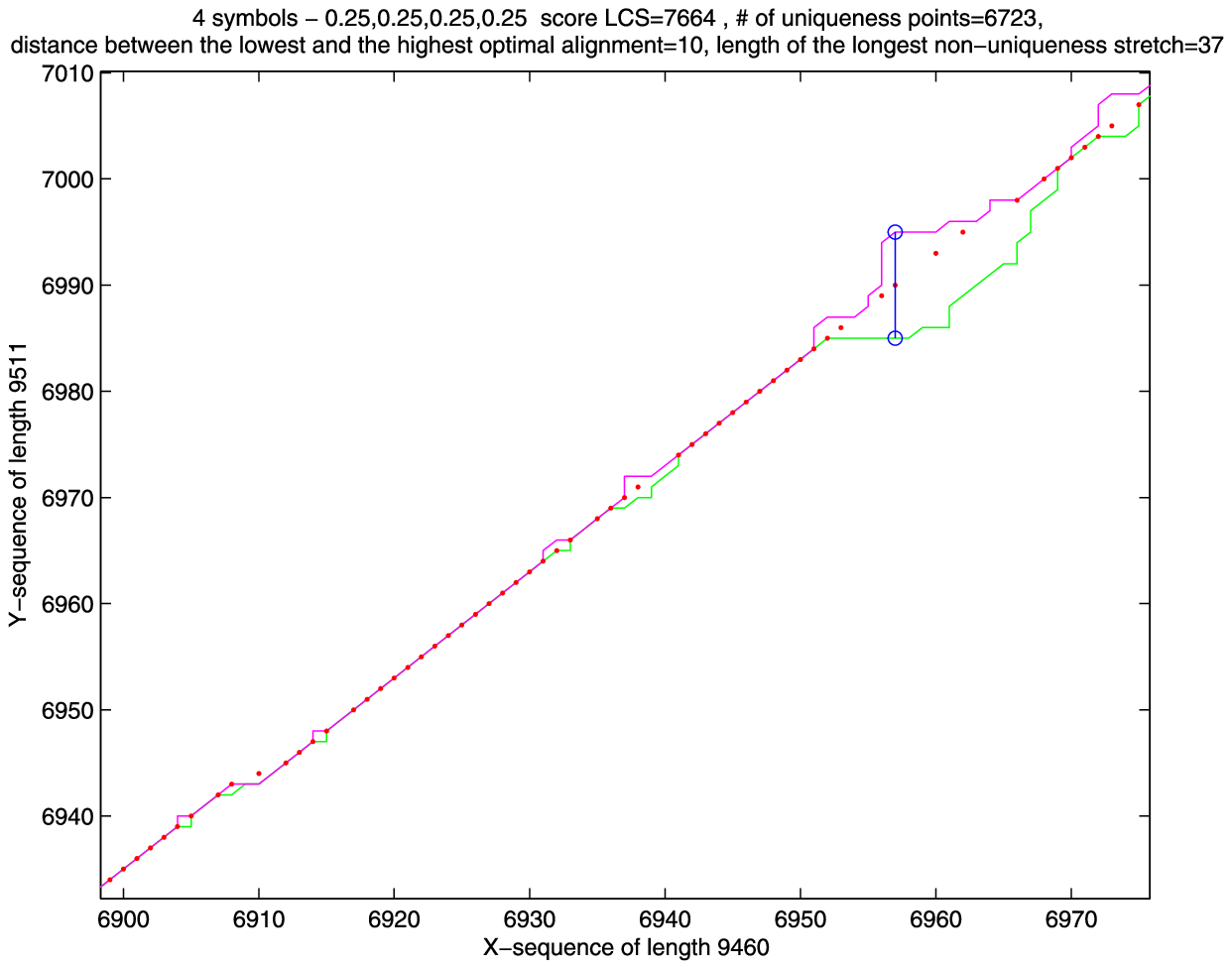}

\caption{The related pairs with the color of the common ancestor (dots) and extremal alignments.}\label{fig6}
\end{figure}

In Figure~\ref{fig6}, there is a zoomsection fragment of two extremal
alignment of the related sequences. Recall the definition of related
pairs -- the corresponding sites have the same ancestor. It does not
necessarily mean that they have the color of the common ancestor,
but often it is so. In the last picture, the related pars \textit{with
the color of the common ancestor} are marked with dots. Note
that in some small region, there are relatively many those pairs, on
same other region, there are less those pairs. The picture (and
other similar simulations) also shows that in the regions with many
these pairs, the extremal alignments coincide with them. This means
that in both sequences, there are parts that relatively less mutated
and the behavior of the extremal alignments indicate the existence
of such region rather well. In the area with relatively few
dots, the extremal alignments fluctuate indicating that in this part
(at least in one sequence) many mutations have been occurred. Hence,
based on these simulations, we can conclude the extremal alignments
are rather good tools for finding the less mutated regions and
obtaining information about the common ancestor.
%
%
\begin{appendix}
\setcounter{equation}{0}
\section*{Appendix}\label{app}
In the following, let $X_1,X_2,\ldots, Y_1,Y_2,\ldots$ be related
sequences. Recall, that our model for related sequences incorporates
the independent case. Recall the convergence (\ref{piirv-a}):
\[
{1\over n} L(X_1,\ldots, X_{\lfloor n a\rfloor};
Y_1,\ldots,Y_n)\to\gamma_{\tt{R}}(a),\qquad  \mbox{a.s.}
\]
%

\begin{lemma}\label{lemma:appendix} For every $0< a < 1$, $\gamma
_{\tt{R}}(a)<\gamma_{\tt{R}},$ for every $a>1$, $\gamma
_{\tt{R}}(a)>\gamma_{\tt{R}}$.
\end{lemma}
\begin{pf} Clearly the function $a\mapsto\gamma_{\tt{R}}(a)$
is nondecreasing in $a$ and
there exists a $K\in\mathbb{N}$ so big that $\gamma_{\tt
{R}}(K)>\gamma_{\tt{R}}$.

Fix $0<a<1$ and choose $\epsilon>0$ be so small that ${1-a\over
\epsilon}>K$. For every $m\in\mathbb{Z}$, $a,b>0$ let
\[
L_{m:an,bn}:=L(X_{m+1},\ldots,X_{\lfloor
na\rfloor};
Y_{m+1},\ldots,Y_{\lfloor nb\rfloor}).
\]
Let $ c:=1-a$.
By superadditivity,
\begin{eqnarray*}
L_{n(1+c),n}\geq L_{n(1-\epsilon),n(1-\epsilon)}+ L_{\lfloor
n(1-\epsilon)\rfloor: n(1+c), n}.
\end{eqnarray*}
Let $m=\lfloor
n(1-\epsilon)\rfloor$. Since, $c>K\epsilon$ and for every $c\geq0$,
\[
\bigl\lfloor n(1-\epsilon)\bigr\rfloor+ \bigl\lfloor n(c+\epsilon)\bigr\rfloor
\leq \bigl\lfloor\bigl\lfloor n(1-\epsilon)\bigr\rfloor+n(c+\epsilon)\bigr\rfloor
\leq \bigl\lfloor n(1+c)\bigr\rfloor,
\]
it holds
\begin{eqnarray*}
L_{m:n(1+c),n}\geq L_{m:m+n(c+\epsilon),m+n\epsilon}\geq L_{m: m+
\lfloor n\epsilon\rfloor K,m+\lfloor n\epsilon\rfloor}.
\end{eqnarray*}
Clearly
%
\begin{eqnarray}
\label{koond1} \lim_n{1\over n}
L_{n(1-\epsilon),n(1-\epsilon)}=(1-\epsilon)\lim_{u\to\infty}{L_u
\over u}=(1-
\epsilon)\gamma_{\tt{R}} \qquad \mbox{a.s.}
\end{eqnarray}
Let us now show that
%
\begin{equation}
\label{koond2} \lim_n{1\over n}L_{m: m+ \lfloor n\epsilon\rfloor K,m+\lfloor
n\epsilon\rfloor}=
\lim_{u} {\epsilon\over u} L_{Ku,u} =\gamma
_{\tt{R}}(K)\epsilon\qquad  \mbox{a.s.}
\end{equation}
For independent
sequences, (\ref{koond2}) follows from (\ref{LDV2ind}). Indeed, for
i.i.d. sequences, the random variables $L_{m:m+\lfloor n\epsilon
\rfloor K,m+\lfloor n\epsilon\rfloor} \mbox{ and }L_{\lfloor
n\epsilon\rfloor K,\lfloor n\epsilon\rfloor}$ are identically
distributed. By (\ref{LDV2ind}), thus, for any $\Delta>0$ (and
ignoring $\lfloor\cdot\rfloor$, for simplicity)
\begin{eqnarray*}
P \bigl(\bigl|L_{m:m+ n\epsilon K,m+n\epsilon }-\gamma(K)\epsilon n\bigr|>\Delta n \bigr)&=&P
\bigl(\bigl|L_{n\epsilon K, n\epsilon
}-\gamma(K)\epsilon n\bigr| >\Delta n \bigr)
\\
&=&P \biggl(\biggl|L_{n\epsilon K, n\epsilon} -{1\over K}\gamma(K)K\epsilon n\biggr| >
\Delta n \biggr)
\\
& =&P \biggl(\biggl|L_{n\epsilon K, n\epsilon} -\gamma\biggl({1\over K}\biggr)K
\epsilon n\biggr| >\Delta n \biggr)
\\
& =& P \biggl(L_{uk, k} -\gamma(u)k >\Delta{uk \over\epsilon}
\biggr)
\\
&\leq&2\exp\biggl[-{\Delta^2 u^2\over
2\epsilon^2(1+u)}{k}\biggr]=2\exp\biggl[-
{\Delta^2 \over\epsilon
(K+1)}n\biggr].
\end{eqnarray*}
Here $u={1\over K}$ and $k=\epsilon n K$.
In the third
equality, the relation $K\gamma({1\over K})=\gamma(K)$ is used. For
related sequences, the random variables $L_{m:m+\lfloor n\epsilon
\rfloor K,m+\lfloor n\epsilon\rfloor}$ and $L_{\lfloor n\epsilon
\rfloor K,\lfloor n\epsilon\rfloor}$ are not necessarily
identically distributed, hence another argument should be used. Let
$\overline a(m):=a^x(m)\vee a^y(m)$ and $\underline
a(m):=a^x({m})\wedge a^y(m)$. Let $Y_{m+k^y}$ (resp., $X_{m+k^x}$) be
the smallest element in $Y$ (resp., in $X$) that has ancestor at
least $\overline a(m)$. Similarly, let $X_{m-l^x}$ (resp.,
$Y_{m-l^y}$) be the smallest element in $X$ (resp., in $Y$) that has
ancestor at least $\underline a(m)$. If $a^x(m)\geq a^y(m)$, then
$k^x=l^y=0$ and if $a^x(m)\leq a^y(m)$, then $k^y=l^x=0$. Hence,
%
\begin{eqnarray}
\label{app1}
&&L(X_{m-l^x+1},\ldots, X_{m+ K\lfloor n\epsilon
\rfloor};Y_{m-l^y+1},
\ldots,Y_{m+ \lfloor n\epsilon\rfloor})\nonumber\\[-8pt]\\[-8pt]
&&\quad \geq L_{m:
m+ \lfloor n\epsilon\rfloor K,m+\lfloor n\epsilon\rfloor}\nonumber
\\
\label{app2}&&\quad \geq L(X_{m+k^x+1},\ldots, X_{m+ K\lfloor
n\epsilon\rfloor};Y_{m+k^y+1},
\ldots,Y_{m+ \lfloor
n\epsilon\rfloor}).
\end{eqnarray}
Note that the random variables
$X_{m+k^x+1},X_{m+k^x+2},\ldots$ and
$Y_{m+k^y+1},Y_{m+k^y+2},\ldots$ depend on i.i.d. random vectors
$U_{\overline{a}(m)+1},U_{\overline{a}(m)+2},\ldots\,$, where $U_i$ is
defined as in (\ref{vectors}). Similarly
$X_{m-l^x+1},X_{m-l^x+2},\ldots$ and
$Y_{m-l^y+1},Y_{m-l^y+2},\ldots$ depend on i.i.d. random vectors
$U_{\underline{a}(m)+1},\allowbreak U_{\underline{a}(m)+2},\ldots\,$. Hence, the
random variables
\[
L(X_{m+k^x+1},\ldots,X_{m+k^x+K\lfloor
n\epsilon\rfloor};Y_{m+k^y+1},
\ldots,Y_{m+k^y+\lfloor
n\epsilon\rfloor})
\]
and
\[
L(X_{1},\ldots,X_{K\lfloor
n\epsilon\rfloor};Y_{1},
\ldots,Y_{\lfloor n\epsilon\rfloor})
\]
have
the same distribution so that (as in the independent case) by
(\ref{LDV2})
\begin{eqnarray*}
&&P \bigl(\bigl|L(X_{m+k^x+1},\ldots,X_{m+ k^x+K\lfloor
n\epsilon\rfloor};Y_{m+k^y+1},
\ldots,Y_{m+k^y+\lfloor
n\epsilon\rfloor})-\epsilon\gamma_{\tt{R}}(K)n\bigr|>\Delta n \bigr)
\\
&&\quad =P \bigl(\bigl|L(X_{1},\ldots,X_{K\lfloor
n\epsilon\rfloor};Y_{1},
\ldots,Y_{\lfloor
n\epsilon\rfloor})-\epsilon\gamma_{\tt{R}}(K)n\bigr|>\Delta n \bigr)\leq
4\exp\biggl[-{1\over32}{\Delta^2 \over\epsilon K^2}n\biggr].
\end{eqnarray*}
Thus,
as $n$ grows,
\[
{1\over n}L(X_{m+k^x+1},\ldots,X_{m+k^x+K\lfloor
n\epsilon\rfloor};Y_{m+k^y+1},
\ldots,Y_{m+k^y+\lfloor
n\epsilon\rfloor})\to\epsilon\gamma_{\tt{R}}(K) \quad \mbox{a.s.}
\]
The random variables $k:=k^x\vee k^y$ and $l:=l^x\vee l^y$ satisfy
\[
k\vee l\leq\overline a(m) - \underline a(m).
\]
Note that
\[
\sum_{i=1}^{a^x(m)}D_i^x=
\sum_{i=1}^{a^y(m)}D_i^y=m.
\]
Now, using Hoeffding inequality for i.i.d. random variables $D_i^x$
and $D_i^y$, it is easy to see that
\[
{a^x(m)\over m}\to{1\over
p}\qquad  \mbox{a.s.},\qquad
{a^y(m)\over m}\to{1\over p}\qquad  \mbox{a.s.}
\]
Therefore,
\[
{k\vee l\over m}\leq{\overline a(m) -
\underline a(m)\over m}\to0\qquad  \mbox{a.s.}
\]
so that ${k(n)\over
n}\to0$ a.s. and ${l(n)\over n}\to0$ a.s. Since
\begin{eqnarray*}
&& \bigl|L(X_{m+k^x+1},\ldots, X_{m+ K\lfloor n\epsilon
\rfloor};Y_{m+k^y+1},
\ldots,Y_{m+ \lfloor
n\epsilon\rfloor})
\\
&&\quad {}-L(X_{m+k^x+1},\ldots,X_{m+k^x +K\lfloor
n\epsilon\rfloor};Y_{m+k^y+1},
\ldots,Y_{m+k^y+\lfloor
n\epsilon\rfloor}) \bigr|\leq k(n),
\end{eqnarray*}
from ${k(n)\over n}\to0$
a.s., we get that
%
\begin{eqnarray}
\label{lim1}
&&\lim_n {1\over n}L(X_{m+k^x+1},
\ldots, X_{m+ K\lfloor n\epsilon
\rfloor};Y_{m+k^y+1},\ldots,Y_{m+ \lfloor
n\epsilon\rfloor})\nonumber\\[-8pt]\\[-8pt]
&&\quad =\lim
_{u} {1\over u}\epsilon L_{Ku,u} =
\gamma _{\tt{R}}(K)\epsilon \qquad \mbox{a.s.}\nonumber
\end{eqnarray}
By similar argument,
%
\begin{eqnarray}
\label{lim2} &&\lim_n {1\over n}L(X_{m-l^x+1},
\ldots, X_{m+ K\lfloor n\epsilon
\rfloor};Y_{m-l^y+1},\ldots,Y_{m+ \lfloor
n\epsilon\rfloor})\nonumber\\[-8pt]\\[-8pt]
&&\quad =\lim
_{u} {1\over u}\epsilon L_{Ku,u} =
\gamma _{\tt{R}}(K)\epsilon \qquad \mbox{a.s.}\nonumber
\end{eqnarray}
The inequalities (\ref{lim1}) and (\ref{lim2}) together with
(\ref{app1}) and(\ref{app2}) imply (\ref{koond2}). The convergences
(\ref{koond1}) and (\ref{koond2}) imply
%
\begin{equation}
\label{lim+} \lim_{n}{1\over n}L_{n(1+c),n}=
\gamma_{\tt{R}}(1+c)>\gamma _{\tt{R}}\qquad  \mbox{a.s.}
\end{equation}
The limit in (\ref{lim+}) exists by Proposition~\ref{prop:piirv},
the inequality $\gamma_{\tt{R}}(1+c)>\gamma_{\tt{R}}$ follows from
(\ref{koond1}) and (\ref{koond2}), since $\epsilon\gamma_{\tt
{R}}(K)+(1-\epsilon)\gamma_{\tt{R}}>\gamma_{\tt{R}}$. This
proves that $\gamma_{\tt{R}}(a)>\gamma_{\tt{R}}$, when
$a>1$.

Finally,
\begin{eqnarray*}
 {1\over2n}L_{2n,2n}\geq {1\over2n}L_{n(1+c),n}+
{1\over
2n}L(X_{\lfloor n(1+c)\rfloor+1},\ldots,X_{2n};
Y_{n+1},\ldots, Y_{2n}).
\end{eqnarray*}
Since ${1\over2n}L_{2n,2n}\to\gamma_{\tt{R}}$,
a.s. and, using the same argument as proving (\ref{koond2}), we get
\[
\lim_{n}{1\over2n}L(X_{\lfloor n(1+c)\rfloor+1},
\ldots,X_{2n}; Y_{n+1},\ldots, Y_{2n})=
{1\over2}\lim_n {1\over n}
L_{(1-c)n,n}={\gamma_{\tt{R}}(1-c)\over2}\qquad  \mbox{a.s.},
\]
by
(\ref{lim+}), we have $\gamma_{\tt{R}}\geq{\gamma_{\tt
{R}}(1+c)\over
2}+{\gamma_{\tt{R}}(1-c)\over2}>{\gamma_{\tt{R}}+\gamma
_{\tt{R}}(1-c)\over2}.$ This implies that $\gamma_{\tt
{R}}(a)=\break \gamma_{\tt{R}}(1-c)<\gamma_{\tt{R}}$.
\end{pf}
%

The following corollary generalizes Proposition~\ref{STind} for
related sequences. Moreover, we allow the sequences to be unequal
length. Hence, we consider the case $X=X_1,\ldots,X_n$,
$Y=Y_1,\ldots, Y_m$, $n\leq m \leq n(1+\Delta)$, where $\Delta\geq
0$. The case $\Delta=0$ corresponds to the case $m=n$. Recall
the random variables $S:=j^h_1-1$ and $T:=n-i^h_k$, that obviously are
the functions of
$X$ and $Y$. The proof of the following corollary is very similar to
that one of Proposition~\ref{STind}.
%
\begin{corollary}\label{cor:st} Let $1>c>\Delta$. Then there exists
constant $d(c)>0$,
so that, for $n$ big enough, $P(T>cn)\leq\exp[-dn]$, $P(S>cn)\leq
\exp[-dn]$.
\end{corollary}
\begin{pf} As in the proof of Proposition~\ref{STind}, note that
for any $\bar{\gamma}$,
\[
\{S>cn\}\subset\{L_{n,m-cn}=L_{n,m}\}\subset\{L_{n,m-cn}
\geq\bar {\gamma}n\}\cup \{L_{n,m}\leq\bar{\gamma}n\}.
\]
By Lemma~\ref{lemma:appendix},
$\gamma_{\tt{R}}>\gamma_{\tt{R}}(1+\Delta-c)$. Let
$\bar{\gamma}:={1\over2}(\gamma_{\tt{R}}+\gamma_{\tt
{R}}(1+\Delta-c))$. Let $\epsilon:=\gamma_{\tt{R}}-\bar{\gamma
}$. Since
$L_{n,m-cn}\leq L_{n,(1+\Delta-c)n}$ and $L_n=L_{n,n}\leq L_{n,m}$,
Corollary~\ref{cor} states that for\vadjust{\goodbreak} $n$ big enough,
\begin{eqnarray*}
P(S>cn)&\leq& P (L_{n,(1+\Delta-c)n}\geq\bar{\gamma}n )+P (L_{n}\leq
\bar{\gamma}n )
\\
&=&P \bigl(L_{n,(1+\Delta-c)n}\geq\bigl(\gamma_{\tt{R}}(1+\Delta -c)+
\epsilon\bigr)n \bigr)+P \bigl(L_{n}\leq({\gamma}_{\tt{R}}-
\epsilon )n \bigr)\\
&\leq&8\exp\biggl[-{p\over32}(1+\Delta-c)
\epsilon^2 n\biggr].
\end{eqnarray*}
This concludes the proof.
\end{pf}
\end{appendix}


\section*{Acknowledgements}
Supported by the Estonian Science Foundation
Grant nr. 9288; SFB 701 of Bielefeld University and targeted
financing project SF0180015s12.



\printhistory


\begin{thebibliography}{23}

\bibitem{Alexander}
\begin{barticle}[mr]
\bauthor{\bsnm{Alexander},~\bfnm{Kenneth~S.}\binits{K.S.}}
(\byear{1994}).
\btitle{The rate of convergence of the mean length of the longest common
  subsequence}.
\bjournal{Ann. Appl. Probab.}
\bvolume{4}
\bpages{1074--1082}.
\bid{issn={1050-5164}, mr={1304773}}
\bptok{imsref}%
\end{barticle}
\endbibitem

\bibitem{macrounilargefluct}
\begin{barticle}[mr]
\bauthor{\bsnm{Amsalu},~\bfnm{Saba}\binits{S.}},
  \bauthor{\bsnm{Matzinger},~\bfnm{Heinrich}\binits{H.}} \AND
  \bauthor{\bsnm{Popov},~\bfnm{Serguei}\binits{S.}}
(\byear{2007}).
\btitle{Macroscopic non-uniqueness and transversal fluctuation in optimal
  random sequence alignment}.
\bjournal{ESAIM Probab. Stat.}
\bvolume{11}
\bpages{281--300}.
\bid{doi={10.1051/ps:2007014}, issn={1292-8100}, mr={2320822}}
\bptok{imsref}%
\end{barticle}
\endbibitem

\bibitem{Watermanphase}
\begin{barticle}[mr]
\bauthor{\bsnm{Arratia},~\bfnm{Richard}\binits{R.}} \AND
  \bauthor{\bsnm{Waterman},~\bfnm{Michael~S.}\binits{M.S.}}
(\byear{1994}).
\btitle{A phase transition for the score in matching random sequences allowing
  deletions}.
\bjournal{Ann. Appl. Probab.}
\bvolume{4}
\bpages{200--225}.
\bid{issn={1050-5164}, mr={1258181}}
\bptok{imsref}%
\end{barticle}
\endbibitem

\bibitem{Baeza1999}
\begin{barticle}[mr]
\bauthor{\bsnm{Baeza-Yates},~\bfnm{R.~A.}\binits{R.A.}},
  \bauthor{\bsnm{Gavald{\`a}},~\bfnm{R.}\binits{R.}},
  \bauthor{\bsnm{Navarro},~\bfnm{G.}\binits{G.}} \AND
  \bauthor{\bsnm{Scheihing},~\bfnm{R.}\binits{R.}}
(\byear{1999}).
\btitle{Bounding the expected length of longest common subsequences and
  forests}.
\bjournal{Theory Comput. Syst.}
\bvolume{32}
\bpages{435--452}.
\bid{doi={10.1007/s002240000125}, issn={1432-4350}, mr={1693383}}
\bptok{imsref}%
\end{barticle}
\endbibitem

\bibitem{barder}
\begin{barticle}[mr]
\bauthor{\bsnm{Barder},~\bfnm{Stanislaw}\binits{S.}},
  \bauthor{\bsnm{Lember},~\bfnm{J{\"u}ri}\binits{J.}},
  \bauthor{\bsnm{Matzinger},~\bfnm{Heinrich}\binits{H.}} \AND
  \bauthor{\bsnm{Toots},~\bfnm{M{\"a}rt}\binits{M.}}
(\byear{2012}).
\btitle{On suboptimal {LCS}-alignments for independent {B}ernoulli sequences
  with asymmetric distributions}.
\bjournal{Methodol. Comput. Appl. Probab.}
\bvolume{14}
\bpages{357--382}.
\bid{doi={10.1007/s11009-010-9206-7}, issn={1387-5841}, mr={2912343}}
\bptok{imsref}%
\end{barticle}
\endbibitem

\bibitem{BonettoLCS}
\begin{barticle}[mr]
\bauthor{\bsnm{Bonetto},~\bfnm{Federico}\binits{F.}} \AND
  \bauthor{\bsnm{Matzinger},~\bfnm{Heinrich}\binits{H.}}
(\byear{2006}).
\btitle{Fluctuations of the longest common subsequence in the asymmetric case
  of 2- and 3-letter alphabets}.
\bjournal{ALEA Lat. Am. J. Probab. Math. Stat.}
\bvolume{2}
\bpages{195--216}.
\bid{issn={1980-0436}, mr={2262762}}
\bptok{imsref}%
\end{barticle}
\endbibitem

\bibitem{china}
\begin{bbook}[mr]
\bauthor{\bsnm{Chao},~\bfnm{Kun-Mao}\binits{K.M.}} \AND
  \bauthor{\bsnm{Zhang},~\bfnm{Louxin}\binits{L.}}
(\byear{2009}).
\btitle{Sequence Comparison:
Theory and Methods}.
\bseries{Computational Biology}.
\blocation{London}: \bpublisher{Springer}.
\bid{mr={2723076}}
\bptok{imsref}%
\end{bbook}
\endbibitem

\bibitem{Sankoff1}
\begin{barticle}[mr]
\bauthor{\bsnm{Chvatal},~\bfnm{V{\'a}cl{\'a}v}\binits{V.}} \AND
  \bauthor{\bsnm{Sankoff},~\bfnm{David}\binits{D.}}
(\byear{1975}).
\btitle{Longest common subsequences of two random sequences}.
\bjournal{J. Appl. Probability}
\bvolume{12}
\bpages{306--315}.
\bid{issn={0021-9002}, mr={0405531}}
\bptok{imsref}%
\end{barticle}
\endbibitem

\bibitem{Durbin}
\begin{bbook}[auto:STB|2013/12/09|07:59:19]
\bauthor{\bsnm{Durbin},~\bfnm{R.}\binits{R.}},
  \bauthor{\bsnm{Eddy},~\bfnm{S.}\binits{S.}},
  \bauthor{\bsnm{Krogh},~\bfnm{A.}\binits{A.}} \AND
  \bauthor{\bsnm{Mitchison},~\bfnm{G.}\binits{G.}}
(\byear{1998}).
\btitle{Biological Sequence Analysis: Probabilistic Models of Proteins and
  Nucleic Acids}.
\blocation{Cambridge}: \bpublisher{Cambridge Univ. Press}.
\bptok{imsref}%
\end{bbook}
\endbibitem

\bibitem{hansen}
\begin{barticle}[mr]
\bauthor{\bsnm{Hansen},~\bfnm{Niels~Richard}\binits{N.R.}}
(\byear{2006}).
\btitle{Local alignment of {M}arkov chains}.
\bjournal{Ann. Appl. Probab.}
\bvolume{16}
\bpages{1262--1296}.
\bid{doi={10.1214/105051606000000321}, issn={1050-5164}, mr={2260063}}
\bptnote{check year}%
\bptok{imsref}%
\end{barticle}
\endbibitem

\bibitem{hirmo}
\begin{bmisc}[auto:STB|2013/12/09|07:59:19]
\bauthor{\bsnm{Hirmo},~\bfnm{Erik}\binits{E.}},
  \bauthor{\bsnm{Lember},~\bfnm{J{\"u}ri}\binits{J.}} \AND
  \bauthor{\bsnm{Matzinger},~\bfnm{Heinrich}\binits{H.}}
  (\byear{2012}).
\bhowpublished{Detecting the homology of DNA-sequence based on the variety of
  optimal alignments: A case study. Available at}
  \href{http://arxiv.org/abs/1210.3771}{arXiv:1210.3771}.
\bptok{imsref}%
\end{bmisc}
\endbibitem

\bibitem{LCIS}
\begin{barticle}[mr]
\bauthor{\bsnm{Houdr{\'e}},~\bfnm{Christian}\binits{C.}},
  \bauthor{\bsnm{Lember},~\bfnm{J{\"u}ri}\binits{J.}} \AND
  \bauthor{\bsnm{Matzinger},~\bfnm{Heinrich}\binits{H.}}
(\byear{2006}).
\btitle{On the longest common increasing binary subsequence}.
\bjournal{C. R. Math. Acad. Sci. Paris}
\bvolume{343}
\bpages{589--594}.
\bid{doi={10.1016/j.crma.2006.10.004}, issn={1631-073X}, mr={2269870}}
\bptok{imsref}%
\end{barticle}
\endbibitem

\bibitem{kiwi}
\begin{barticle}[mr]
\bauthor{\bsnm{Kiwi},~\bfnm{Marcos}\binits{M.}},
  \bauthor{\bsnm{Loebl},~\bfnm{Martin}\binits{M.}} \AND
  \bauthor{\bsnm{Matou{\v{s}}ek},~\bfnm{Ji{\v{r}}{\'{\i}}}\binits{J.}}
(\byear{2005}).
\btitle{Expected length of the longest common subsequence for large alphabets}.
\bjournal{Adv. Math.}
\bvolume{197}
\bpages{480--498}.
\bid{doi={10.1016/j.aim.2004.10.012}, issn={0001-8708}, mr={2173842}}
\bptok{imsref}%
\end{barticle}
\endbibitem

\bibitem{notsym}
\begin{barticle}[mr]
\bauthor{\bsnm{Lember},~\bfnm{J{\"u}ri}\binits{J.}} \AND
  \bauthor{\bsnm{Matzinger},~\bfnm{Heinrich}\binits{H.}}
(\byear{2009}).
\btitle{Standard deviation of the longest common subsequence}.
\bjournal{Ann. Probab.}
\bvolume{37}
\bpages{1192--1235}.
\bid{doi={10.1214/08-AOP436}, issn={0091-1798}, mr={2537552}}
\bptok{imsref}%
\end{barticle}
\endbibitem

\bibitem{rate}
\begin{barticle}[mr]
\bauthor{\bsnm{Lember},~\bfnm{J{\"u}ri}\binits{J.}},
  \bauthor{\bsnm{Matzinger},~\bfnm{Heinrich}\binits{H.}} \AND
  \bauthor{\bsnm{Torres},~\bfnm{Felipe}\binits{F.}}
(\byear{2012}).
\btitle{The rate of the convergence of the mean score in random sequence
  comparison}.
\bjournal{Ann. Appl. Probab.}
\bvolume{22}
\bpages{1046--1058}.
\bid{doi={10.1214/11-AAP778}, issn={1050-5164}, mr={2977985}}
\bptok{imsref}%
\end{barticle}
\endbibitem

\bibitem{VollmerReport}
\begin{bmisc}[auto:STB|2013/12/09|07:59:19]
\bauthor{\bsnm{Lember},~\bfnm{J{\"u}ri}\binits{J.}},
  \bauthor{\bsnm{Matzinger},~\bfnm{Heinrich}\binits{H.}} \AND
  \bauthor{\bsnm{Vollmer},~\bfnm{Anna}\binits{A.}}
(\byear{2007}).
\bhowpublished{Path properties of LCS-optimal alignments. SFB 701
  Preprintreihe, Univ. Bielefeld (07 - 77)}.
\bptok{imsref}%
\end{bmisc}
\endbibitem

\bibitem{periodicLCS}
\begin{barticle}[mr]
\bauthor{\bsnm{Matzinger},~\bfnm{Heinrich}\binits{H.}},
  \bauthor{\bsnm{Lember},~\bfnm{J{\"u}ri}\binits{J.}} \AND
  \bauthor{\bsnm{Durringer},~\bfnm{Clement}\binits{C.}}
(\byear{2007}).
\btitle{Deviation from mean in sequence comparison with a periodic sequence}.
\bjournal{ALEA Lat. Am. J. Probab. Math. Stat.}
\bvolume{3}
\bpages{1--29}.
\bid{issn={1980-0436}, mr={2324746}}
\bptok{imsref}%
\end{barticle}
\endbibitem

\bibitem{siegmund}
\begin{barticle}[mr]
\bauthor{\bsnm{Siegmund},~\bfnm{David}\binits{D.}} \AND
  \bauthor{\bsnm{Yakir},~\bfnm{Benjamin}\binits{B.}}
(\byear{2000}).
\btitle{Approximate {$p$}-values for local sequence alignments}.
\bjournal{Ann. Statist.}
\bvolume{28}
\bpages{657--680}.
\bid{doi={10.1214/aos/1015951993}, issn={0090-5364}, mr={1792782}}
\bptok{imsref}%
\end{barticle}
\endbibitem

\bibitem{steele86}
\begin{barticle}[mr]
\bauthor{\bsnm{Steele},~\bfnm{J.~Michael}\binits{J.M.}}
(\byear{1986}).
\btitle{An {E}fron-{S}tein inequality for nonsymmetric statistics}.
\bjournal{Ann. Statist.}
\bvolume{14}
\bpages{753--758}.
\bid{doi={10.1214/aos/1176349952}, issn={0090-5364}, mr={0840528}}
\bptok{imsref}%
\end{barticle}
\endbibitem

\bibitem{Waterman-estimation}
\begin{barticle}[auto:STB|2013/12/09|07:59:19]
\bauthor{\bsnm{Waterman},~\bfnm{Michael~S.}\binits{M.S.}}
(\byear{1994}).
\btitle{Estimating statistical significance of sequence alignments}.
\bjournal{Phil. Trans. R. Soc. Lond. B}
\bvolume{344}
\bpages{383--390}.
\bptok{imsref}%
\end{barticle}
\endbibitem

\bibitem{watermanintrocompbio}
\begin{bbook}[auto:STB|2013/12/09|07:59:19]
\bauthor{\bsnm{Waterman},~\bfnm{Michael~S.}\binits{M.S.}}
(\byear{1995}).
\btitle{Introduction to Computational Biology}.
\blocation{London}: \bpublisher{Chapman \& Hall}.
\bptok{imsref}%
\end{bbook}
\endbibitem

\bibitem{Vingron}
\begin{barticle}[mr]
\bauthor{\bsnm{Waterman},~\bfnm{Michael~S.}\binits{M.S.}} \AND
  \bauthor{\bsnm{Vingron},~\bfnm{Martin}\binits{M.}}
(\byear{1994}).
\btitle{Sequence comparison significance and {P}oisson approximation}.
\bjournal{Statist. Sci.}
\bvolume{9}
\bpages{367--381}.
\bid{issn={0883-4237}, mr={1325433}}
\bptok{imsref}%
\end{barticle}
\endbibitem

\end{thebibliography}
\end{document}